\documentclass[11pt,twoside]{article}
\usepackage{latexsym}
\usepackage{amssymb,amsbsy,amsmath,amsfonts,amssymb,amscd}
\usepackage{mathrsfs}
\usepackage{epsfig, graphicx}
\usepackage{graphics,graphicx,epsfig,subfigure}
\usepackage{color,verbatim}
\usepackage{authblk}
\usepackage[font=small,format=plain,labelfont=bf,up,textfont=up]{caption} 
\usepackage{subfigure}
\newcommand{\commentout}[1]{}
\newcommand{\R}{\mathbb{R}}
\newcommand{\N}{\mathbb{N}}

\newcommand {\eps}  {\varepsilon}

\newcommand {\vp} {\varphi}

\newcommand {\Chi} {{\bf \raise 2pt \hbox{$\chi$}} }

\newcommand {\ext} { {\rm ext} }

\newcommand {\cak} { {\mathcal K} }

\newcommand {\p}   {\partial}
\newcommand{\ds}{\displaystyle}
\newcommand{\ud}{\, \mathrm{d}}
\def\ddef{\stackrel{\rm def}{=}}

\newcommand{\beq}{\begin{equation}}
\newcommand{\beqa}{\begin{eqnarray}}
\newcommand{\bea} {\begin{array}{ll}}
\newcommand{\beqan}{\begin{eqnarray*}}
\newcommand{\eeq}{\end{equation}}
\newcommand{\eeqa}{\end{eqnarray}}
\newcommand{\eeqan}{\end{eqnarray*}}
\newcommand{\eea} {\end{array}}
\newtheorem{theorem}{Theorem}[section]
\newtheorem{lem}[theorem]{Lemma}

\newtheorem{remark}[theorem]{Remark}
\newtheorem{prop}[theorem]{Proposition}
\newtheorem{coro}[theorem]{Corollary}
\newcommand{\qed}{{ \hfill
                       {\unskip\kern 6pt\penalty 500
                       \raise -2pt\hbox{\vrule\vbox to 6pt{\hrule width 6pt
                       \vfill\hrule}\vrule} \par}  \medskip }}

\title{\Large \bf
 From  Vlasov--Poisson  and Vlasov--Poisson--Fokker--Planck Systems to 
 Incompressible Euler Equations: the case with finite charge}
 \author[1]{Julien Barr\'e\thanks{ {\tt julien.barre@unice.fr}}}
  \author[1]{David Chiron\thanks{ {\tt chiron@unice.fr}}}
\author[1,2]{Thierry Goudon\thanks{ {\tt thierry.goudon@inria.fr}}}
 \author[3]{Nader Masmoudi\thanks{ {\tt masmoudi@cims.nyu.edu}}}
\affil[1]{Univ. Nice Sophia Antipolis, CNRS, Labo. J.-A. Dieudonn\'e, UMR 7351\\ 

Parc Valrose, F-06108 Nice, France}
\affil[2]{Inria,  Sophia Antipolis M\'editerran\'ee Research Centre, Project COFFEE}
\affil[3]{Courant Institute for Math. Sciences, New York University\\

251 Mercer St.,
      New York, NY 10012, U.S.A.}

\begin{document}
\maketitle

\abstract{We study the asymptotic regime of strong electric fields that leads from the Vlasov--Poisson 
system to the Incompressible Euler equations.
We also deal with the Vlasov--Poisson--Fokker--Planck system which induces dissipative effects. 
The originality consists in considering a situation with a finite total charge confined by a 
strong external field. In turn, the limiting equation is set in a bounded domain, the shape 
of which is determined by the external confining potential. 
The analysis extends to the situation where the limiting density is non--homogeneous and 
where the Euler equation is replaced by the Lake Equation, also called Anelastic Equation.}

\medskip
\noindent
  \small{\bf{Key words. }}\small{Plasma physics. Vlasov--Poisson system. 
  Vlasov--Poisson--Fokker--Planck system. Incompressible Euler equations. 
  Lake equations.
  Quasi--neutral regime. Modulated energy. Relative entropy. 
  }

   \medskip
   \noindent
 \small{\bf{2010 MSC Subject Classification.} }\small{
   82D10, 
   35Q35, 
82C40. 
}


\section{Introduction}

\subsection{The Vlasov-Poisson equation in a confining potential} 

\begin{figure}
\begin{center}
\includegraphics[width=0.8\linewidth,height=0.35\linewidth]{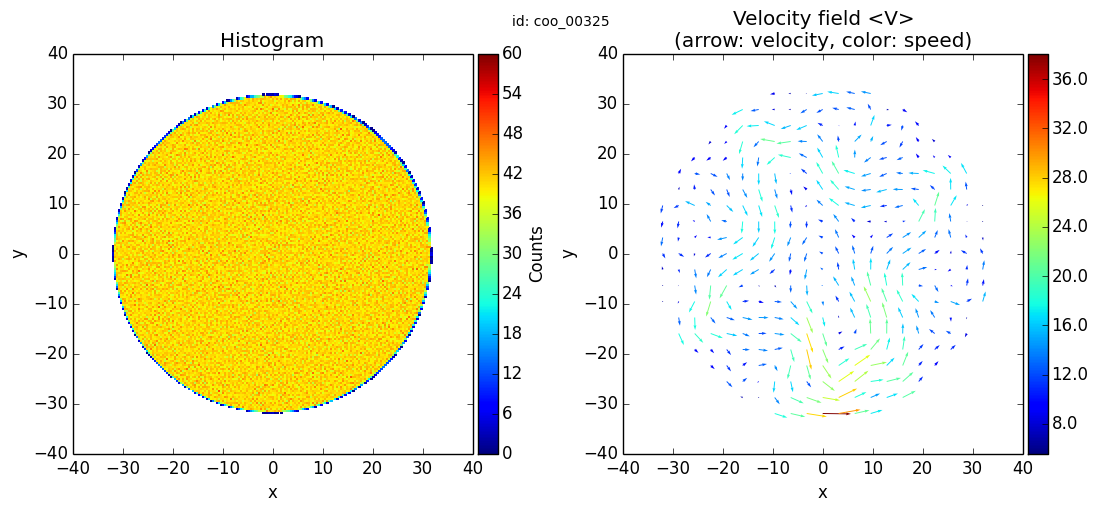}
\end{center}
\caption{Snapshot of a 2D simulation of confined charged
  particles. Particles are subjected to the combination of a harmonic
  and isotropic external potential, a strong Coulomb repulsion, a
  friction and a noise. An external force has been added from the left
  to the right in the lower half of the cloud, in order to set the
  particles in motion. Left: instantaneous locally averaged density
  field. The density is almost uniform inside a ball, and almost zero
  outside. Right: instantaneous locally averaged velocity field. (By
  courtesy of A. Olivetti \cite{Oliv}.)}
\label{simu}
\end{figure}

We are interested in the behavior as $\eps$ tends to 0 of the solutions of the following Vlasov equation
\begin{equation}
\tag{V}
\label{VP1}
\partial_t f_\eps+v\cdot\nabla_x f_\eps - \left(\ds\frac1\eps \nabla_x \Phi_\ext 
+ \nabla_x \Phi_\eps \right) \cdot\nabla_v f_\eps=0,
\end{equation}
where the potential $ \Phi_\eps $ is defined self-consistently by the Poisson equation
\begin{equation}
\tag{P}
\label{VP2}
\Delta_x\Phi_\eps= - \ds\frac1\eps \rho_\eps ,\qquad \rho_\eps(t,x)=\ds\int f_\eps(t,x,v)\ud v ,
\end{equation}
and where $ \eps^{-1} \Phi_\ext $ is a strong external potential applied to the system. 
The problem holds in the entire space: $x\in \mathbb R^N$, $v\in\mathbb R^N$ and it is 
completed by an initial data with \emph{finite} charge
\begin{equation}
\label{alacharge}
f_\varepsilon\Big|_{t=0}=f_\eps^{\mathrm{init}},\qquad 
\ds\iint f_\eps^{\mathrm{init}}\ud v\ud x=\mathfrak m\in (0,\infty).
\end{equation}
Notice that $\Phi_\eps$ is of size $ \eps^{-1} $ and we shall consider
the applied potential $\frac1\eps \Phi_\ext $ also of size $ \eps^{-1}
$.  The problem  is motivated by the study of \emph{non neutral
  plasmas} (see \cite{Dubinreview} for a review): these are
collections of particles all with the same sign of charge, for
instance pure electron, or pure ion plasmas. There are several methods
to confine such a plasma, among which the Paul trap, which uses an
oscillating electric field. The Penning trap, which uses a combination
of static electric and magnetic fields, is also standard, but
\eqref{VP1} is not directly relevant to this situation since there is
no magnetic field in it. A non neutral plasma picture has also been used to
describe trapped neutral atoms~\cite{Mendonca2008}, in the regime
where multiple diffusion of quasi resonant photons induces an
effective interaction force between atoms which is formally similar to
a Coulomb force~\cite{Walker1990}. In this case, the system is however
dissipative; a standard way to take this effect into account is to add
to \eqref{VP1} a Fokker-Planck operator acting on
velocities~\cite{Dalibard1985}. We will also discuss this situation.
  In these physical examples, the small $\eps$ limit is indeed relevant in
many experimental situations. 
Figure~\ref{simu} 
corresponds to a numerical simulation of such an experiment.
It strongly suggests the existence of a limiting fluid
model where the density is  nothing but the characteristic function of  a ball.
Our goal is to justify that,  indeed, a simpler model, purely of hydrodynamic type,
can be used to describe the particles in this asymptotic limit.

In fact, we shall see that the
  limiting model holds in a domain the shape of which depends on the
  external potential $ \Phi_\ext $.  But, to start with,
  we can consider a quadratic and isotropic potential, say:
  \begin{equation}\label{iso}
  \Phi_\ext(x)=\ds\frac{1}{2N}|x|^2\end{equation}
  where we remind the reader that $N$ stands for the space dimension.
  It corresponds to the case displayed in Figure~\ref{simu}.
 The confining potential $
\eps^{-1}\Phi_\ext $ tends to strongly localize  in space the particle density. 
On the support of  the limiting density $ \rho $, the
electric force $ \eps^{-1} \nabla_x \Phi_\ext + \nabla_x \Phi_\eps $
should be of order one.  By \eqref{VP2}, this  imposes   that 
 $\Delta \Phi_\ext + \eps \Delta
\Phi_\eps= \eps\nabla_x\cdot ( \eps^{-1} \nabla_x \Phi_\ext + \nabla_x \Phi_\eps )= 
\Delta \Phi_\ext -\rho_\eps =1-\rho_\eps $ is 
   of order  $\mathscr O(\eps)$ on the support of $\rho$
 for the potential
\eqref{iso}.
Clearly, due to the condition of  finite charge \eqref{alacharge}, the limiting 
density cannot be constant uniformly on the whole space. 
The intuition is that the limiting density has the same radial
symmetries as both the 
external potential \eqref{iso} and the Poisson kernel, see \eqref{Gamma} below.
 Actually, we shall prove some convergence
of $\rho_\eps $ to
\begin{equation}
\label{Guess}
n_{\mathrm{e}} (x) = \mathbf 1_{B(0,R)}(x),
\end{equation}
where $\mathbf{1}_U$ denotes the characteristic function of the set
$U$.  The radius $R$ depends on the total mass  $\mathfrak m$ so that the charge
constraint \eqref{alacharge} is fulfilled.  
In order to find a hydrodynamic description of the particles, it is 
convenient to associate to the particle distribution function $f_\eps$ the
following macroscopic quantities
\[
\begin{array}{ll}
  \text{Current:}\quad& J_\eps(t,x) \ddef \ds\int v\ f_\eps(t,x,v)\ud v,
  \\
  \text{Kinetic pressure:}\quad& \mathbb P_\eps(t,x) \ddef \ds\int v\otimes v\ f_\eps(t,x,v)\ud v.
\end{array}
\]
It turns out that the current looks like
\begin{equation}\label{Guess2}
  J_\eps(t,x)=\rho_\eps (t,x) V_\eps(t,x) \xrightarrow[\eps \rightarrow 0]{} 
 n_{\mathrm{e}} (x) V(t,x) =\mathbf{1}_{B(0,R)}(x)V(t,x),
\end{equation}
where $V$ solves the Incompressible Euler system in $B(0,R)$:
\begin{equation}
\tag{IE}\label{Euler}
\left\{ \begin{array}{l}
\partial_t V+\nabla_x\cdot (V\otimes V)+\nabla_x p=0, \\
\nabla_x\cdot V=0,
\end{array} \right.
\end{equation}
with an appropriate initial condition, and no flux boundary condition
on $\partial B(0,R)$.  In \eqref{Euler}, the pressure $p$ appears as
the Lagrange multiplier associated with the constraint that $V$ is
divergence free.  This incompressibility condition comes from charge
conservation: integrating \eqref{VP1} with respect to the velocity
variable $v$, we get
\begin{equation}\label{ch_cons}
  \partial_t\rho_\eps+\nabla_x\cdot J_\eps=0.
\end{equation}
Letting $\eps$ go to 0, with \eqref{Guess} and \eqref{Guess2}, we
deduce that $V$ is solenoidal. Obtaining the evolution equation for
$V$ is more intricate.

The analysis of such asymptotic problems goes back to \cite{Bre},
where a specific modulated energy method was introduced. It has been
revisited in \cite{Mas}, still by using a modulated energy method, but
which is able to account for oscillations present within the
system. Accordingly, more general initial data can be dealt with in
\cite{Mas}. However, these results hold either on the torus $\mathbb
T^N$, or in the whole space with data having infinite charge, that is
$ \iint f (x,v)\ud v\ud x = \infty $. A case with finite charge, but a
different Poisson equation which leads to a compressible
hydrodynamic limit, has been considered in \cite{HanKwan}, again
with a modulated energy.  Our goal in this article is twofold:
\begin{itemize}
\item To prove the convergence to \eqref{Euler} in the case of a
  trapped system, with finite charge.  Even though our proof also
  relies on a modulated energy functional, there are new difficulties:
  the shape of the domain on which the limiting equation (IE) holds is
  determined by the external potential $ \Phi_{\rm ext} $, and a
  careful treatment of the boundary is needed.
\item To prove the convergence to the analog of \eqref{Euler} in the case of 
a trapped dissipative system.
\end{itemize}
Both improvements are relevant for experiments on non neutral plasmas
or large magneto-optical traps.

\subsection{Statement of the results}

In what follows we shall deal with a smooth solution $(t,x)\mapsto
V(t,x)\in \mathbb R^N$ (possibly defined on a small enough time
interval $[0,T]$) of the incompressible Euler equation \eqref{Euler}
set on the ball $B(0,R)$, completed with no-flux boundary condition
\begin{equation}\label{CLEuler}
V(t,x)\cdot \nu(x)\Big|_{|x|=R}=0,
\end{equation}
where $\nu(x)$ denotes the outward unit vector at $x\in \partial
B(0,R)$ (namely $\nu(x)=x/|x|$).  We work with solutions $V$ that
belongs to 
$L^\infty(0,T;H^s(B(0,R)))$, for a certain $s>0$ large enough.

\begin{theorem}[\cite{RT}-\cite{Temam2}] 
  \label{Eulerexiste}Let $V^{\mathrm{init}}: B(0,R)\rightarrow \mathbb
  R^N$ be a divergence free vector field in $H^s$, with $s>1+N/2$,
  satisfying the no flux condition $ V^{\mathrm{init}} \cdot \nu = 0 $
  on $\p B (0,R) $. There exists 
  $ T > 0 $ and a unique solution $V\in L^\infty (0,T ;H^s(B(0,R)))$ of 
  \eqref{Euler} with the no flux condition \eqref{CLEuler}. Moreover, we have
\[
\ds \sup _{0\leq t\leq T}\Big( \|V(t)\|_{H^s} + \|\partial_t
V(t)\|_{H^{s-1}} + \|\nabla_xp(t)\|_{H^s} + \|\partial_t\nabla_x
p(t)\|_{H^{s-1}}\Big)\leq C(T ) \] for some positive constant $C(T )$
depending on $T$ and the initial datum.
\end{theorem}

\noindent
If $ N =1 $, the only divergence free vector field $ V^{\mathrm{init}}
$ satisfying \eqref{CLEuler} is $ V^{\mathrm{init}} \equiv 0 $ and
then the solution given in Theorem \ref{Eulerexiste} is $ V \equiv 0 $.\\

For further purposes, we need to consider an extension
$\mathscr V$ of the solution $V$ to \eqref{Euler} with
\eqref{CLEuler}, defined on the whole space and compactly supported.
Namely we require $ \mathscr V\in 
L^\infty(0,T;H^s(\mathbb R^N)) $ to satisfy
\begin{equation}
\label{CondExtension}
\mathscr V\Big|_{B(0,R)}= V , \qquad 
\mathscr V\Big|_{ \mathbb{R}^N \setminus B(0,2R)}= 0,\qquad
\mathscr V(t,x)\cdot \nu(x)\Big|_{|x|=R}=0 .
\end{equation}
For the construction of such an extension, we refer to \cite[Chapter I:
Theorem 2.1 p. 17 \& Theorem 8.1 p. 42]{LM}. For an
extension which is in addition divergence--free, see
Lemma~\ref{extension} in the appendix.

In order to state our first result, we need to introduce an auxiliary 
potential function $ \Phi_{\rm e} $. Suppose that \eqref{Guess} indeed
holds true.  Then, by using  \eqref{VP2} and  $ \Delta \Phi_\ext = 1 $
for the potential \eqref{iso}, we infer, for $\eps \rightarrow 0 $,
$$
 \Delta (\Phi_\ext + \eps \Phi_\eps ) = \Delta \Phi_\ext - \rho_\eps 
 \to \Delta \Phi_\ext - {\bf 1}_{B(0,R)} 
 = {\bf 1}_{\R^N \setminus B(0,R)} .
$$
Moreover, since we want the electric force 
$ \eps^{-1} \nabla_x \Phi_\ext + \nabla_x \Phi_\eps = \eps^{-1} ( \nabla_x \Phi_\ext + \eps \nabla_x \Phi_\eps ) $ 
to be of order one on the ball $
B(0,R) $, this imposes $ \Phi_\ext + \eps \Phi_\eps $ to be close to a
constant, say zero, on the ball $ B(0,R) $.  It is therefore natural
to look for a solution $ \Phi_{\rm e} $ to the Poisson problem
\begin{equation}
\label{defPhie0}
\Delta \Phi_{\mathrm e}(x)=1-n_{\mathrm{e}}(x)= \mathbf 1_{\R^N \setminus B(0,R)} ,
\quad \quad \quad \Phi_{\mathrm e} = 0 \ {\rm in} \ B(0,R).
\end{equation}
In this specific case, we can find an explicit radially symmetric solution:
\begin{equation}\label{defPhie}
  \Phi_{\mathrm{e}}
  (x)=\mathbf 1_{\R^N \setminus B(0,R)}  \times \left\{\begin{array}{ll}
        \ds\frac {|x|^2}{2N}+\ds\frac{R^N}{N(N-2)|x|^{N-2}}-\ds\frac{R^2}{2(N-2)}
        \quad & 
      \text{ if $N>2$},\\[.4cm]
        \ds\frac {|x|^2- R^2}{4}-\ds\frac{R^2}{2}\ln(|x|/R) 
        \quad &  \text{ if $N=2$,}\\[.4cm]
      \ds\frac12 ( |x| - R )^2\ &  \text{ if $N=1$.}
\end{array}
\right.\end{equation}
With $\Phi_{\mathrm{e}}$ and $n_{\mathrm{e}}$ in hand, we split the Poisson
equation \eqref{VP2} as follows, where $ n_{\rm e} $ is defined in \eqref{Guess},
\[
\Delta_x
\Phi_\eps(t,x)=\ds\frac{1-n_{\mathrm{e}}(x)}{\eps}+\ds\frac{n_{\mathrm{e}}(x)-\rho_\eps(t,x)}{\eps}-\frac{1}{\eps}\Delta \Phi_\ext=
\frac{1}{\eps}\Delta_x \Phi_{\mathrm
  e}(x)+\ds\frac{1}{\sqrt\eps}\Delta _x \Psi_\eps(t,x)
-\frac{1}{\eps}\Delta \Phi_\ext,
\]
namely, we have
\begin{equation} 
\label{VP2psi}
\Phi_\eps(x)+\frac{1}{\eps} \Phi_\ext=\ds\frac1\eps\Phi_{\mathrm e}(x) 
+\ds\frac{1}{\sqrt\eps}\Psi_\eps(t,x),\qquad
\Delta _x \Psi_\eps(t,x)=\ds\frac{1}{\sqrt\eps}(n_{\mathrm e}(x)-\rho_\eps(t,x)),
\end{equation}
where
$\Psi_\eps$ represents the fluctuations of the potential.
According to \cite{Bre}, we introduce a modulated energy:
\begin{equation*} 
\mathscr H_{\mathscr V,\eps} \ddef 
\ds\frac12\ds\iint |v-\mathscr V|^2\ f_\eps\ud v\ud x+
\ds\frac12\ds\int |\nabla_x\Psi_\eps|^2\ud x +\ds\frac1\eps\ds\iint
\Phi_{\mathrm e}\ f_\eps\ud v\ud x .
\end{equation*}
When the external potential is given by \eqref{iso}, we shall establish the following 
statement\footnote{Throughout the paper, we denote by $\mathscr M^1(X)$ the space of 
bounded measures on $X\subset \mathbb R^D$.
It identifies with the dual space of the separable space $C^0_0(X)$ of the
continuous functions that vanish at infinity.}.

\begin{theorem}
\label{thresultradial} 
Let $ V^{\mathrm{init}} \in H^s(B(0,R))$ satisfy $\nabla_x\cdot V^{\mathrm{init}}=0$ and the no flux
condition \eqref{CLEuler}. Denote by $V$ the solution, on 
$ [ 0, T ] $, to \eqref{Euler} with the no flux condition \eqref{CLEuler} given in Theorem
\ref{Eulerexiste}. Consider $\mathscr V$ a smooth extension of $V$ satisfying the conditions 
\eqref{CondExtension}.  Let $f_\eps^{\mathrm{init}}:\mathbb R^N\times
\mathbb R^N\rightarrow [0,\infty)$ be a sequence of integrable
functions that satisfy the following requirements
\begin{equation}
\label{query}\left\{ 
\begin{array}{l}
\ds\iint f_\eps^{\mathrm{init}}\ud v\ud x=\mathfrak{m},
\\
\ds\lim_{\eps\rightarrow 0}  \left\{
\ds\frac12\ds\iint
 \lvert v-\mathscr V^{\mathrm{init}} \rvert^2\ f_\eps^{\mathrm{init}}\ud v\ud x+ \ds\frac12\ds\int
|\nabla_x\Psi_\eps^{\mathrm{init}}|^2\ud x
+\ds\frac1\eps\ds\iint
\Phi_{\mathrm e}\ f_\eps^{\mathrm{init}}\ud v\ud x\right\}=0.
\end{array} \right.
\end{equation}
Then, the associated solution $f_\eps$ of the Vlasov--Poisson equation 
\eqref{VP1}--\eqref{VP2} satisfies, 
as $ \eps\rightarrow 0 $,
\begin{itemize}
\item[i)] $\rho_\eps$ converges to $n_{\mathrm e}$ in 
$C^0([0,T];\mathscr M^1(\mathbb R^N)-\text{weak}-\star)$;
\item[ii)] $ \mathscr H_{\mathscr V,\eps}$ converges to $0 $ uniformly
  on $[ 0, T ]$;
\item[iii)] $J_\eps$ converges to $J$ in $\mathscr M^1([0,T]\times
  \mathbb R^N)$ weakly-$\star $, the limit $J$ lies in
  $L^\infty(0,T;L^2(\mathbb R^N))$ and satisfies $J\big|_{[0,T]\times
    B(0,R)}=V$, $\nabla_x\cdot J=0$ and $J\cdot \nu(x)\big|_{\p B(0,R)}=0$.
\end{itemize}\end{theorem}

\begin{remark} 
(i) Here, we were not very precise about the type of solutions to 
the Vlasov--Poisson system \eqref{VP1}--\eqref{VP2} we are considering. 
We refer to  \cite{Pfa,PLLP} for the construction of  global 
regular solutions to the system and some extra conditions to 
ensure the propagation of regularity. There are also weaker notions 
of solutions (weak solutions or renormalized solutions) to which 
our theorem can apply.  We refer the reader to the introduction of \cite{PLLP} 
for a discussion about these solutions.

(ii)  The second part of the hypothesis \eqref{query} imposes that the
  initial modulated energy is small; this is a strong hypothesis on
  the initial data.  
  When the problem is set on the torus, or on the whole space with infinite charge, it can be relaxed, 
  see \cite{Mas}.
  In the present framework, going beyond \eqref{query}  would certainly require 
  a fine description of boundary layers on $\{|x|=R\}$.
  Assuming  \eqref{query}, point \emph{ii)} of the theorem then ensures that
  the modulated energy remains small at later times.  As typical
  initial data satisfying \eqref{query}, we can take
$$
f_\eps^{\mathrm{init}} (x , v) = \frac{n_{\mathrm e}( x) - \delta_\eps
  \Delta \chi( x)}{ \sigma_\eps^N } G \left( \frac{v - \mathscr
    V^{\mathrm{init}}(x)}{\sigma_\eps} \right) ,
$$
where $ \chi \in C^\infty_c( B(0,R) ) $ and where $ G $ is a nonnegative function that 
belongs to the Schwartz space  and satisfies $ \int G \ud v = 1 $ (for
instance, $G$ is a normalized Gaussian $ G(v) = ( 2\pi)^{-N/2} \exp( -
\lvert v \rvert^2 /2 ) $). Then, we choose $ \sigma_\eps
\rightarrow 0 $ as $\eps \rightarrow 0 $ and $ \delta_\eps = o (
\sqrt{\eps}) $ (so that $ \ds n_{\mathrm e} - \delta_\eps \Delta \chi
\ge 0 $ for $ \eps $ small enough).  Indeed, we easily obtain $ \iint
\Phi_{\mathrm e}\ f_\eps^{\mathrm{init}}\ud v\ud x = 0 $, $ \iint
|v-\mathscr V^{\mathrm{init}}|^2\ f_\eps^{\mathrm{init}}\ud v\ud x =
\sigma_\eps^2 ( \int |v|^2 G(v) \ud v ) (\int n_{\mathrm e} \ud x)
\rightarrow 0 $ and $ \Psi_\eps = \delta_\eps / \sqrt{\eps} \chi $,
hence $ \int |\nabla_x\Psi_\eps^{\mathrm{init}}|^2\ud x =
\delta_\eps^2 / \eps \int |\nabla_x\chi|^2\ud x \rightarrow 0 $.
\end{remark}

\begin{figure}
\begin{center}
\includegraphics[width=0.5\linewidth,height=0.5\linewidth]{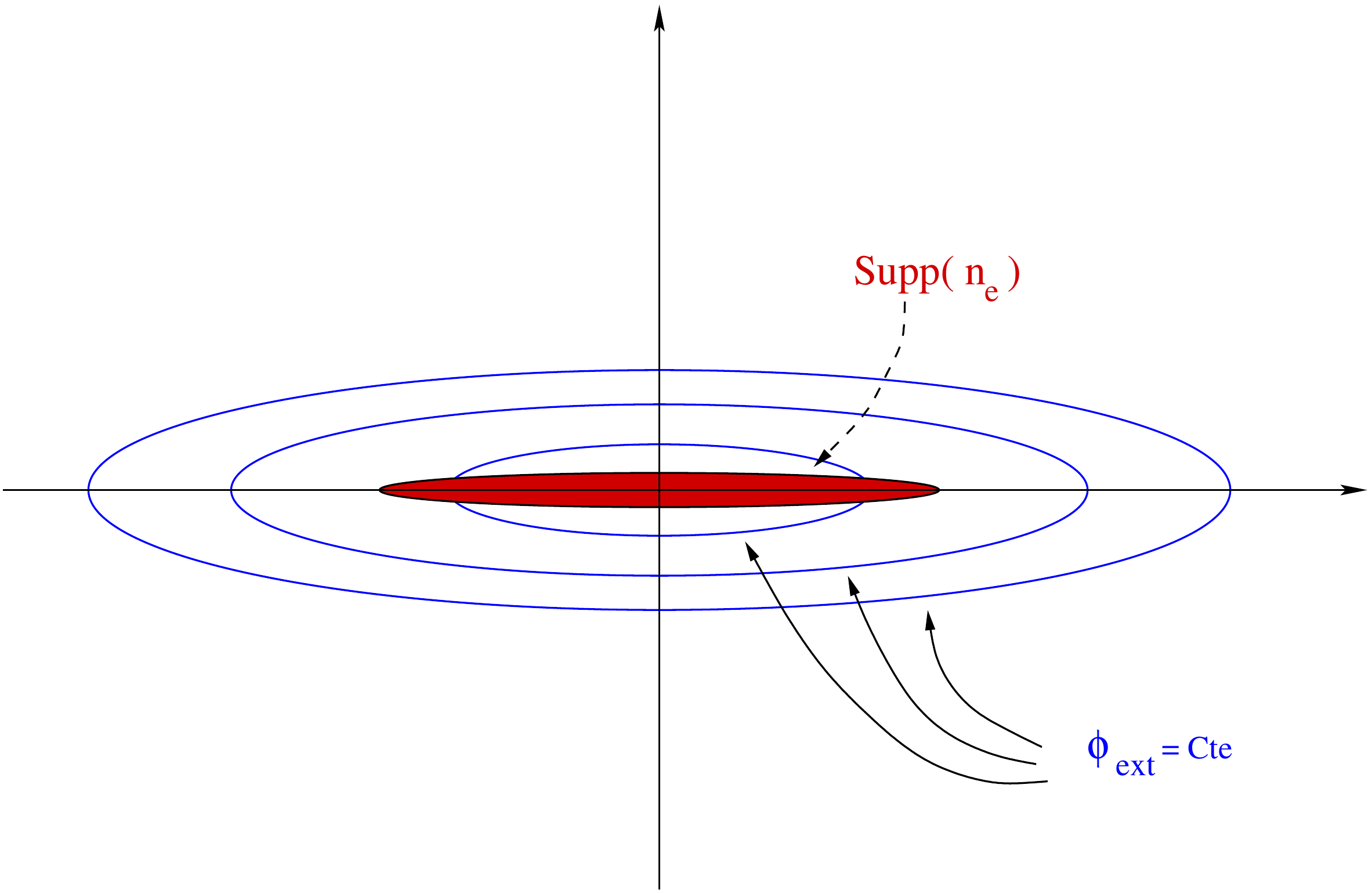}
\end{center}
\caption{Some level sets of $ \Phi_{\rm ext} $ and the support of $ n_{\rm e}$.}
\label{patator}
\end{figure}

We wish to extend this analysis by dealing with more general external
potentials.  We distinguish two situations depending on the expression of the 
external potential:
\begin{itemize}
\item The quadratic potential
\begin{equation}\label{quadra}
 \Phi_{\rm ext} (x) = \frac12 \sum_{j=1}^N \frac{x_j^2}{\lambda_j^2} ,
\end{equation}
with $ \lambda_j > 0 $, $1 \le j \le N $, in dimension $ N \ge 2 $ is
typical to model non neutral plasmas \cite{Dubinreview} or
magneto-optical traps experiments.  In this case $\Delta \Phi_{\rm
  ext}$ is still a constant, that therefore determines the value of
the (uniform) particle density $n_{\rm e}$ on its support.  But the
problem has lost its symmetries and the shape of the support becomes
non trivial.  We shall see that $\rho_\eps$ tends to a uniform
distribution $n_{\rm e}$, supported in an ellipsoid.  However, we
point out that the support of $n_{\rm e}$ does not coincide with a
level set of $\Phi_{\rm ext}$.  An example with $N=2$ is given in
Figure~\ref{patator}. The potential $\Phi_{\rm e}$ can be computed
rather explicitly, and Theorem \ref{thresultradial} generalizes
directly. See Section \ref{sec:patate} for a precise statement.

\item In the case of a non quadratic potential, under suitable
  hypotheses on $\Phi_{\rm ext}$, the limiting density $n_{\rm e}$
  still has a compact support $\cak$ and is still given on $\cak$
  by $n_{\rm e}=\Delta \Phi_{\rm ext}$. However, $ n_{\rm e} $ is 
  clearly no longer constant on $\cak$.  The identification of $\cak$ and 
  $n_{\rm e}$ relies on variational techniques, with connection to the
  obstacle problem.  It is still possible to prove the analog of
  Theorem~\ref{thresultradial}, but, since $n_{\rm e}$ becomes non
  homogeneous, instead of \eqref{Euler} the limiting equations are now
  the so-called Lake Equations, see e.~g. \cite{Lever-Oli-Titi_PhysD}:
\begin{equation}
\tag{LE}\label{lake}
\left\{ \begin{array}{l}
\partial_t V + V \cdot \nabla_x V + \nabla_x p = 0 , \\
\nabla_x\cdot (n_{\rm e} V)=0.
\end{array} \right.
\end{equation} 
Such model --- also referred to as the Anelastic Equations --- arise
in the modelling of atmospheric flows \cite{OgPh}; we refer the reader
to \cite{Mas2} for the justification of a derivation from the
compressible Navier-Stokes system.  As a matter of fact, we can
observe that the first equation in \eqref{lake} may be writen in the
following conservative form $ \partial_t (n_{\rm e} V) + \nabla_x
\cdot ( n_{\rm e} V \otimes V ) + n_{\rm e} \nabla_x p = 0 $.  The
construction of $\Phi_{\rm e}$ and $\cak$, and a precise statement of
the corresponding convergence theorem can be found in Section
\ref{sec:lake}.
\end{itemize}

Motivated by actual experiments, we will also generalize the 
results to the case where a Fokker-Planck operator acting on
velocities is added to Eq. \eqref{VP1}. Our starting point then becomes:
\begin{equation}
\tag{VFP}
\label{FPV1}
\partial_t f_\eps+v\cdot\nabla_ x f_\eps-\nabla_x\Phi_\eps\cdot\nabla_v f_\eps= 
Lf_\eps ,
\end{equation}
with 
\[
Lf=\nabla_v\cdot (vf+\theta \nabla_vf) =
\theta \nabla_v\cdot\left(M_{0,\theta}\nabla_v\ds\left(\frac{f}{M_{0,\theta}} \right)\right),\qquad
M_{0,\theta} (v)=\ds\frac{1}{(2\pi\theta)^{N/2}}e^{-\lvert v \rvert^2 /(2\theta)},
\] 
for some $\theta>0$.  Equation \eqref{FPV1} is still coupled to the
Poisson equation \eqref{VP2}.  Using a modified modulated energy, we
are able to show in this case that solutions $f_\varepsilon$ of
\eqref{FPV1} and \eqref{VP2} with well-prepared data converge when
$\varepsilon$ and $\theta$ tends to $0$ in the sense of Theorem
\ref{thresultradial} to $n_{\rm e} V$, where $V$ is now the solution
of the Lake Equation with friction
\begin{equation}
\label{lakefriction}
\left\{ \begin{array}{l}
\partial_t V + V \cdot \nabla_x V + \nabla_x p + V =0, \\
\nabla_x\cdot (n_{\rm e} V)=0 .
\end{array} \right.
\end{equation}
On the boundary, we still have the no-flux condition \eqref{CLEuler}. 
For the sake of completeness, the necessary analog of
Theorem~\ref{Eulerexiste} for the systems \eqref{lake} and
\eqref{lakefriction} is sketched in appendix \ref{sec:eq_lake}.
See Section \ref{sec:VFP} for a precise statement on the asymptotic
behavior of \eqref{FPV1} and its proof.


\section{The limit density $n_{\rm e}$ and total potential $\Phi_{\rm
    e}$}

As said above, we have a clear intuition and explicit formulae for the
equilibrium distribution $n_{\rm e}$ and the potential $\Phi_{\rm e}$
in the specific case of the isotropic external potential \eqref{iso}.
Let us discuss in further details how $\Phi_{\rm {ext}}$ determines
$n_{\rm e}$ and its support, and how the auxiliary potential
$\Phi_{\rm e}$, which plays a crucial role in the analysis through the
decomposition \eqref{VP2psi}, can be defined.

We remind the reader the definition of  the 
 fundamental solution, hereafter denoted $\Gamma$, of $(-\Delta)$ (mind
the sign) in the whole space $\mathbb{R}^N$:
\begin{equation}\label{Gamma}
\Gamma(x) \ddef \left\{ 
\begin{array}{ll}
\ds \frac{1}{N(N-2)|B_{\mathbb{R}^N} (0,1) | \cdot |x|^{N-2}} \quad & {\rm if} \ N>2,\\
\ds -\frac{ \ln |x|}{2\pi} &{\rm if}\ N=2 ,\\
\ds -\frac{|x|}{2} &{\rm if}\ N=1.
\end{array} \right.
\end{equation}

\subsection{The case of a general quadratic potential}
\label{sec:patate}

Let us consider in this section the case of a quadratic potential
\eqref{quadra}.  We have $ \Delta \Phi_{\rm ext} = \sum_{j=1}^N
\frac{1}{\lambda_j^2} > 0 $ which is constant in space.  It gives the
value of the equilibrium density on its support since we still expect
$ \rho_\eps \to {\bf 1}_{\cak}  \Delta \Phi_{\rm ext}$.  But it
remains to determine this support $\mathrm{Supp}(n_{\rm e})= \cak
\subset \R^N $ on which we have the volume constraint
$$
 \mathfrak{m} = \int n_{\rm e} \ud x= | \cak | \sum_{j=1}^N \frac{1}{\lambda_j^2} 
$$
coming from \eqref{alacharge}.  Note that a quick computation reveals
that $\cak$ can be neither radially symmetric, nor a level set of
$\Phi_{\rm ext}$.\\

In order to extend Theorem~\ref{thresultradial} for a potential as in
\eqref{quadra}, we need to construct a domain $ \cak \subset \R^N $
and a function $ \Phi_{\rm e} : \R^N \to \R $ such that
\begin{equation}\label{defPhie0patate}
\Delta \Phi_{\mathrm e}(x) = 
\left( \sum_{j=1}^N \frac{1}{\lambda_j^2} \right) {\bf 1}_{ \R^N \setminus \cak } ,
\quad \quad \quad \Phi_{\mathrm e} = 0 \quad {\rm in} \ \cak.
\end{equation}
The starting point is the observation that given $ a = (a_1 , ... , a_N ) \in
( \R_+^*)^N $, then the characteristic function of the ellipsoid $$
\cak_a = \{ x\in \R^N ; \ \sum_{j=1}^N x_j^2 / a_j^2 \le 1 \} $$
generates an electric potential which is {\it quadratic} inside the
ellipsoid.  This can be found for instance in \cite[Chapter
VII, $\S$ 6]{Kellogg}; the computation there is for $N=3 $, but the extension to
the case $ N \ge 3 $ is straightforward, and the two-dimensional case 
is treated by using  arguments from complex analysis in 
 \cite{Furman_1994}.
 
 For $x \in \R^N $, we denote by $ \sigma_a (x) $ the largest solution
 of the equation $$ \ds \sum_{j=1}^N \frac{x_j^2}{a_j^2 + \varsigma }
 = 1 $$ (with $ \varsigma \in \R $ as unknown).  Consequently, $ x \in
 \cak_a $ holds if and only if $ \sigma_a (x) \le 0 $.  By convention,
 $ \sigma_a (0)= -\infty $. This quantity can be seen as an equivalent
 of the radial coordinate in the ellipsoidal coordinate system.  It
 allows us to construct a solution to \eqref{defPhie0patate} where $
 \cak $ is an ellipsoid, the coefficients of which depend on the
 mass $ \mathfrak{m} $ and the $ \lambda_j $'s.

\begin{prop} 
\label{potardellipse} 
Let $ a = (a_1 , ... , a_N ) \in ( \R_+^*)^N $.\\ 
(i)  \cite{Kellogg} If $ N \ge 3 $, then
$$
\Gamma \star {\bf 1}_{\cak_a} (x) = \frac14 \left( \prod_{j=1}^N a_j
\right) \times \left\{\begin{array}{ll} \ds \int_{\sigma_a ( x)
    }^{+\infty} \left( 1 - \sum_{j=1}^N \frac{x_j^2}{a_j^2 + s
      }\right) \, \left( \ds \prod_{j=1}^N ( a_j^2 + s )
    \right)^{-1/2} \, \ud s & \quad {\it if} \ \sigma_a (x) \ge 0 ,
    \\
    \ds \int_{0 }^{+\infty} \left( 1 - \sum_{j=1}^N \frac{x_j^2}{a_j^2
        + s }\right) \, \left( \ds \prod_{j=1}^N ( a_j^2 + s )
    \right)^{-1/2} \, \ud s & \quad {\it if} \ \sigma_a (x) \le 0 .
\end{array}\right.
$$
$ (ii ) $ If $ N= 2 $, then 
$$
 \Gamma \star {\bf 1}_{\cak_a} (x) = \frac14 \left( a_1 a_2 \right) \times
 \left\{\begin{array}{ll}
\ds 
- \ln \left( \sigma_a(x) + \frac{a_1^2 + a_2^2}{2} + 
\sqrt{( a_1^2 + \sigma_a(x) )( a_2^2 + \sigma_a(x) )} \right) 
& \\ 
\quad \quad \ds - \int_{\sigma_a ( x) }^{+\infty} \sum_{j=1}^2 \frac{x_j^2}{a_j^2 + s } 
\, \frac{\ud s }{\sqrt{ ( a_1^2 + s )( a_2^2 + s )}}
& \quad {\it if} \ \sigma_a (x) \ge 0 ,
\\
\ds - \ln \left( \frac12 ( a_1 + a_2 )^2 \right) 
- \int_{0 }^{+\infty}\sum_{j=1}^2 \frac{x_j^2}{a_j^2 + s } 
\, \frac{\ud s }{\sqrt{ ( a_1^2 + s )( a_2^2 + s )}} 
& \quad {\it if} \ \sigma_a (x) \le 0 .
\end{array}\right.
$$
\end{prop}

\begin{remark} An alternative point of view for the two
  dimensional case is to work with the electric field instead of the
  potential. We refer to \cite{Furman_1994} for expressions of the
  electric field generated by ellipses in $N=2$. In the case of a
  uniform charge distribution, the electric field is linear inside the
  ellipse, with the same coefficients for the quadratic terms as those
  coming from the expression in (ii).
\end{remark}

We define the mapping $ \mathcal{Z} : (\R_+^*)^N \to (\R_+^*)^N $ by
\begin{equation}
\label{defGrandZ}
 \mathcal{Z}_j ( \alpha ) = \int_{0 }^{+\infty} \frac{1}{ \alpha_j + s } 
 \, \left( \ds \prod_{j=1}^N ( \alpha_j + s ) \right)^{-1/2} \, \ud s > 0 .
\end{equation}
From Proposition \ref{potardellipse}, we know that the potential
generated by $ {\bf 1}_{\cak_a} $ is quadratic inside $ \cak_a $,
up to an additive constant.  The coefficients of the quadratic terms
are the $ - ( \prod_{j=1}^N a_j ) \mathcal{Z}_k ( \alpha ) /4 $, $ 1
\le k \le N $.  The idea to make the connexion with the external
potential $ \Phi_\ext $ is now to adapt the $a_j$'s so that the
quadratic terms in $ \Gamma \star {\bf 1}_{\cak_a} $ (inside
$\cak_a$) cancel out the quadratic terms of $ \Phi_\ext $, so that $
(\Delta \Phi_\ext)\Gamma \star {\bf 1}_{\cak_a} + \Phi_\ext $ is
constant in $ \cak_a $. We observe that $ \prod_{j=1}^N a_j $ is
related to the total charge of the ellipsoid $ \cak_a$ since
$$ \mathfrak{m} = \int n_{\rm e} = | \cak_a |\sum_{j=1}^N \lambda_j^{-2} 
= | B_{\R^N} (0,1) | \left( \prod_{j=1}^N a_j \right) \sum_{j=1}^N \lambda_j^{-2} .
$$
We shall thus need to solve equations in $a$ of the form $ \mathcal{Z}
( a_1^2 , ... , a_N^2 ) = z $, where $z \in (\R_+^*)^N$ is given. 
Therefore, we are interested in showing that $ \mathcal{Z} : (\R_+^*)^N \to (\R_+^*)^N $
is a smooth diffeomorphism. When $N=2$, explicit computations may be
carried out.

\begin{prop} 
\label{diffeo2} 
Assume $ N = 2 $. Then, for any $ \alpha \in ( \R_+^*)^2 $
$$ \mathcal{Z} (\alpha ) = ( \mathcal{Z}_1 (\alpha ) , \mathcal{Z}_2 (\alpha ) ) 
= \left( \frac{2}{ \alpha_1 + \sqrt{\alpha_1 \alpha_2}}, \frac{2}{
    \alpha_2 + \sqrt{\alpha_1 \alpha_2}}\right) .
$$
Moreover, $ \mathcal{Z} : (\R_+^*)^2 \to (\R_+^*)^2 $ is a smooth diffeomorphism 
and its  inverse is given by
$$
\mathcal{Z}^{-1} ( z ) = ( (\mathcal{Z}^{-1})_1 (z ) ,
(\mathcal{Z}^{-1})_2 (z ) ) = \left( \frac{2 z_2}{ z_1 ( z_1 + z_2) },
  \frac{2 z_1}{ z_2 ( z_1 + z_2) } \right) .
$$
\end{prop}

\noindent {\it Proof.} The explicit formula for $\mathcal{Z} (\alpha )$ comes by  computing the
Abelian integral $$ \ds \int_{0 }^{+\infty} \frac{ \ud s}{ (\alpha_1 +
  s)^{3/2} ( \alpha_2 + s )^{1/2}} = \int_{0 }^{+\infty} \frac{\ud}{
  \ud s} \left( \frac{2}{\alpha_2 - \alpha_1 } \sqrt{ \frac{\alpha_1 +
      s}{\alpha_2 + s}} \right) \ \ud s = \frac{2}{ \alpha_1 +
  \sqrt{\alpha_1 \alpha_2}}, $$ for $ \alpha_1 \not = \alpha_2 $, and
the formula holds true when $ \alpha_1 = \alpha_2 $ as well. The formula for
the inverse then follows by direct substitution. \qed

For $ N \ge 3 $, we no longer have simple expressions for $\mathcal{Z}
$.  However, we shall prove that $ \mathcal{Z} : (\R_+^*)^N \to
(\R_+^*)^N $ is a smooth diffeomorphism by using the fact that $
\mathcal{Z} : (\R_+^*)^N \to (\R_+^*)^N $ is a gradient vector field
associated with a strictly concave function.

\begin{prop} 
\label{diffeo3etplus} 
Assume $ N \ge 2 $ and let us define the function $ \zeta : (\R_+^*)^N
\to \R $ by:
$$
\zeta ( \alpha ) = \left\{\begin{array}{ll} - \ds \int_0^{+\infty}
    \left( \prod_{k=1}^N ( \alpha_k + s ) \right)^{-1/2} \, \ud s &
    \quad {\it if} \ N \ge 3
    \\
    4 \ln ( \sqrt{\alpha_1} + \sqrt{\alpha_2} ) & \quad {\it if} \ N =
    2 .
\end{array}\right.
$$
Then, $ \zeta : (\R_+^*)^N \to \R $ is smooth, strictly concave and it satisfies 
$ \nabla \zeta = \mathcal{Z} $. Furthermore, $ \nabla \zeta =
\mathcal{Z} : (\R_+^*)^N \to (\R_+^*)^N $ is a smooth diffeomorphism
and for any $ z \in (\R_+^*)^N $, $ \mathcal{Z}^{-1} ( z ) $ is the
unique minimizer for
\begin{equation}
\label{transfoLegendre}
 \inf_{ \alpha \in (\R_+^*)^N } \left( z \cdot \alpha - \zeta ( \alpha ) \right ).
\end{equation}
\end{prop}

In \eqref{transfoLegendre} we recognize the minimization problem that defines the Legendre transform 
of $ \zeta $. This gives a way to compute numerically $ \mathcal{Z}^{-1} ( z ) $ through the
minimization of a convex function.
\\

\noindent {\it Proof.} 
The smoothness of $ \zeta $ is clear and $ \nabla \zeta = \mathcal{Z}
$ follows from direct computations.  If $ N = 2 $, the strict
concavity of $ \zeta $ is straightforward and the fact that $ \nabla
\zeta = \mathcal{Z} : (\R_+^*)^2 \to (\R_+^*)^2 $ is a smooth
diffeomorphism comes from Proposition \ref{diffeo2}: for any $ z \in
(\R_+^*)^2 $, $ \mathcal{Z}^{-1} (z) $ is a critical point of the
strictly convex (since $ \zeta $ is strictly concave) function $
\alpha \mapsto z \cdot \alpha - \zeta ( \alpha ) $, hence is the
unique minimizer of that function. We assume now $ N \ge 3 $.  Then,
for each $ s \in \R_+ $, the function $$ \varpi_s : \alpha\in (\R_+^*)^N 
 \mapsto \left( \prod_{k=1}^N ( \alpha_k + s ) \right)^{-1/2} $$ is
logarithmically strictly convex since $ \ln \circ \varpi_s ( \alpha )
= (-1/2) \sum_{k=1}^N \ln( \alpha_k + s ) $ and $ {\rm Hess}( \ln
\circ \varpi_s , \alpha ) = (1/2) {\rm Diag} ( (\alpha_1 + s )^{-2} ,
... , (\alpha_N + s)^{-2} ) $.  Consequently, $ - \zeta ( \alpha ) =
\int_0^{+\infty} \varpi_s ( \alpha ) \ud s $ is a strictly convex
function of $\alpha $. Let us show that the Jacobian determinant of $
\mathcal{Z} $ never vanishes, that is $ {\rm Hess}( \zeta , \alpha ) $
is everywhere negative definite.  For that purpose, for $v \in \R^N $,
we write $ - v^T {\rm Hess}( \zeta , \alpha ) v = \int_0^{+\infty} v^T
{\rm Hess}(\varpi_s , \alpha ) v \ud s $, and thus it suffices to show
that $ {\rm Hess}( \varpi_s , \alpha ) $ is positive definite for any
$ s \ge 0 $. Now, we write $ \varpi_s ( \alpha ) = \exp( \ln \circ
\varpi_s ( \alpha ) ) $, thus $ \p^2_{j,k} \varpi_s ( \alpha )= \exp(
\ln \circ \varpi_s( \alpha ) ) [ \p^2_{j,k} ( \ln \circ \varpi_s ) (
\alpha ) + \p_j (\ln \circ \varpi_s)( \alpha ) \p_k ( \ln \circ \varpi_s )(\alpha) ] $.
Therefore, if $v \not = 0 $, we obtain
\begin{align*}
v^T {\rm Hess}( \varpi_s , \alpha ) v 
= & \, \varpi_s ( \alpha ) 
 \left[ v^T {\rm Hess}( \ln \circ \varpi_s , \alpha ) v 
 + \left( \sum_{j=1}^N v_j \p_j ( \ln \circ \varpi_s)( \alpha ) \right)^2 
 \right] \\
\ge & \, 
\varpi_s ( \alpha ) v^T {\rm Hess}( \ln \circ \varpi_s , \alpha ) v 
= \varpi_s ( \alpha ) \sum_{j=1}^N \frac{v_j^2}{2(\alpha_j + s)^2} > 0 ,
\end{align*}
as wished. 

Let us now fix $ z \in (\R_+^*)^N $ and consider the minimization
problem \eqref{transfoLegendre}. In view of the negativity of $\zeta$,
this infimum $ \mu $ belongs to $ [ 0, +\infty ) $. Since $ \zeta $ is
strictly concave, this problem has at most one minimizer. Let us show
that it has at least one by considering a minimizing sequence $ (
\alpha^n )_{ n\ge 0 } \in (\R_+^*)^N $. We claim that the sequence $ (
\alpha^n )_{ n\ge 0 } $ is bounded. Indeed, we have $ z \cdot \alpha^n
- \zeta ( \alpha^n ) \to \mu \in \R_+ $, and since $ \zeta \le 0 $,
this implies $ z \cdot \alpha^n = \zeta ( \alpha^n ) + \mu + o(1) \le
\mu + o(1) $. Using that all components of $z$ are positive, the claim
follows. As a consequence, we may assume, up to a subsequence, that
there exists $ \alpha = ( \alpha_1,...,\alpha_N) \in \R_+^N $ such that 
$ \alpha^n \to \alpha $ as $ n \to +\infty $. In particular, 
$ \zeta ( \alpha^n ) = z \cdot \alpha^n - \mu + o(1) $ converges. We now 
prove that at most two components of $ \alpha $ vanish.  For otherwise,
Fatou's lemma would yield
\begin{align*}
 +\infty = & \, 
 \int_0^{+\infty} \left( \prod_{k=1}^N ( \alpha_k + s ) \right)^{-1/2} \, \ud s 
 = \int_0^{+\infty} \liminf_{n \to +\infty } 
 \left( \prod_{k=1}^N ( \alpha_k^n + s ) \right)^{-1/2} \, \ud s 
 \\
 \le & \, \liminf_{n \to +\infty } \int_0^{+\infty}  
 \left( \prod_{k=1}^N ( \alpha_k^n + s ) \right)^{-1/2} \, \ud s 
 = \liminf_{n \to +\infty } ( - \zeta ( \alpha^n ) ) ,
\end{align*}
contradicting the convergence of $ ( \zeta ( \alpha^n ) )_{n\in \mathbb N} $.  It
remains to show that $ \alpha $ has no zero component to ensure that $
\mu + o(1) = z \cdot \alpha^n - \zeta ( \alpha^n ) \to z \cdot \alpha
- \zeta ( \alpha ) $ so that $ \alpha \in ( \R_+^*)^N $ is actually a
minimizer for \eqref{transfoLegendre}.  We then assume that $ \alpha_1
= 0 $, for instance, and show that for sufficiently small $ \delta > 0 $, $ z \cdot
( \delta , \alpha_2 , ... , \alpha_N) - \zeta ( \delta , \alpha_2 ,
... , \alpha_N ) < z \cdot ( 0 , \alpha_2 , ... , \alpha_N) - \zeta (
0 , \alpha_2 , ... , \alpha_N ) $.  This reaches a contradiction for
$n$ large enough. We thus compute
\begin{align*}
  D(\delta) = & \, \Big( z \cdot ( \delta , \alpha_2 , ... , \alpha_N)
  - \zeta ( \delta , \alpha_2 , ... , \alpha_N ) \Big) - \Big( z \cdot
  ( 0 , \alpha_2 , ... , \alpha_N) - \zeta ( 0 , \alpha_2 , ... ,
  \alpha_N ) \Big)
  \\
  = & \, z_1 \delta + \int_0^{+\infty} ( \delta + s )^{-1/2} \left(
    \prod_{k=2}^N ( \alpha_k + s ) \right)^{-1/2} \, \ud s -
  \int_0^{+\infty} s^{-1/2} \left( \prod_{k=2}^N ( \alpha_k + s )
  \right)^{-1/2} \, \ud s
  \\
  = & \, \delta \Big( z_1 - \int_0^{+\infty} \frac{ 1}{ s^{1/2}
    (\delta + s )^{1/2} [s^{1/2} + (\delta + s )^{1/2}] } \left(
    \prod_{k=2}^N ( \alpha_k + s ) \right)^{-1/2} \, \ud s \Big) .
\end{align*}
As $ \delta \to 0 $, we have, by monotone convergence, 
$$
\int_0^{+\infty} \frac{ \left( \prod_{k=2}^N ( \alpha_k + s )
  \right)^{-1/2}}{ s^{1/2} (\delta + s )^{1/2} [s^{1/2} + (\delta + s
  )^{1/2}] } \, \ud s \to \int_0^{+\infty} \frac{ 1}{ 2 s^{3/2} }
\left( \prod_{k=2}^N ( \alpha_k + s ) \right)^{-1/2} \, \ud s = +
\infty ,
$$
hence for $ \delta $ sufficiently small, $ D(\delta) < 0 $, as
claimed.  Therefore, $\alpha \in ( \R_+^*)^N $ and $ \alpha $ is a
minimizer for \eqref{transfoLegendre}.  It then follows that $ \nabla
\zeta (\alpha) = z $ as wished. \qed

We may now construct a solution to \eqref{defPhie0patate}.

\begin{coro} {\bf (Construction of the function $\Phi_{\rm e} $ for quadratic potentials).} 
\label{omega_patate}
  Let $N \ge 2 $ and assume that
$$
 \Phi_{\rm ext} (x) = \frac12 \sum_{j=1}^N \frac{x_j^2}{\lambda_j^2} ,
$$
with $ \lambda_j > 0 $, $1 \le j \le N $. Let us also fix $
\mathfrak{m} > 0 $.  Then, there exists a unique $ a \in (\R_+^*)^N $
such that
\begin{equation}
\label{Kondi}
| \cak_a | \left( \sum_{j=1}^N \frac{1}{\lambda_j^2} \right) = \mathfrak{m} 
\quad \quad \quad {\it and} \quad \quad \quad 
\frac{\mathfrak{m}}{2 |B_{\R^N}(0,1)| \sum_{k=1}^N \lambda_k^{-2} }\  \mathcal{Z} (a_1^2, ... , a_N^2) 
= \left( \frac{1}{\lambda_j^2} \right)_{ 1 \le j \le N}.
\end{equation}
Therefore, there exists a constant $ \kappa $,
depending only on the $ \lambda_j $'s, $N$ and $ \mathfrak{m} $ such
that the function $$ \ds \Phi_{\rm e} = \Phi_{\rm ext} + \left(
  \sum_{j=1}^N \frac{1}{\lambda_j^2} \right) \Gamma \star {\bf
  1}_{\cak_a} + \kappa $$ is convex and satisfies
\begin{equation}
\label{Bonnard}
 - \Delta \Phi_{\rm e} = \left( \sum_{j=1}^N \frac{1}{\lambda_j^2} \right) 
 {\bf 1}_{\R^N \setminus \cak_a}
 \text{ with, furthermore,  
 $\Phi_{\rm e} = 0$ in $\cak_a$ and 
$  \Phi_{\rm e} > 0$  in  $\R^N \setminus  \cak_a$} .
\end{equation}
\end{coro}

\noindent {\it Proof.} We define $ \lambda > 0 $ such that
$  \lambda^{-2} = \sum_{j=1}^N \lambda_j^{-2} $ and the constant $ \kappa $ by the
formulas $ 4 \kappa = - \lambda^{-2} \left( \prod_{j=1}^N a_j \right)
\int_0^{+\infty} \left( \prod_{j=1}^N ( a_j^2 + s ) \right)^{-1/2} \ \ud s $
if $ N \ge 3 $ and $ 4 \kappa = - \lambda^{-2} ( a_1 a_2 ) \ln ( (a_1 + a_2 )^2 /2 ) $ if $ N = 2 $.
The existence (and uniqueness) of $a$ satisfying the conditions \eqref{Kondi} then ensures that 
 $\Phi_{\rm e} = 0$ in $\cak_a$ and $ - \Delta \Phi_{\rm e} = \lambda^{-2} {\bf 1}_{\R^N \setminus \cak_a} $. 
Then, from the formulas in Proposition
\ref{potardellipse}, we get, in $ \{ \sigma_a > 0 \} $,
\begin{align}
\label{expresso}
 \frac{4 \lambda^{2}}{\prod_{j=1}^N a_j } \Phi_{\rm e} (x)
 = & \, 
 \Gamma \star {\bf 1}_{\cak_a} (x) + \frac{4 \lambda^{2}}{\prod_{j=1}^N a_j } \Phi_\ext (x)
 + \frac{4 \lambda^{2}}{\prod_{j=1}^N a_j } \kappa
 \nonumber \\
 = & \,
 \int_{\sigma_a(x)}^{+\infty}
 \left( 1-\sum_{j=1}^N \frac{x_j^2}{a_j^2+s} \right) \left( \prod_{j=1}^N (a_j^2 +s) \right)^{-1/2} \ \ud s
 \nonumber \\
 & \, + \sum_{k=1}^N x_k^2 \int_0^{+\infty} \left( \prod_{j=1}^N ( a_j^2 + s ) \right)^{-1/2} \
 \frac{\ud s }{a_j^2 + s}
 - \int_0^{+\infty} \left( \prod_{j=1}^N ( a_j^2 + s ) \right)^{-1/2} \ \ud s
 \nonumber \\
 = & \,
 \int_0^{\sigma_a(x)} \left( \sum_{j=1}^N \frac{x_j^2}{a_j^2 + s } - 1 \right)
    \left( \prod_{j=1}^N ( a_j^2 + s ) \right)^{-1/2} \ \ud s .
\end{align}
The last integral is positive if $ \sigma_a (x) > 0 $ since,
when $ 0 \le s < \sigma_a(x) $,
$ \sum_{j=1}^N \frac{x_j^2}{a_j^2 + s } - 1 > \sum_{j=1}^N \frac{x_j^2}{a_j^2 + \sigma_a(x) } - 1 = 0 $. 
In order to see that $ \Phi_{\rm e} $ is convex, we notice that $ \Phi_{\rm e} \equiv 0 $ 
in $ \{ \sigma_a \le 0 \} = \mathcal{K}_a $ and that, from \eqref{expresso}, we have, 
when $ \sigma_a(x) > 0 $, and for any direction $ \omega \in \mathbb S^{N-1} $, 
\begin{align*}
\p^2_\omega \Phi_{\rm e} (x) = 
& \, \frac{\prod_{j=1}^N a_j }{2 \lambda^{2}} 
\int_0^{\sigma_a(x)} \left( \sum_{j=1}^N \frac{\omega_j^2}{a_j^2 + s } \right)
    \left( \prod_{j=1}^N ( a_j^2 + s ) \right)^{-1/2} \ \ud s 
    \\ 
+ &\, 
 \frac{\prod_{j=1}^N a_j }{ \lambda^{2}} \left( \sum_{j=1}^N \frac{x_j \omega_j}{a_j^2 + \sigma_a(x) }\right)^2 
 \left( \sum_{j=1}^N \frac{x_j^2}{ ( a_j^2 + \sigma_a(x) )^2 }\right)
 \left( \prod_{j=1}^N ( a_j^2 + \sigma_a(x) ) \right)^{-1/2},
\end{align*}
since $ \sum_{j=1}^n x^2_j / ( a_j^2 + \sigma_a(x) ) = 1 $, which is indeed $ > 0 $. \qed

Clearly, the ellipsoid $\cak_a$ is not a level set of the external
potential $\Phi_\ext$ (except when all the $\lambda_j$'s are all equal). 
It is interesting to study the limiting case of
a very asymmetric external potential. For instance in $N=2$,
we consider a trapping potential \eqref{quadra} with a large aspect ratio
$A=\lambda_1 /\lambda_2\gg 1$. Direct computations (using Proposition \ref{diffeo2}) lead to 
\[
a_1 = \sqrt{\frac{\mathfrak{m}}{\pi}}
\frac{\lambda_1}{\sqrt{1+\frac{\lambda_2^2}{\lambda_1^2}}}~;~ a_2 =
\sqrt{\frac{\mathfrak{m}}{\pi}}
\frac{\lambda_2}{\sqrt{1+\frac{\lambda_1^2}{\lambda_2^2}}}.
\]
Hence
\[
\frac{a_1}{a_2} = \frac{\lambda_1^2}{\lambda_2^2} = A^2
\]
Thus, the aspect ratio of the particles' cloud is much larger than the
aspect ratio of the external potential: this is an effect of the
strong repulsion, see Figure~\ref{patator} for a typical picture. A
similar phenomenon occurs in higher dimensions. For $N=3$ with
cylindrical symmetry, explicit formulae corresponding to our
$\mathcal{Z}$ function are given for instance in
\cite{Wineland1985}. It is easy to check that for a strongly oblate
external potential (``pancake shape''), the aspect ratio of the cloud
is again of the order of the square of the aspect ratio of the
external potential.  We can now state the analog of
Theorem~\ref{thresultradial} for a general quadratic $\Phi_\ext$.
\begin{theorem}
\label{thresultpatate} 
Let $\Phi_\ext$ be any quadratic potential \eqref{quadra} to which we
associate, by virtue of Corollary~\ref{omega_patate}, the ellipsoid
$\cak_a$ and the potential $\Phi_{\rm e}$.  Let $ V^{\mathrm{init}}
\in H^s(\stackrel{\circ}{\cak}_a)$ satisfy $\nabla_x\cdot
V^{\mathrm{init}}=0$ in $ \stackrel{\circ}{\cak}_a$ and the no flux
condition \eqref{CLEuler} on $\partial \cak_a$.  Denote by $V$ the
solution on $ [ 0 , T ] $ to \eqref{Euler} with the no flux condition
\eqref{CLEuler} given in Theorem \ref{Eulerexiste} and consider
$\mathscr V^{\mathrm{init}}$ a smooth extension of $V$ in $\mathbb{R}^N $ 
satisfying the
following conditions, where $R> 0 $ is such that $\cak_a\subset B(0,R)$,
\[
\mathscr V\Big|_{\cak_a}=V,\qquad
\mathscr V\Big|_{\mathbb R^N\setminus B(0,2R)}=0,
\qquad
\mathscr V(t,x)\cdot \nu(x)\Big|_{\partial\cak_a}=0 .
\]
Let $f_\eps^{\mathrm{init}}:\mathbb R^N\times
\mathbb R^N\rightarrow [0,\infty)$ be a sequence of integrable
functions that satisfy \eqref{query}.
Then, the associated solution $f_\eps$ of the Vlasov--Poisson equation 
\eqref{VP1}--\eqref{VP2} satisfies, 
as $ \eps\rightarrow 0 $,
\begin{itemize}
\item[i)] $\rho_\eps$ converges to $n_{\mathrm e} = \left( \sum_{j=1}^N \lambda_j^{-2} \right) 
 {\bf 1}_{\cak_a}$ in 
$C^0( 0,T ;\mathscr M^1(\mathbb R^N)-\text{weak}-\star)$;
\item[ii)] $ \mathscr H_{\mathscr V,\eps}$ converges to $0 $ uniformly
  on $[ 0, T ]$;
\item[iii)] $J_\eps$ converges to $J$ in $\mathscr M^1([0,T]\times
  \mathbb R^N)$ weakly-$\star $, the limit $J$ lies in
  $L^\infty(0,T;L^2(\mathbb R^N))$ and satisfies $J\big|_{[0,T]\times
    \cak_a}=V$, $\nabla_x\cdot J=0$ and $J\cdot \nu(x)\big|_{\partial\cak_a}=0$.
    \end{itemize}
\end{theorem}

The existence of a smooth extension $\mathscr V$ of $V$ on
$\mathbb{R}^N $ satisfying the above mentioned constraints follows
from~\cite[Chapter I: Theorem 2.1 p. 17 \& Theorem 8.1 p. 42]{LM}. See
Lemma \ref{extension} in the appendix for a divergence-free extension.

\subsection{The case of a general potential}
\label{sec:lake}

We wish now to extend the above results to a general confining
potential.  When $\Phi_{\mathrm{ext}}$ is not quadratic, the
equilibrium density $n_{\rm e}$ cannot be expected to be constant on
its support. In turn, the limiting equation will be more 
complicated than the Incompressible Euler system.  Besides, the
determination of the domain $ \cak = \{ \Phi_{\rm e} = 0 \} $ is a non 
trivial issue, and its geometry might be quite involved \cite{Cours_Serfaty}. 
In the following we write $\Omega=\mathring{\cak}$ for the interior of $\cak$.

The pair $(n_{\rm e},\Omega)$ should be thought of through energetic consideration. 
As it will be detailed below, the total energy of the  system \eqref{VP1}-\eqref{VP2} is
\[
\iint \frac{|v|^2}{2}f_\eps \ud v\ud x+\frac{1}{2\eps} \int \Phi_\eps \rho_\eps\ud x +
\frac{1}{\eps} \int \Phi_\ext \rho_\eps\ud x.
\]
It is natural to investigate solutions whose energy does not diverge
when $\eps$ tends to $0$. Hence we are interested in configurations
close to the ground state $n_{\mathrm{e}}$ defined by the variational
problem where only the electrostatic part of the energy is involved:
namely, we wish to minimize
\[
\mathcal{E}[\rho] \ddef \int \Phi_\ext(x) \ud \rho(x) + \frac12 \iint
\Gamma(x-y) \ud \rho(y) \ud \rho(x) ,
\]
for a fixed $\mathfrak m>0$
over the convex subset  
$\mathscr{M}_{\mathrm{ext}}^+(\mathfrak m)$
made of nonnegative Borel 
measures $\rho$ of total mass $\mathfrak m>0$ such that $\int \Phi_{\mathrm{ext}}\ud \rho $ is 
finite.
This problem, which is often referred to as the generalized Gauss variational problem, is quite 
classical and the basis of the theory dates back to \cite{Fro}. 
We refer the reader to \cite[Chapter 1]{ST}  
for the case $ N=2 $,  and to \cite[Theorem 1.2]{ChaGoZi}  when $ N \ge 3 $ 
for the existence of a minimizer under suitable assumptions on $ \Phi_\ext $.
In what follows, we shall assume that $ \Phi_{\rm ext} $ fulfils the following requirements:
\begin{itemize}
\item[h1)] $ \Phi_{\rm ext} : \R^N \to \R_+ $ is 
continuous, nonnegative and satisfies $ \Phi_{\rm ext} (x) \to + \infty $ 
as $ |x| \to + \infty $,
\item[h2)] If $ N=2 $ or $N=1$, we have 
$ \ds \lim_{|x| \to +\infty} ( \Phi_{\rm ext} +  \mathfrak{m} \Gamma ) (x) = +\infty $.
\end{itemize}
The following statement collects from  \cite{ChaGoZi, ST, Cours_Serfaty} the 
results we shall need.

\begin{theorem} 
\label{Minimizor} We assume that the potential  $ \Phi_{\rm ext} $ satisfies the  hypotheses h1) and h2).\\
$( i) $ The functional $ \mathcal{E} $ is strictly convex on $ \mathscr{M}_{\mathrm{ext}}^+(\mathfrak m) $.\\
$ (ii) $ The problem
\begin{equation}
\label{eq:var}
\inf \left\{\mathcal{E}[\rho]~;~\rho \in \mathscr{M}_{\mathrm{ext}}^+(\mathfrak m)
 \right\}
\end{equation}
 has a unique minimizer $ n_{\rm e} $ which has a compact support of positive capacity.
Moreover, there exists a constant $ C_* $ such that
\begin{equation}
 \label{Robin}
 \left\{ 
\begin{array}{ll}
 \Gamma \star n_{\rm e}  + \Phi_{\rm ext} \ge C_* \quad {\it quasi\ everywhere} ,\\
 \Gamma \star n_{\rm e}  + \Phi_{\rm ext} = C_* \quad {\it quasi\ everywhere\ on} \ {\rm Supp}(n_{\rm e}).
 \end{array} \right.
\end{equation}
(iii) Conversely, assume that $ \rho_0 \in \mathscr{M}_{\mathrm{ext}}^+ (\mathfrak{m} ) $ 
and $ C_0 $ are such that 
$$
\left\{ 
\begin{array}{ll}
 \Gamma \star \rho_0 + \Phi_{\rm ext} \ge C_0 \quad {\it quasi\ everywhere} ,\\
 \Gamma \star \rho_0 + \Phi_{\rm ext} = C_0 \quad {\it quasi\ everywhere \ on} \ {\rm Supp}(\rho_0) .
\end{array} \right.
$$
Then, $ \rho_0 $ is the minimizer for \eqref{eq:var}: $\rho_0 =n_{\rm e} $.
\end{theorem}

We then define the potential $\Phi_{\rm e} \ddef 
\Gamma \star n_{\rm e} + \Phi_\ext - C_{*} $. 
The constant $ C_* $ in \eqref{Robin} is called the modified Robin constant 
and {\it quasi everywhere} (q.~e.) means up a set of zero capacity (which is 
a bit stronger than to be Lebesgue-negligible); see 
\cite[Definition 2.11]{Cours_Serfaty}. If $N=1$, \eqref{Robin} holds pointwise.\\

\noindent
{\it Proof.}
The statements for $ N \ge 3 $ can be found in \cite[Theorem 1.2]{ChaGoZi}. When $ N = 2 $, we refer 
the reader to  \cite[Theorem 1.3 for (ii) 
and Theorem 3.3 for (iii)]{ST}. If $N=2$, the strict convexity
(i) is not explicited in \cite{ST}. Thus we give proofs of (i) for
$N=2$, and (i)-(iii) for $N=1$.

The argument for (i) is that if $ \rho_0
$, $ \rho_1 \in \mathscr{M}_{\mathrm{ext}}^+(\mathfrak m) $ and $ \theta \in ( 0,1 )
$, then
\begin{align*}
\mathcal{E}[ (1-\theta) \rho_0 & \, + \theta \rho_1 ] 
- (1-\theta) \mathcal{E}[ \rho_0 ] 
- \theta \mathcal{E}[ \rho_1 ] 
\\
= & \, \frac12 \iint \Gamma(x-y) 
\ud [ (1-\theta) \rho_0 + \theta \rho_1 ] (y) \ud [ (1-\theta) \rho_0 + \theta \rho_1 ](x) \\
 & \,
- \frac12 (1-\theta) \iint \Gamma(x-y) \ud \rho_0(y) \ud \rho_0 (x) 
- \frac12 \theta \iint \Gamma(x-y) \ud \rho_1(y) \ud \rho_1(x) 
\\ 
= & \, - \frac12 \theta ( 1 - \theta ) \iint \Gamma(x-y) 
\ud [ \rho_0 - \rho_1](y) \ud [ \rho_0 - \rho_1 ] (x).
\end{align*}
Unless $ \rho_0 = \rho_1 $, the last integral is shown to be positive
if $ N \ge 3 $  in  \cite[Lemma 3.1]{ChaGoZi}. The case 
$N=2 $ is dealt with in    \cite[Lemma 1.8]{ST}, under the restriction that $
\rho_0 - \rho_1 $ has compact support.
Actually, the method used in \cite{ChaGoZi}, which consists in writing 
$ \Gamma (x) $ as an integral of Gaussians $ {\sf e}^{ -|x|^2 /2t } $, can be 
extended to the case $ N= 2 $ as we check now. The starting point is the 
equality (see \cite[equation (12)]{BenArousGuionnet})
$$
\ln \frac{1}{r} = \int_0^{+\infty} \frac{1}{2t} \left( 
    {\sf e}^{-r^2/2t} - {\sf e}^{-1/2t}\right) \, \ud t .
$$
Therefore, denoting $ \rho \ddef \rho_0 - \rho_1 $ and $ r \ddef |x-y|
$ and using the dominated convergence theorem (on each of the sets $
\{ |x-y| < 1 \} $ and $ \{ |x-y| \ge 1 \} $), we obtain 
\begin{align*}
  \iint \ln \frac{1}{ |x-y|} \ud \rho(y) \ud \rho (x) = & \, \lim_{T
    \to +\infty } \int_{1/T}^{T} \frac{1}{2t} \iint \left( {\sf
      e}^{-r^2/2t} - {\sf e}^{-1/2t}\right) \ud \rho(y) \ud \rho (x)
  \, \ud t
  \\
  = & \, \lim_{T \to +\infty } \int_{1/T}^{T} \frac{1}{2t} \iint {\sf
    e}^{-r^2/2t} \ud \rho(y) \ud \rho (x) \, \ud t
  \\
  = & \, \lim_{T \to +\infty } \int_{1/T}^{T} \frac{1}{4\pi} \iint
  \int {\sf e}^{-t |\xi|^2/2 - i \xi \cdot ( x-y) } \ud \xi \ud
  \rho(y) \ud \rho (x) \, \ud t
  \\
  = & \, \lim_{T \to +\infty } \int_{1/T}^{T} \int \frac{1}{4\pi} {\sf
    e}^{-t |\xi|^2/2 } \lvert \hat{\rho} (\xi) \rvert^2 \, \ud \xi \ud
  t
  \\
  = & \, \int \frac{1}{2\pi |\xi|^2} \lvert \hat{\rho} (\xi)
  \rvert^2\, \ud \xi ,
\end{align*}
where, for the second equality, we use $ \int \ud \rho = 0 $, and
for the third one, we write $ {\sf e}^{-r^2/2t} $ as the Fourier
transform of a two dimensional Gaussian. This clearly shows that $
\iint \Gamma (x-y) \ud \rho(y) \ud \rho (x) $ is positive unless $
\rho = 0 $, ensuring the strict convexity of $ \mathcal{E} $ on 
$\mathscr{M}_{\mathrm{ext}}^+(\mathfrak{m})$.
When $ N=1 $, we argue in a similar way by observing that
$$
- r = \int_0^{+\infty} \frac{1}{\sqrt{2\pi t}} \left( {\sf
    e}^{-r^2/2t} - 1 \right) \, \ud t .
$$
Indeed, $ {\sf e}^{-r^2/2t} - 1 = \int_0^r \p_u ({\sf e}^{-u^2/2t}) \, \ud u 
= - \int_0^r ( u/t ) {\sf e}^{-u^2/2t} \, \ud u $, thus
\begin{align*}
  - \int_0^{+\infty} \frac{1}{\sqrt{t}} \left( {\sf e}^{-r^2/2t} - 1
  \right) \, \ud t = & \, \int_0^{+\infty} \frac{1}{\sqrt{t}} \int_0^r
  ( u/t ) {\sf e}^{-u^2/2t} \, \ud u \ud t = \int_0^r \int_0^{+\infty}
  \frac{u}{t^{3/2}} {\sf e}^{-u^2/2t} \, \ud t \ud u
  \\
  = & \, \int_0^r \int_0^{+\infty} 2 \sqrt{2} {\sf e}^{-\tau^2} \, \ud
  \tau \ud u = \int_0^r
  \sqrt{2\pi} \ud u = r \sqrt{2\pi} ,
\end{align*}
where we have used the change of variable $ \tau = u/\sqrt{2t}$.
Owing to this relation, we can follow the same lines as above:
\begin{align*}
  - \frac12 \iint |x-y| \ud \rho(y) \ud \rho (x) = & \, \lim_{T \to +\infty }
  \int_{1/T}^{T} \frac{1}{2 \sqrt{2\pi t}} \iint \left( {\sf
      e}^{-r^2/2t} - 1 \right) \ud \rho(y) \ud \rho (x) \, \ud t
  \\
  = & \, \lim_{T \to +\infty } \int_{1/T}^{T} \frac{1}{2 \sqrt{2\pi t}}
  \iint {\sf e}^{-r^2/2t} \ud \rho(y) \ud \rho (x) \, \ud t
  \\
  = & \, \lim_{T \to +\infty } \int_{1/T}^{T} \frac{1}{4\pi} \iint
  \int {\sf e}^{-t \xi^2/2 - i \xi ( x-y) } \ud \xi \ud \rho(y) \ud
  \rho (x) \, \ud t
  \\
  = & \, \lim_{T \to +\infty } \int_{1/T}^{T} \int \frac{1}{4\pi} {\sf
    e}^{-t \xi^2/2 } \lvert \hat{\rho} (\xi) \rvert^2 \, \ud \xi \ud
  t
  \\
  = & \, \int \frac{1}{2\pi |\xi|^2} \lvert \hat{\rho} (\xi)
  \rvert^2\, \ud \xi .
\end{align*}
It only remains to prove (ii)-(iii) for $N=1$. This is tackled in 
\cite{Cours_Serfaty}, but with $\Gamma(x) =-\ln|x|$. The very same arguments
apply to the  case $\Gamma(x) =-|x|/2$.
\qed

\begin{remark} To motivate the above computation, one can remark that for 
$N \geq 1$ and under the condition $ \int \ud \rho = 0 $, we have, at least formally, 
$$   \iint \Gamma(x-y) 
\ud \rho (y) \ud  \rho (x)  = \int |\nabla \Delta^{-1} \rho |^2 dx 
= (2\pi)^{-N} \int \frac{1}{|\xi|^2} \lvert \hat{\rho} (\xi) \rvert^2\, \ud \xi . $$ 
\end{remark}  

The minimization of the functional $\mathcal{E}$ is connected to an {\it obstacle problem}. 
This connection is explained in details in \cite[Section 2.5]{Cours_Serfaty}.

\begin{prop}
\label{sautobstacle}
  If $n_{\rm e}$ is the minimizer of Theorem~\ref{Minimizor}, then
  $h=\Gamma \star n_{\rm e}$ is the unique solution to the
  obstacle problem
\[\begin{array}{l}
\text{To find $\phi\in H^1_{\mathrm{loc}}(\mathbb R^N)$ such that }
\\
\ds\int \nabla \phi\cdot\nabla(g-\phi)\ud x\geq 0,
\\
\text{holds for any $g\in H^1_{\mathrm{loc}}(\mathbb R^N)$, 
with $g-\phi $ compactly supported  and  $\phi\geq \psi$ q.~e.}
\end{array}\]
where $\psi(x) \ddef C_{*}-\Phi_\ext(x)$. 
\end{prop}

We then define the coincidence set
$$ 
\cak \ddef \{ \Phi_{\rm e} = 0 \} = \{ \Gamma \star n_{\rm e} = C_{*} - \Phi_\ext \} 
$$
and claim that $\cak$ is compact. Indeed, as $ |x| \to +\infty $, we
have $ - \Gamma \star n_{\rm e}(x) \sim - \mathfrak{m} \Gamma (x) $
and $ \Phi_\ext(x) + \mathfrak{m} \Gamma (x) \gg 1 $ whatever is the
dimension $N$ by h1)-h2), thus $ \cak = \{ \Gamma \star n_{\rm e} =
C_{*} - \Phi_\ext \} $ is bounded. Moreover, by \eqref{Robin}, the set
$ {\rm Supp} ( n_{\rm e} ) \setminus \cak $ has zero capacity. We give
some examples in section \ref{sec:hypo} below where $ {\rm Supp} (
n_{\rm e} ) \subsetneq \cak $, due to the presence of points or
regions where $\Delta \Phi_{\rm ext} $ vanishes. These points are
precisely defined in \cite[Section 3.6]{HanMak} and called 'shallow
points' and it is shown in this paper (see Proposition 3.12 there)
that it is possible to pass from $ {\rm Supp} ( n_{\rm e} ) $ to $
\cak $ by simply adding these 'shallow points'. This fact is
illustrated in section \ref{sec:hypo} below.

For a general potential $\Phi_{\mathrm {ext}}$, the variational
viewpoint and the theory of the obstacle problem provide a definition
for the equilibrium distribution $n_{\rm e}$, the domain $\cak$ (which
is not always the support of $ n_{\rm e} $) and the potential $
\Phi_{\rm e} $.  The regularity of $\Phi_{\mathrm {ext}} $ is not
``transferred'' to the solution $ \Phi_{\rm e} $ or $ \Gamma \star
n_{\rm e} $ beyond $ C^{1,1} $ regularity (see \cite{Frehse},
\cite{LC}) since the Laplacian of these functions is discontinuous.
In addition, the topology of $ \cak $ is difficult to analyse in
general: $\cak $ may have empty interior or may exhibit cusps.
Hence, these regularity issues for both $ \cak $ and $n_{\rm e} $ need to be 
discussed individually.  Let us then list the properties, which very 
likely are far from optimal, that we need to deal with the asymptotic 
regime: there exists $s>1+N/2$ such that
 \begin{itemize}
 \item[H1)] $\cak$ has a non empty interior $\Omega$ and
 $\partial\Omega$ is of class $C^{1}$.
 \item[H2)] $ \Phi_{\rm ext} \in C^{s+3}(\mathbb R^N) $, 
 $ \Delta \Phi_{\rm ext} $ is bounded away from zero on $ \cak $. 
\end{itemize}

The $C^1 $ regularity assumption H1) on $\partial \Omega $ excludes
the presence of cusps in $\cak$. Let us point out two regularity
results derived from the obstacle problem theory.

\begin{prop} 
  Let $ \Phi_{\mathrm {ext}} $ be a potential satisfying $h1)$ and
  $h2)$ and consider $ n_{\rm e} $ the minimizer of \eqref{eq:var}
  given by Theorem \ref{Minimizor} and let $\cak \ddef \{ \Phi_{\rm
    e}=0 \} $.
\begin{itemize}
\item[i)] \cite{Kinder-Niren_77} Assume that $H1)$ and $H2)$ are
  satisfied. Then $ \bar{\Omega} = \cak $ and the boundary $ \partial
  \Omega $ is $ C^{s+1} $.
\item[ii)] \cite{Frehse}, \cite{LC}, \cite{HanMak} Assume that $
  \Phi_{\rm ext} \in C^{1,1} (\R^N ) $.  Then, $ \Gamma \star n_{\rm
    e} \in C^{1,1} $ and $ n_{\rm e} = {\bf 1}_{\Omega} ( \Delta
  \Phi_{\rm ext} ) $ as measures.
 \end{itemize}
\end{prop}

\begin{remark} 
  A consequence of ii) is that, in \eqref{Robin}, we may replace
  'quasi everywhere' by 'everywhere' since all the functions involved
  are continuous.  Note that $H2)$ then implies that $n_{\rm e}$ is
  $C^{s+1}$ and bounded from below on $\cak$.
\end{remark}

\begin{remark} 
  Under the low regularity assumption $ \Phi_{\rm ext} \in C^{1} (\R^N
  ) $ and when $N \ge 2 $, it follows from \cite[Theorem 2]{LC} that $
  \Gamma \star n_{\rm e} \in C^1 $.  This prevents the singular part
  of the measure $ n_{\rm e} $ from being a Dirac mass or a finite sum
  of Dirac masses, since the fundamental solution $ \Gamma $ is
  unbounded (if $N \ge 2 $) near the origin. This however may happen
  in dimension 1 (see the examples in section \ref{sec:hypo} below)
  and in these cases, the relation $ n_{\rm e} = {\bf 1}_{\Omega} (
  \Delta \Phi_{\rm ext} ) $ (as measures) might not be true.
\end{remark}

\noindent {\it Proof.} 
For the first point, notice that $\bar{\Omega} = \cak $ by H1).  In
addition, we also have from \cite[Theorem 1]{Kinder-Niren_77} or
\cite[Chapter 2, Theorem 1.1]{Friedman}, since we assume that $ \Delta
\Phi_{\rm ext} \in C^{s+1}_{\rm loc}(\mathbb R^N) $ and does not
vanish on $ \cak $, that the boundary $ \partial \Omega $ is
automatically of class $ C^{s+1} $ (and even $C^{s+1+ \beta} $ for any
$ \beta \in (0,1 ) $).  If $\Phi_{\rm ext} $ is analytic, then
$ \partial \Omega $ is also analytic.

For the second statement, we first invoke the regularity result of
\cite{Frehse} (see also \cite{LC}, \cite{Friedman}) saying that since
$ \Phi_{\rm ext} \in C^{s+3} (\mathbb R^N) \subset C^{1,1} (\mathbb
R^N) $, then $ \Gamma \star n_{\rm e} $ belongs to $ C^{1,1} (\mathbb
R^N) $. Consequently, in the distributional sense, the compactly
supported measure $ n_{\rm e} = - \Delta ( \Gamma \star n_{\rm e} ) $
belongs to $ L^\infty ( \mathbb R^N) $. Since $ \Gamma \star n_{\rm e}
= C_{*} - \Phi_\ext $ in $ \Omega $ (q. e., hence everywhere by
continuity of the functions), we infer that $ n_{\rm e} = \Delta
\Phi_\ext $ in $ \Omega $.  If we make the assumption H1), $\partial
\Omega $ is of class $ C^1 $ and we then deduce that $ n_{\rm e} = (
\Delta \Phi_{\rm ext} ) {\bf 1}_{\Omega} $ as a measure.  If
assumption H1) is not satisfied, then, as noticed in \cite[Theorem
3.10]{HanMak}, it follows from \cite[Chapter 2, Lemma
A.4]{Kinder-Stampa_88} that $ n_{\rm e} = \Delta \Phi_\ext $ holds
almost everywhere in $ \cak $, which concludes. \qed

With these assumptions H1) and H2), we can establish the following statement, 
where we point out that the limit problem is the Lake Equation \eqref{lake} 
instead of the mere incompressible Euler system, since now the equilibrium 
distribution $n_{\rm e}$ is inhomogeneous. We obviously need a smooth enough 
solution to the Lake Equation \eqref{lake}: we may refer to the works 
\cite{Lever-Oli-Titi_PhysD}, \cite{Oliver_97} (when the domain $\Omega$ is 
simply connected) and \cite{Lever-Oli-Titi_Indi} (without simple connectedness 
assumption on the domain $\Omega$), which rely on a vorticity formulation 
{\it \`a la} Yudovitch and are then restricted to the dimension $ N=2 $. 
We provide in the appendix (see Theorem \ref{Lakeexiste}) a well-posedness 
result analogous to Theorem \ref{Eulerexiste} valid in any dimension and 
without simple connectedness assumption on $ \Omega $.

\begin{theorem}
\label{thresultlake} 
Let $ \Phi_{\mathrm {ext}} $ be a potential satisfying $h1)$ and $h2)$
and consider $ n_{\rm e} $ the minimizer of \eqref{eq:var} given by
Theorem \ref{Minimizor} and let $\cak \ddef \{ \Phi_{\rm e}=0 \} $.
Assume in addition that $H1)$ and $H2)$ are satisfied. Let $
V^{\mathrm{init}} \in H^s(\Omega)$ satisfy $\nabla_x\cdot (n_{\rm
  e}V^{\mathrm{init}})=0$ in $ \Omega$ and the no flux condition
\eqref{CLEuler}.  Denote by $V$ the solution on $ [ 0 ,T ] $ to the
Lake Equation \eqref{lake}, with the no flux condition \eqref{CLEuler}
and initial condition $V^{\mathrm{init}}$, given in
Theorem~\ref{Lakeexiste} and consider $\mathscr V^{\mathrm{init}}
$ a smooth extension of $V^{\mathrm{init}}
$ satisfying the following conditions, where $R>0 $ is such that $
\Omega \subset B(0,R)$,
\[
\mathscr V\Big|_{\Omega}=V,\qquad
\mathscr V\Big|_{\mathbb R^N\setminus B(0,2R)}=0,
\qquad
\mathscr V(t,x)\cdot \nu(x)\Big|_{\partial\Omega}=0 .
\]
Let $f_\eps^{\mathrm{init}}:\mathbb R^N\times
\mathbb R^N\rightarrow [0,\infty)$ be a sequence of integrable
functions that satisfy
\eqref{query}. 
Then, the associated solution $f_\eps$ of the Vlasov--Poisson equation 
\eqref{VP1}--\eqref{VP2} satisfies, 
as $ \eps\rightarrow 0 $,
\begin{itemize}
\item[i)] $\rho_\eps$ converges to $n_{\mathrm e}$ in 
$C^0(0,T;\mathscr M^1(\mathbb R^N)-\text{weak}-\star)$;
\item[ii)] $ \mathscr H_{\mathscr V,\eps}$ converges to $0 $ uniformly
  on $[ 0, T ]$;
\item[iii)] $J_\eps$ converges to $J$ in $\mathscr M^1([0,T]\times
  \mathbb R^N)$ weakly-$\star $, the limit $J$ lies in
  $L^\infty(0,T;L^2(\mathbb R^N))$ and satisfies $J\big|_{[0,T]\times
    \Omega}=V$, $\nabla_x\cdot J=0$ and $J\cdot \nu(x)\big|_{\partial\Omega}=0$.
\end{itemize}
\end{theorem}

The existence of a smooth extension $ \mathscr V $ of $ V $ follows
from~\cite[Chapter I: Theorem 2.1 p. 17 \& Theorem 8.1 p. 42]{LM}.

For convex potentials $ \Phi_{\rm ext} $, the only situation where we
have been able to check the hypotheses H1) and H2) (except the
quadratic potentials for which $ \Delta \Phi_{\rm ext} $ is constant)
is the case of the space dimension $N=1$ (see Proposition
\ref{1Dconvex}) and the case of a radial potential (see Proposition
\ref{radi} below).

\subsection{About hypothesis $H1)$ for convex potentials $ \Phi_{\rm ext} $}
\label{sec:hypo}

For the problem we have in mind, it is natural to assume that the confining 
potential $ \Phi_{\rm ext} $ is smooth and convex. In this case, one may think 
that the coincidence set $ \cak $ or $ {\rm Supp} ( n_{\rm e} ) $ is convex. 
We have not been able to find such a result in the literature for a general convex, 
coercive and smooth enough confining potential $\Phi_{\rm ext} $. 
Actually, the obstacle problem is, in most cases, set on a bounded convex domain $G$ 
with suitable boundary conditions instead of the whole space $ \mathbb R^N $. 

For the obstacle problem in bounded convex domains $G$, 
we can find a convexity result for the coincidence set $\cak $ in 
\cite[Theorem 6.1]{FriedPhil} in the specific assumptions that 
$ \Delta \Phi_{\rm ext} $ is constant and with the boundary condition 
$ \Gamma \star n_{\rm e} = 1 + \psi = 1 + C_* - \Phi_{\rm ext} $ on $ \p G $. 
Just after \cite[Theorem 6.1]{FriedPhil}, an example is given (in a bounded convex 
domain $G$) showing that the assumption $ \Phi_{\rm ext} $ smooth and strictly 
convex (and $ \Gamma \star n_{\rm e} > \psi $ on $ \p G $) is not sufficient to 
guarantee that $\cak $ is convex. Roughly speaking, $ \Delta \Phi_{\rm ext} $ is 
constant for quadratic potentials.

Turning back to the obstacle problem in the whole space $\mathbb R^N $, the only 
convexity result we are aware of is \cite[Corollary 7]{LC}, 
which corresponds to the case where $ \Delta \Phi_{\rm ext } $ is constant. 
Extending this result to space depending functions $ \Delta \Phi_{\rm ext } $ 
is a delicate issue (see however \cite[Chapter 2, section 3]{Friedman}, which 
is not sufficient for our situation).

In the one dimensional case and for a convex potential $ \Phi_{\rm ext}$, 
there is a simple characterization of $ \cak $, as explicited in the following Proposition.

\begin{prop}[The one dimensional case with a convex potential] 
\label{1Dconvex}
Assume that $ N =1 $ and that $ \Phi_{\rm ext} : \R \to \R $ is of class $ C^1 $, 
piecewise $ C^2 $, nonnegative, convex ({\it i.e.} $ \Phi_{\rm ext}' $ is nondecreasing) 
and that $ \Phi_{\rm ext} (x) - \mathfrak{m} |x| / 2 \gg 1 $ for $ |x| \gg 1 $ 
(so that h1) and h2) are satisfied). 
We denote by $ \p \Phi'_{\rm ext} $ the piecewise continuous function 
associated with the second order derivative of $ \Phi_{\rm ext} $. 
Then, the minimizer $ n_{\rm e}$ for \eqref{eq:var} is given by
\begin{equation}
 \label{Mini1D}
  n_{\rm e} = ( \p \Phi'_{\rm ext} ) \Big|_{ ] a_-, a_+ [ } 
\end{equation}
where $ a_+ $ and $ a_- $ are defined by the equations
\begin{equation}
 \label{Eq:bord}
   \frac{\mathfrak{m}}{2} = \Phi_{\rm ext}'(a_+) 
   \quad \quad \quad {\it and} \quad  \quad \quad 
   - \frac{\mathfrak{m}}{2} = \Phi_{\rm ext}'(a_-) .
\end{equation}
Furthermore, $ {\rm Supp} ( n_{\rm e} ) = {\rm Supp} ( \p \Phi'_{\rm ext} ) \cap [a_- , a_+ ] $ 
and $ \{ \Phi_{\rm e} = 0 \} = [a_- , a_+ ] $. 
In addition, the potential $ \Phi_{\rm e} $ is convex.
\end{prop}

\noindent {\it Proof.} As a first observation, notice that \eqref{Eq:bord} 
has at least one (possibly non unique) solution since $ \Phi_{\rm ext}' $ is 
continuous, nondecreasing and tends to $ \geq \mathfrak{m}/2 $ (resp. 
$ \leq -\mathfrak{m}/2 $) in view of our hypothesis. If the limit 
at $ +\infty $  is $ \mathfrak{m}/2 $, 
it follows from the convexity of $ \Phi_{\rm ext} $ that 
$ \Phi_{\rm ext}(x) - \Phi_{\rm ext}(y) \leq \mathfrak{m}/2 $. 

cici

(resp. $ - \infty$) (resp. $ - \mathfrak{m}/2 $)
Since $ \Phi_{\rm ext}' $ is nondecreasing, and if $ b_+ > a_+ $ also 
solves $ \mathfrak{m}/2 = \Phi_{\rm ext}'(b_+) $, this implies 
that, on $[ a_+, b_+] $, $ \Phi_{\rm ext}' \equiv \mathfrak{m}/2 $, thus 
$ \p \Phi'_{\rm ext} \equiv 0 $ and this does not change $ n_{\rm e} $.

Let us use the characterization (iii) in Theorem \ref{Minimizor} and look for 
the measure $ n_{\rm e} $ under the form 
$ n_{\rm e} = ( \p \Phi'_{\rm ext} )\Big|_{ ] a_-, a_+ [ } $, which is 
piecewise continuous. This function $ n_{\rm e} $ satisfies the mass constraint 
if and only if
\begin{equation}
\label{alamasse}
\mathfrak{m} = \int_{a_-}^{ a_+ } \p \Phi'_{\rm ext} \, \ud x
= \Phi_{\rm ext}'(a_+) - \Phi_{\rm ext}'(a_-) .
\end{equation}
Now, let us compute $ \Gamma \star n_{\rm e} + \Phi_{\rm ext}$ in $ [ a_-
, a_+ ] $ and investigate under which condition this function is
constant (in $ [ a_- , a_+ ] $). Elementary computations give, for $ a_- \le x \le a_+ $:
\begin{align*}
\Gamma \star n_{\rm e} (x) = &\, 
- \frac12 \int_{a_-}^{ a_+ } | y-x | (\p \Phi'_{\rm ext} ) (y) \, \ud y \\ 
= &\, - \frac12 \Phi_{\rm ext}' (a_+) ( a_+ - x) + \frac12 \Phi_{\rm ext}' (a_-) ( x - a_- ) 
+\frac12 \int_{a_-}^{a_+} {\rm sgn} ( y-x ) \Phi_{\rm ext}'(y) \, \ud y
  \\
  = &\, - \frac12 \Phi_{\rm ext}' (a_+) ( a_+ - x) 
  + \frac12 \Phi_{\rm ext}' (a_-) ( x - a_- ) + \frac12 \Phi_{\rm ext}(a_+) 
  + \frac12 \Phi_{\rm ext}(a_-) - \Phi_{\rm ext}(x) .
\end{align*}
As a consequence, $ \Gamma \star n_{\rm e} + \Phi_{\rm ext} $ is constant
in $ [ a_-, a_+] $ if and only if $ \Phi_{\rm ext}'(a_+) + \Phi_{\rm ext}'(a_-) = 0 $. 
Combining this with the mass constraint 
$ \Phi_{\rm ext}'(a_+) - \Phi_{\rm ext}'(a_-) = \mathfrak{m} $ yields the relation 
\eqref{Eq:bord}. It then follows that, on $[ a_- , a_+ ] $, 
\begin{align*}
 \Gamma \star n_{\rm e} + \Phi_{\rm ext} = 
 C_* \ddef & \, \frac12 \left(  \Phi_{\rm ext}(a_+) + \Phi_{\rm ext}(a_-) 
 - a_+ \Phi_{\rm ext}' (a_+) -  a_- \Phi_{\rm ext}' (a_-)  \right ) 
 \\
 = & \, \frac12 \left( \Phi_{\rm ext}(a_+) + \Phi_{\rm ext}(a_-) \right ) 
 - \frac{\mathfrak{m}}{4} ( a_+ - a_- ) . 
\end{align*}
It remains to check that $  \Gamma \star n_{\rm e} + \Phi_{\rm ext} \ge C_* $ 
in $ \R $. To see this, note that $ \Phi_{\rm e} \ddef \Gamma \star n_{\rm e} + \Phi_{\rm ext} - C_* $ 
is convex since its (distributional) second order derivative is equal to the piecewise continuous 
function $ \p_x \Phi_{\rm ext}' {\bf 1}_{\R \setminus [ a_- , a_+ ] } $, and $ \Phi_{\rm e} \equiv 0 $ 
on $ [ a_- , a_+ ] $, hence is $ \ge 0 $ everywhere. This finishes the proof. \qed

Let us give some examples illustrating Proposition \ref{1Dconvex}.\\

\noindent {\bf Example 1 (1D):} If $ \Phi_{\rm ext} $ is of class $ C^2 $ 
and $ \Phi_{\rm ext}'' $ is positive on $ \R $, then 
$ n_{\rm e}(x) = \Phi_{\rm ext}''(x) {\bf 1}_{[ a_-, a_+] } (x) $ and is 
absolutely continuous with respect to the Lebesgue measure. 
We then have $ {\rm Supp}(n_{\rm e}) = [a_-, a_+ ] $.\\

\begin{figure}
\begin{center}
\includegraphics[width=0.8\linewidth,height=0.5\linewidth]{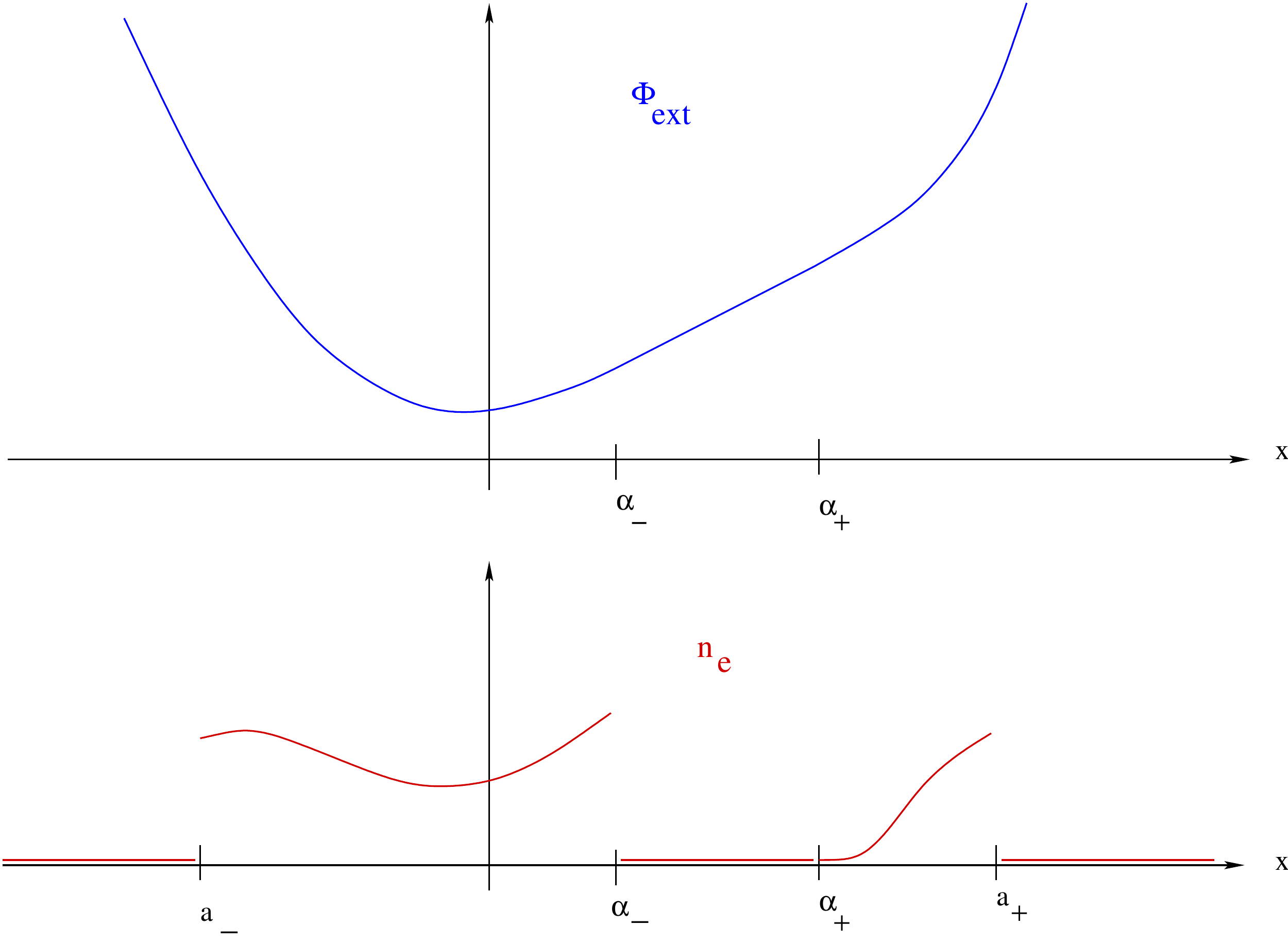}
\end{center}
\caption{The potential $ \Phi_{\rm ext} $ and the corresponding measure 
$ n_{\rm e} $ for example 2}
\label{pot2}
\end{figure} 

\noindent {\bf Example 2 (1D):} the potential $ \Phi_{\rm ext} $ is $ C^1 $, piecewise $C^2$, 
but is affine on the interval $ [ \alpha_- , \alpha_+ ] $ (hence it is not strictly convex), 
where its slope belongs to $ ] - \mathfrak{m}/2 , + \mathfrak{m}/2 [ $ 
(see figure \ref{pot2}). In addition, the second order derivative $ \Phi_{\rm ext}'' $ 
is discontinuous at $ \alpha_- $ and continuous at $ \alpha_+ $ and $ \p \Phi'_{\rm ext} $ is 
positive except on $ [ \alpha_- , \alpha_+ ] $. In this case, we may still define 
$ a_\pm $ as the unique solutions to $ \Phi_{\rm ext}' (a_\pm ) = \pm \mathfrak{m}/2 $, and 
we have $ {\rm Supp}(n_{\rm e}) = [a_-, \alpha_- ] \cup [\alpha_+, a_+ ] \subsetneq 
[a_- , a_+ ] = \{ \Phi_{\rm e} = 0 \} $ and this is then a disconnected set. 
If the slope in the region $ [ \alpha_- , \alpha_+ ] $ where $ \Phi_{\rm ext} $ 
is affine does not belong to $ ] - \mathfrak{m}/2 , + \mathfrak{m}/2 [ $, then 
the support of $ n_{\rm e} $ is an interval as in Example 1.\\

Example 1 fits the hypotheses of Theorem \ref{thresultlake}, but not Example 2 
since $ n_{\rm e} $ is not bounded away from zero (near $ \alpha_+ $). 
In particular, for Example 2, we have to face new difficulties in solving the Cauchy 
problem (see Theorem \ref{Lakeexiste}) for the Lake Equation \eqref{lake}. 
If in the one dimensional situation one can easily check that the support of $n_{\rm e} $ 
(instead of $ \cak $) is smooth, in a similar higher dimensional case, the 
regularity of Supp $(n_{\rm e} ) $ is certainly not easy to analyse since 
we can not rely on the results in \cite[Theorem 1]{Kinder-Niren_77} or 
\cite[Chapter 2, Theorem 1.1]{Friedman}. All these issues motivate hypothesis H2).

Let us give now examples which do not fit the regularity hypotheses required in 
Proposition \ref{1Dconvex}. These expressions are justified through the characterization 
(iii) in Theorem \ref{Minimizor} and simple computation of $ \Gamma \star n_{\rm e} $.\\

\noindent {\bf Example 3 (1D):} Take the potential $\Phi_{\rm ext} (x) = |x| $. 
Then, hypothesis h2) exactly means $ \mathfrak{m} < 1 $. In that case, 
we have $ n_{\rm e} = \mathfrak{m} \delta_0 $ and 
$ {\rm Supp}(n_{\rm e}) = \{ 0 \} = \{ \Phi_{\rm e} = 0 \} $.\\

\noindent {\bf Example 4 (1D):} Take two reals $ a < b $ and a convex potential 
$ \Phi_{\rm ext} $ which is affine on $ ] - \infty , a] $, on $ [a,b] $ and on $ [ b ,+\infty [ $. 
Assume also that h2) is satisfied, that is 
$ \mathfrak{m} < \min (\Phi'_{\rm ext} ( +\infty ) , - \Phi'_{\rm ext} ( -\infty )) $. 
Then, $ \ds n_{\rm e} = 
 \min \left( \frac12\left( \mathfrak{m} 
 + \frac{\Phi_{\rm ext}(b) - \Phi_{\rm ext} (a) }{b-a}\right)_+ , \mathfrak{m}\right)\delta_a
+ \min \left( \frac12 \left( \mathfrak{m} - \frac{\Phi_{\rm ext}(b) 
- \Phi_{\rm ext} (a) }{b-a} \right)_+ , \mathfrak{m}\right) \delta_b $. 
As a consequence:\\
- if $ \ds - \mathfrak{m} < \frac{\Phi_{\rm ext}(b) - \Phi_{\rm ext} (a) }{b-a} < \mathfrak{m} $, 
then 
$ {\rm Supp}(n_{\rm e}) = \{ a,b \} $ and $ \{ \Phi_{\rm e} = 0 \} = [a,b ] $;\\
- if $ \ds \frac{\Phi_{\rm ext}(b) - \Phi_{\rm ext} (a) }{b-a} \leq - \mathfrak{m} $, 
then $ n_{\rm e} = \mathfrak{m} \delta_b $ 
and $ {\rm Supp}(n_{\rm e}) = \{ b \} = \{ \Phi_{\rm e} = 0 \} $;\\
- if $ \ds \frac{\Phi_{\rm ext}(b) - \Phi_{\rm ext} (a) }{b-a} \geq \mathfrak{m} $, 
then $ n_{\rm e} = \mathfrak{m} \delta_a $ and 
$ {\rm Supp}(n_{\rm e}) = \{ a \} = \{ \Phi_{\rm e} = 0 \} $.\\

\noindent {\bf Example 5 (1D):} Consider the potential $\Phi_{\rm ext} (x) = |x| + x^2 /2 + \max( x-1, 0) $: \\
- if $ \mathfrak{m} \leq 2 $, then $ \ds n_{\rm e} = \mathfrak{m} \delta_0 $, 
$ \Gamma \star n_{\rm e} (x) = - \mathfrak{m} |x| /2 $ and 
$ {\rm Supp}(n_{\rm e}) = \{ 0 \} = \{ \Phi_{\rm e} = 0 \} $;\\
- if $ 2 \le \mathfrak{m} \le 4 $, then 
$ \ds n_{\rm e} = 2 \delta_0 + {\bf 1}_{[ - \mathfrak{m}/2 +1,\mathfrak{m}/2 -1]}$ 
and $ {\rm Supp}(n_{\rm e}) = [ - \mathfrak{m}/2 +1,\mathfrak{m}/2 -1] = \{ \Phi_{\rm e} = 0 \} $;\\
- if $ 4 \le \mathfrak{m} \le 6 $, then 
$ \ds n_{\rm e} = 2 \delta_0 + (\mathfrak{m}/2 - 2 ) \delta_1 + {\bf 1}_{[ - \mathfrak{m}/2 + 1,1]} $ and 
$ {\rm Supp}(n_{\rm e}) = [ - \mathfrak{m}/2 +1,1] = \{ \Phi_{\rm e} = 0 \} $;\\
- if $ \mathfrak{m} \ge 6 $, then 
$ \ds n_{\rm e} = 2 \delta_0 + \delta_1 + {\bf 1}_{[ - \mathfrak{m}/2 + 1, \mathfrak{m}/2 - 2] } $ 
and $ {\rm Supp}(n_{\rm e}) = [ - \mathfrak{m}/2 + 1, \mathfrak{m}/2 - 2] = \{ \Phi_{\rm e} = 0 \} $. \\


Examples 3, 4 and 5 show that the single convexity hypothesis on $ \Phi_{\rm ext} $ 
does not guarantee that $ n_{\rm e} $ is a restriction of the nonnegative measure 
$ \p_x^2 \Phi_{\rm ext} $ (in the distributional sense). It appears in these 
examples that $ n_{\rm e} $ is nondecreasing with respect to the mass $ \mathfrak{m} $, 
and thus that we always have $ n_{\rm e} \le \p_x^2 \Phi_{\rm ext} $ in the distributional 
sense. It is an open problem to determine whether this holds true in higher dimensions. 
Here again, these issues motivate the regularity assumptions on $ \Phi_{\rm ext} $ in H2).

The other situation where we may verify hypothesis H1) is the radial case 
(see \cite[Corollary 1.4]{ChaGoZi} for a related result in dimension $N\ge 3 $ for 
$ C^2 $ potentials $ \Phi_{\rm ext} $). 
Let $ \vp_{\rm ext} : \R_+ \rightarrow \R_+ $ be a nondecreasing function of 
class $C^1$ and piecewise $ C^2 $. Consider now the potential $ \Phi_{\rm ext} : \R^N \to \R $ 
given by $ \Phi_{\rm ext} (x) = \vp_{\rm ext} ( | x | ) $. It is then clear 
that $ \vp_{\rm ext} $ is convex if and only if $ \Phi_{\rm ext} $ is convex.

\begin{prop}[The radial case with a convex potential]
\label{radi}
Assume that $ N \ge 2 $ and that $ \Phi_{\rm ext} : \R^N \to \R $ is as above. 
Then, the minimizer $ n_{\rm e} $ for \eqref{eq:var} is given by
\begin{equation}
 \label{Miniradi}
  n_{\rm e} (x) = {\bf 1}_{ B (0, R ) } (x) \Delta \Phi_{\rm ext} (x),
\end{equation}
where $ R $ is defined by the equation
\begin{equation}
 \label{Eq:bordradi}
   \mathfrak{m} = \int_{B(0,R)} \Delta \Phi_{\rm ext}(x) \, \ud x 
   \quad \quad {\it or,\ equivalently, } \quad \quad 
   N | B(0,1) | R^{N-1} \vp_{\rm ext}' (R) = \mathfrak{m} .
\end{equation}
Furthermore, $ {\rm Supp} ( n_{\rm e} ) = \bar{B}(0, R ) \setminus B(0, R_{\rm min} ) $, 
where $ R_{\rm min} \ddef \max \{ \vp_{\rm ext}' = 0 \} \le R $. In addition, the potential 
$ \Phi_{\rm e} $ is convex.
\end{prop}

\noindent {\it Proof.} The existence of $R$ is clear. We may have non uniqueness 
only in the case where $ \Phi_{\rm ext } $ is constant on a ball $ B(0,R_0) $ (of positive 
radius), since $ \vp_{\rm ext} '$ is nondecreasing. The potential $\Phi_{\rm e}  $ may be searched for
under the form of a radial function, and we find the expressions
\begin{equation*}
  \Phi_{\rm e}  (x) = ( \vp_{\rm ext} (R ) - \vp_{\rm ext} (|x| ) 
  + \vp_{\rm ext}'(R) \Gamma( R) ) \mathbf 1_{B(0,R)} 
  + \vp_{\rm ext}'(R) \Gamma( x) \mathbf 1_{\R^N \setminus B(0,R)} ,
\end{equation*}
where $ \Gamma( R) $ stands for $ \Gamma (y) $ for any $y \in \partial B(0,R) $.\qed


Let us give some examples illustrating Proposition \ref{radi}.\\

\noindent {\bf Example 1 (radial):} If $ \vp_{\rm ext} $ is of class $ C^2 $ 
and $ \vp_{\rm ext}'' $ is positive on $ \R_+ $, then 
$ n_{\rm e}(x) = {\bf 1}_{ B(0,R) } (x) \Delta \Phi_{\rm ext} (x) $ and is 
absolutely continuous with respect to the Lebesgue measure. 
We then have $ {\rm Supp}(n_{\rm e}) = \bar{B}(0,R) $.\\

\noindent {\bf Example 2 (radial):} The potential $ \vp_{\rm ext} $ is $ C^1 $, piecewise $C^2$, 
but is constant on the interval $ [ 0, R_0 ] $ (hence it is not strictly convex). 
It does not matter whether the second order derivative of $ \vp_{\rm ext} $ is continuous or not 
at $R_0 $. We define $ R \ge R_0 > 0 $ by the relation 
$ \mathfrak{m} = \int_{B(0,R)} \Delta \Phi_{\rm ext} \, \ud x $, 
or, equivalently, $ N | B(0,1) | R^{N-1} \vp_{\rm ext}' (R) = \mathfrak{m} $. 
Then, $ n_{\rm e} = {\bf 1}_{ B(0,R) \setminus B(0,R_0)} \Delta \Phi_{\rm ext} $, 
$ {\rm Supp}(n_{\rm e}) = \bar B (0,R) \setminus B(0,R_0) \subsetneq 
\bar B (0,R) = \{ \Phi_{\rm e} = 0 \} $ and this set is then neither starshaped 
nor simply connected. Here again, if $ \vp_{\rm ext} \in C^2 $, this potential does 
not fit hypothesis H2) since $ \Delta \Phi_{\rm ext} $ is not bounded 
away from 0 near $ R_0 $.\\

Let us give now examples which do not fit the regularity hypotheses required in 
Proposition \ref{1Dconvex}. These expressions are justified through the characterization 
(iii) in Theorem \ref{Minimizor} and simple computation of $ \Gamma \star n_{\rm e} $.\\

\noindent {\bf Example 3 (radial):} Take the potential $ \vp_{\rm ext} (r) = r $, that is 
$ \Phi_{\rm ext} (x) = \lvert x \rvert $. Then, $ \Delta \Phi_{\rm ext} = (N -1) /r > 0 $, 
$ n_{\rm e} = ( N -1 ) |x|^{-1} {\bf 1}_{B(0,R)} $, with $ N | B(0,1) | R^{N-1} = \mathfrak{m} $, 
and $ {\rm Supp}(n_{\rm e}) = \bar{B}(0,R) = \{ \Phi_{\rm e} = 0 \} $.\\

\noindent {\bf Example 5 (radial):} Consider the potential $ \vp_{\rm ext} (r) = r + \max( r-1, 0) $: \\
- if $ \mathfrak{m} \leq N | B(0,1) | $, then $ \ds n_{\rm e} = (N-1) |x|^{-1} {\bf 1}_{B(0,R)} $, 
with $ R = ( \mathfrak{m} / N | B(0,1) | )^{1/N-1} $ and 
$ {\rm Supp}( n_{\rm e}) = \bar B(0,R) = \{ \Phi_{\rm e} = 0 \} $;\\
- if $ N | B(0,1) | \le \mathfrak{m} \le 2 N | B(0,1) | $, then 
$ \ds n_{\rm e} = (N-1) |x|^{-1} {\bf 1}_{B(0,1)} 
+ ( \mathfrak{m} - N | B(0,1) | ) \delta_{ \partial B(0,1)} $ 
and $ {\rm Supp}(n_{\rm e}) = \bar B (0,1) = \{ \Phi_{\rm e} = 0 \} $;\\
- if $ \mathfrak{m} \ge 2 N | B(0,1) | $, then 
$ \ds n_{\rm e} = (N-1) |x|^{-1} {\bf 1}_{B(0,1)} + N | B(0,1) | \delta_{ \partial B(0,1)} 
+ 2 (N-1) |x|^{-1} {\bf 1}_{ B(0,R) \setminus B(0,1)} $, 
where $ R \ge 1 $ is such that $ 2 N | B(0,1) | ( R^{N-1} - 1 ) + N | B(0,1) | = \mathfrak{m} $, 
and $ {\rm Supp}(n_{\rm e}) = \bar B (0,R) = \{ \Phi_{\rm e} = 0 \} $.\\

Since we assume $ \vp_{\rm ext} $ convex and with $0$ as a minimum point, it follows that 
$ \vp_{\rm ext} $ has a right-derivative at $0$, hence the singularity in $ 1 / |x| $ at the 
origin for $n_{\rm e} $ is the worst we can have. The radial Example 5 also shows that 
we may have Dirac masses on a sphere (of positive radius).

Our next results guarantees that $ \mathcal{K} $ has non empty interior when the 
confining potential $ \Phi_{\rm ext} $ is $C^1 $ and convex.

\begin{prop}
\label{internonvide} 
We assume that $0$ is a minimum point of $ \Phi_{\rm ext} $ and that the potential 
$ \Phi_{\rm ext} $ is of class $ C^1 $ and convex. 
Then, there exists $ r_0 > 0 $ such that $ B_{r_0} (0) \subset \cak $. 
In particular, $ \cak $ has non empty interior.
\end{prop}

\noindent {\it Proof.} We follow the argument of \cite[Chapter 5, Theorem 6.2]{Kinder-Stampa_88}, 
where we work on $ h \ddef \Gamma \star n_{\rm e} $ and shall use that it is a solution to 
the obstacle problem given in Proposition \ref{sautobstacle} with the obstacle 
$ \psi = C_* - \Phi_{\rm ext} $. 

We first consider the case $ N \ge 3 $ and notice that Supp$(n_{\rm e})$ has a positive capacity: 
we fix some $ a \in$ Supp$(n_{\rm e})$ such that $ C_* = h(a) + \Phi_{\rm ext} (a) $ (see \eqref{Robin}). 
Now, since $N \ge 3 $, we observe that (with $c_N > 0 $) 
$ h(a) = \Gamma \star n_{\rm e} (a) = c_N |\cdot|^{2 - N} \star n_{\rm e} > 0 $ and that 
$0$ is actually a global minimum point of $ \Phi_{\rm ext} $, thus 
$ C_* > \Phi_{\rm ext} (a) \ge \Phi_{\rm ext} (0) $ and it follows that 
$ \psi (0 ) = C_* - \Phi_{\rm ext} (0) > 0 $. On the other hand, 
$ h(x) \sim \mathfrak{m} \Gamma(x) $ tends to $0 < \psi (0) $ at infinity, 
thus there exists an $ R_0 > 0 $ such that $ h (x) \le \psi(0) /2 $ when $ |x| \ge R_0 $. 
For $x_0$ that will be close to $0$, we let 
$v (x) \ddef \psi(x_0) + (x-x_0)\cdot \nabla \psi (x_0 ) $ be the affine tangent 
to $ \psi $ at $x_0 $. Since $ \psi $ is concave ($\Phi_{\rm ext} $ is convex), we have 
$ \psi \le v $ in $\mathbb R^N $. Furthermore, if $ x_0$ is sufficiently close to $0$ 
(depending on $R_0$), then $ \nabla \psi (x_0 ) $ is small (since $ \psi $ is $C^1$ 
and achieves a minimum at $0$) and thus $ v > \psi(0) /2 > 0 $ 
on $ \partial B (0, R_0 ) $. Since $ \Delta v \equiv 0 $, we may now apply 
\cite[Chapter 4, Theorem 8.3]{Kinder-Stampa_88} to infer $ h \le v $ in $ B (0, R_0 ) $ 
(this is a maximum type principle proved using the comparison function 
$ g \ddef \min( h , v) {\bf 1}_{B(0,R_0)} + h {\bf 1}_{\mathbb R^N \setminus B(0,R_0)} $ 
in the formulation of the obstacle problem given in Proposition \ref{sautobstacle}). In particular, 
$ \psi (x_0 ) \le h (x_0 ) \le v (x_0 ) = \psi(x_0 ) $, which means that, as wished, 
$ x_0 \in \cak $. 

Let us now turn to the dimensions $N=2$ and $N=1 $. Then, it may happen that $ \psi(0) \le 0 $, 
but since $ h(x) \sim \mathfrak{m} \Gamma(x) $ tends to $ - \infty < \psi (0) $ at infinity, 
the previous argument still applies.\qed

If one is able to prove that $ \cak $ is convex and assuming that $ \Phi_{\rm ext} $ 
satisfies H2), then H1) is automatically true. Indeed, any point of $ \partial \cak $ 
has then a positive density and we may then apply the regularity result of 
L. Caffarelli (see {\it e.g.}, \cite[Chapter 2, Theorem 3.10]{Friedman}) which ensures 
that $ \p\Omega $ is of class $C^1$ (hence $C^{s+1}$ by H2)).

We conclude with a result from \cite[Theorem 3.24]{HanMak} on the topology of $ \cak $ valid 
only in space dimension two (the proof uses complex analysis).

\begin{prop} [\cite{HanMak}]
\label{topo} 
We assume $N=2$. Suppose that $ \Phi_{\rm ext} $ is of class $ C^2 $ and that its 
Hessian is everywhere positive definite. Then, $ {\rm supp}( n_{\rm e} ) $ is simply connected, 
and equal to the closure of its interior. Moreover, if $ \Phi_{\rm ext} $ is $C^{2,\alpha}$ 
for some $ \alpha \in ]0,1[ $, then $ \p \cak $ is a $ C^{1,\beta} $ Jordan curve, for 
some $ \beta \in ] 0, 1[ $.
\end{prop}

The above result does not prevent cusps in $\p \cak $, but just says that the 
boundary $ \p \cak $ possesses a $ C^{1,\beta} $ parametrization.

\section{Asymptotic analysis}

This section is devoted to the analysis of the asymptotic regime $\eps\to 0$. 
We shall point out the difficulties and necessary adaptations between the case of 
quadratic potentials, Theorem \ref{thresultradial} and Theorem \ref{thresultpatate},
and the general case,  Theorem \ref{thresultlake}.
 For the
existence theory of the Vlasov--Poisson equation, we refer the reader
to \cite{Ars} for weak solutions and more recently to \cite{PLLP, Pfa}
where strong solutions and regularity issues are discussed. Further
details and references can be found in the survey \cite{Gla}.

\subsection{A useful estimate on $\Phi_{\rm e}$}\label{sec:lake2}

Before we turn to the analysis of the asymptotic regime $\eps\to 0$, 
it is convenient to set up an estimate that describes 
the behavior of $\Phi_{\rm e} $ close to the neighborhood of $\partial\cak$.
In the isotropic case, $\Phi_{\mathrm {ext}}$ being given by \eqref{iso}, 
the potential $\Phi_{\rm e}$ is defined by  \eqref{defPhie}, 
and we observe that there exists $ C>0$ such that
\begin{equation} 
\label{tropcool}
 0\leq ( |x|- R )\ | \nabla_x \Phi_{\mathrm{e}}(x) | \leq C\, \Phi_{\mathrm{e}}(x) 
\end{equation}
holds for any $x$ with $ | x| \ge R$. More generally, for a quadratic 
potential \eqref{quadra}, we can establish the following property, based on 
the formulas in Section \ref{sec:patate}.

\begin{lem}
\label{controlpatate}
Let $\Phi_{\mathrm e}$ be the quadratic potential defined as in Corollary \ref{omega_patate}. 
Let $\mathscr V:\mathbb R^N\to\mathbb R^N$ be smooth, compactly
supported and such that $\mathscr V\cdot \nu\big|_{\partial\cak_a}=0$.
Then, there exists a positive constant $C$, depending only on $ N $, $
\Phi_\ext $ and $\mathscr V$ such that we have, for any $ x \in
\mathbb R^N $,
\begin{equation}\label{encorepluscool}
|  \mathscr{V} \cdot \nabla \Phi_{\rm e} (x)| \le C \Phi_{\rm e} (x).
\end{equation}
\end{lem}

\noindent {\it Proof.} Since $ \mathscr{V} $ is compactly supported
and $ \Phi_{\rm e} $ is positive in $ \{ \sigma_a > 0 \} $,  we just need
to prove the inequality for $x$ close to $ \p \cak_a $, that is for
$ \sigma_a (x) $ small. We still define $ \lambda > 0 $ so that
$  \lambda^{-2} = \sum_{j=1}^N \lambda_j^{-2} $. From \eqref{expresso},
and by Taylor expansion of the integral, we infer that for $ 0 < \sigma_a (x) \ll 1 $
and $ 1 \le k \le N $,
$$
\lambda^2 \p_k \Phi_{\rm e} (x)
=  \frac{x _k}{2}\ds \left( \prod_{j=1}^N a_j \right) 
\left( \sigma_a(x) \frac{1}{a_k^2 } \, \left( \ds \prod_{j=1}^N a_j^2 \right)^{-1/2} 
+ \mathcal{O} ( \sigma_a^2 (x) ) \right)
= \frac{ x _k \sigma_a(x) }{2 a_k^2 } + \mathcal{O} ( \sigma_a^2 (x) ) .
$$
Let $\mathscr X(x)$ stands for the vector with components $x_k/a_k^2$. 
In particular, for $ 0 < \sigma_a (x) \ll 1 $, we get
\begin{align}
\label{unitaire}
 | \nabla \Phi_{\rm e}(x)| = \mathcal{O} (\sigma_a (x) )
 \quad \quad \quad {\rm and} \quad \quad \quad
 \frac{\nabla \Phi_{\rm e} (x)}{ | \nabla \Phi_{\rm e} (x)| }
  = \frac{ \mathscr X(x)}{ |\mathscr X(x)| }
  + \mathcal{O} ( \sigma_a (x) ) ,
\end{align}
where the unit vector field
$x\mapsto \frac{\mathscr X(x)}{|\mathscr X(x)|}  $
is smooth near $ \p \cak_a $ and is the (outward) normal  on
$ \p \cak_a $. Now, observe that
\begin{align*}
 0 = & \, \p_k \left( \sum_{j=1}^N \frac{x_j^2}{a_j^2 + \sigma_a(x) } \right)
 = 2 \frac{x_k}{a_k^2 + \sigma_a(x) }
 - \left( \sum_{j=1}^N \frac{x_j^2}{ ( a_j^2 + \sigma_a(x))^2 } \right) \p_k \sigma_a(x)
 \\
 = & \,
 2 \frac{x_k}{a_k^2 }
 - \p_k \sigma_a(x) \left( \sum_{j=1}^N \frac{x_j^2}{ a_j^4 } \right) + \mathcal{O}( \sigma_a(x) ).
\end{align*}
Therefore, for $ 0 < \sigma_a (x) \ll 1 $ and $ 1 \le k \le N $, we have
\begin{align*}
\lambda^2 \p_k \Phi_{\rm e} (x)
= \frac{1}{4} \sigma_a(x) \p_k \sigma_a(x)
\left( \sum_{j=1}^N \frac{x_j^2}{ a_j^4 } \right) + \mathcal{O} ( \sigma_a^2 (x) )
= \frac{1}{8} \p_k \left( \sigma_a^2 (x) \left( \sum_{j=1}^N \frac{x_j^2}{ a_j^4 } \right) \right)
+ \mathcal{O} ( \sigma_a^2 (x) ) .
\end{align*}
As a consequence,
\begin{equation}
\label{minotor}
 \lambda^2 \Phi_{\rm e} (x) =
 \frac{1}{8} \sigma_a^2 (x) \left( \sum_{j=1}^N \frac{x_j^2}{ a_j^4 } \right)
 + \mathcal{O} ( \sigma_a^3 (x) ) \ge \frac{\sigma_a^2 (x)}{C },
\end{equation}
holds for some $C>0$. Going back to  \eqref{unitaire}, we arrive at
\begin{align*}
\mathscr{V}(x) \cdot \nabla \Phi_{\rm e}(x)
= &\, \mathscr{V} (x)\cdot \left( \frac{\nabla \Phi_{\rm e}(x)}{| \nabla \Phi_{\rm e}(x)| } \right)
\times
| \nabla \Phi_{\rm e}(x)| 
= \mathscr{V} (x)\cdot \left( \frac{ \mathscr X(x)}{ |\mathscr X(x)|}
  + \mathcal{O} ( \sigma_a (x) ) \right)
\times \mathcal{O} ( \sigma_a (x) ) \\
= & \, \left( \mathcal{O} ( \sigma_a (x) ) + \mathcal{O} ( \sigma_a (x) ) \right)
\times \mathcal{O} ( \sigma_a (x) )
= \mathcal{O} ( \sigma_a^2(x) ) = \mathcal{O} ( \Phi_{\rm e}(x) ),
\end{align*}
by \eqref{minotor} and since $ \mathscr{V} \cdot \frac{ \mathscr X}{ |\mathscr X|} $
vanishes when $ \sigma_a = 0 $ in view of the no flux condition satisfied by $\mathscr V$. 
This finishes the proof. \qed

In the more general setting considered in Theorem \ref{thresultlake}, the result is the 
following and simply relies on the use of a local chart.

\begin{lem}
\label{distaucarre} 
We assume that $ \p \Omega $ is of class $ C^1 $ and that H2) is satisfied. 
Then, there exists a constant $C$ such that, for any $ x \in \mathbb R^N $,
\begin{equation}
\label{megacool}
| \mathscr{V} \cdot \nabla \Phi_{\rm e} (x)| \le C \Phi_{\rm e} (x).
\end{equation}
\end{lem}

\noindent {\it Proof.} We have already seen that $ \p \Omega $ is actually of class $C^{s+1} $. 
Since $ \Phi_{\rm e} $ is positive in $ \mathbb R^N \setminus \mathcal{K} $ and $ \mathscr{V} $ has 
compact support, by a compactness argument, it suffices to show that \eqref{megacool} holds 
near any point $ a \in \p \Omega $. Possibly translating and rotating the axis, we assume 
$ a = 0 $ and that the inward normal to $ \Omega $ at $ a=0 $ is $ e_1 = ( 1, 0, ... , 0 ) $. 
We let $ x_1 = \Theta ( x_\perp ) $, where $ x_\perp = ( x_2 , ... , x_N ) $, 
be a $ C^2 $ parametrization of $\p \Omega $ near $ 0 $, with $ \nabla \Theta (0) = 0 $, 
hence $ \Theta ( x_\perp ) = \mathcal{O} ( \lvert x_\perp \rvert^2 ) $. 

We now consider the function $ \vp : \R^N \to \R $ defined by 
$ \vp( y) \ddef \Phi_{\rm e} ( y_1 + \Theta (y_\perp), y_\perp ) $, where 
$ y_\perp = ( y_2, ... , y_N ) \in \R^{N-1} $. Then, $ \vp( y) = 0 $ when $ y_1 \ge 0 $, 
hence, for $ 2 \le j \le N $ and $ 1 \le k \le N $ and if $ y_1 = 0 $, 
$ \p_k \vp (y) = \p^2_{j,k} \vp (y) = 0 $; moreover, 
$ \p^2_{1,1} \vp (0,y_\perp) = \Delta \Phi_{\rm ext} ( 0,y_\perp) $ in view of the equality 
$ \Delta \Phi_{\rm ext} (x) = \Delta_x \Phi_{\rm e} (x) = 
( \Delta_y \vp - (\Delta_\perp \Theta ) \p_1 \vp + \sum_{j=2}^N ( \p_j \Theta )^2 \p^2_{1,j} \vp ) 
( x_1 - \Theta (x_\perp ), x_\perp ) $ in $ \{ x_1 \le \Theta (x_\perp ) \} $. 

It follows from these relations that, by the Taylor formula and by using $ \Delta \Phi_{\rm ext} ( 0) > 0 $ 
and $y_1 = x_1 - \Theta (x_\perp) \leq 0 $,
$$
 \vp(y) = \vp(y) - \vp( 0, y_\perp) - y_1 \p_1 \vp( 0, y_\perp) 
 = y_1^2 \int_0^1 (1-t) \p^2_1 \vp( ty_1, y_\perp) \, \ud t 
 \ge \frac{y_1^2}{C} ,
$$
and we deduce 
\begin{equation}
\label{minor}
 \Phi_{\rm e} ( x ) \ge \frac{(x_1 - \Theta (x_\perp))^2}{C} .
\end{equation}
Still by the Taylor formula, we have, for $ 2 \le j \le N $,
$$
 \p_j \vp (y) = y_1^2 \int_0^1 ( 1 - t ) \p^3_{1,1,j} \vp ( t y_1 , y_\perp ) \, \ud t 
 = \mathcal{O} (y_1^2 )
$$
and
$$
 \p_1 \vp (y) = y_1 \p^2_1 \vp( 0, y_\perp ) 
 + y_1^2 \int_0^1 ( 1 - t ) \p^3_{1} \vp ( t y_1 , y_\perp ) \, \ud t 
 = y_1 \Delta \Phi_{\rm ext} ( 0,y_\perp) + \mathcal{O} (y_1^2 ) .
$$
Now, we write $ \p_1 \Phi_{\rm e} (x) = \p_1 \vp (y) $ (with $ y = ( x_1 - \Theta (x_\perp) , x_\perp ) $) 
and $ \nabla_\perp \Phi_{\rm e} (x) = \nabla_\perp \vp (y) - \p_1 \vp (y) \nabla_\perp \Theta (y_\perp) $, 
thus
$$
 \mathscr{V} (x) \cdot \nabla \Phi_{\rm e} (x) 
 = \mathscr{V}_1 (x) \p_1 \Phi_{\rm e} (x) + \mathscr{V}_\perp (x) \cdot \nabla_\perp \Phi_{\rm e} (x) 
 = \mathscr{V}_1 (x) \p_1 \vp (y) + \mathscr{V}_\perp (x) \cdot \nabla_\perp \vp (y) 
 - \p_1 \vp (y) \mathscr{V}_\perp (x) \cdot \nabla_\perp \Theta (y_\perp) . 
$$
Note that $ \nabla_\perp \vp (y) = \mathcal{O} (y_1^2 ) $. Furthermore, 
since $ \mathscr{V} \cdot \nu = 0 $ on $ \p \Omega = \{x_1 = \Theta (x_\perp ) \} $ and 
$ \nu(x) = (1 , - \nabla_\perp \Theta(x_\perp) )/ \lvert  (1 , - \nabla_\perp \Theta(x_\perp) ) \rvert $, 
we deduce
\begin{align*}
 \mathscr{V} (x) \cdot \nabla \Phi_{\rm e} (x) = & \, \mathcal{O} (y_1^2 ) 
 + \p_1 \vp (y) \Big( [ \mathscr{V}_1 (x) - \mathscr{V}_\perp (x) \cdot \nabla_\perp \Theta (x_\perp) ] 
 - [ \mathscr{V}_1 (\Theta (x_\perp ), x_\perp) - \mathscr{V}_\perp (\Theta (x_\perp ), x_\perp) \cdot \nabla_\perp \Theta (x_\perp) ] \Big) 
 \\ 
 = & \, \mathcal{O} (y_1^2 ) + \mathcal{O} (| y_1| ) \times \mathcal{O} (| x_1 - \Theta(x_\perp) | ) 
 = \mathcal{O} (y_1^2 ) .
\end{align*}
We conclude by using \eqref{minor}. \qed

\subsection{Basic a priori estimates}

Now that we have in hand the limiting density profile $n_{\mathrm{e}}$
and the associated potential field $\Phi_{\mathrm{e}}$, we derive some
basic a priori estimates from \eqref{VP1}--\eqref{VP2}. 

Using the splitting of Poisson equation as in \eqref{VP2psi},
\eqref{VP1} recasts as
\[
\partial_t f_\eps+v\cdot\nabla_x f_\eps
-\ds\frac1\eps\nabla_x\Phi_{\mathrm e}\cdot\nabla_v f_\eps-
\ds\frac{1}{\sqrt\eps}\nabla_x\Psi_{\eps}\cdot\nabla_v f_\eps =0.
\]
Let us compute the time variation of the following energies:
\begin{itemize}
\item Kinetic energy
\[\ds\frac{\ud}{\ud t}\ds\iint \ds\frac{\lvert v \rvert^2}{2}\ f_\eps\ud v\ud x= -\ds\frac{1}{\eps}\ds\iint v\cdot \nabla_x \Phi_{\mathrm e}\ f_\eps\ud v\ud x
-\ds\frac{1}{\sqrt \eps}\ds\iint v\cdot \nabla_x\Psi_{\eps}\ f_\eps\ud
v\ud x,
\]
\item Leading order potential energy
\[\ds\frac{\ud}{\ud t}\ds\iint \Phi_{\mathrm e}\ f_\eps\ud v\ud x= \ds\iint v\cdot  \nabla_x\Phi_{\mathrm e}\ f_\eps\ud v\ud x,
\]
\item Fluctuations potential energy
\[\begin{array}{lll}\ds\frac{\ud}{\ud t}\ds\frac12\ds\int |\nabla_x\Psi_{\eps}|^2\ud x&=& 
\ds\int \nabla_x\Psi_{\eps}\cdot \partial_t   \nabla_x\Psi_{\eps}\ud x
=-\ds\int \Psi_{\eps} \partial_t\Big(  \ds\frac{n_{\mathrm e}- \rho_\eps}{\sqrt\eps}\Big)\ud x
\\&=&-
\ds\int \Psi_{\eps}  \ds\frac{1}{\sqrt\eps}\nabla_x\cdot\Big(\ds\int vf_\eps\ud v\Big)\ud x
=
\ds\frac{1}{\sqrt\eps}\ds\iint v\cdot \nabla_x \Psi_{\eps} \ f_\eps\ud v\ud x.
\end{array}\]
\end{itemize}
By summing these relations, we conclude with the following claim
(which applies for all three cases for $\Phi_{\mathrm{ext}}$).

\begin{prop}\label{bdd}
The solution $(f_\eps,\Phi_\eps=\frac1\eps\Phi_{\mathrm e}+\frac{1}{\sqrt\eps}\Psi_\eps)$ of 
\eqref{VP1}--\eqref{VP2} satisfies the following energy conservation equality
\[\ds\frac{\ud}{\ud t}\left\{
\ds\iint \ds\frac{\lvert v \rvert^2}{2}\ f_\eps\ud v\ud x 
+ \ds\frac1\eps\ds\iint \Phi_{\mathrm e}\ f_\eps\ud v\ud x
+ \ds\frac12\ds\int |\nabla_x\Psi_{\eps}|^2\ud x\right\}=0.
\]
Furthermore, the total charge is conserved
\[
\ds\iint f_\eps(t,x,v)\ud v\ud x=\ds\iint f_\eps(0,x,v)\ud v\ud x=\mathfrak{m}.\]
\end{prop}

\subsection{Convergence of the density and the current}

We assume a uniform bound on the energy at the initial time, namely
\begin{equation}
 \label{bornener}
 \sup_{0 < \eps < 1} \ds\iint \frac12 \ds\lvert v \rvert^2\ f_\eps^{\mathrm{init}} \ud v\ud x 
 + \ds\frac1\eps\ds\iint \Phi_{\mathrm e}\ f_\eps^{\mathrm{init}}\ud v\ud x
 + \ds\frac12\ds\int |\nabla_x\Psi_{\eps}^{\mathrm{init}} |^2\ud x < \infty ,
\end{equation}
where $ \Psi_{\eps}^{\mathrm{init}} $ solves the Poisson equation
\eqref{VP2psi}. Then, Proposition \ref{bdd} ensures that the energy
remains uniformly bounded for positive times. Thus, possibly at the
price of extracting subsequences, we can suppose that
\[
f_\eps\rightharpoonup f \text{ weakly-$\star$ in $\mathscr
  M^1([0,T]\times\mathbb R^N\times\mathbb R^N)$}, \quad
\rho_\eps = \ds \int f_\eps \ud v \rightharpoonup \rho \text{ weakly-$\star$ in $\mathscr
  M^1([0,T]\times\mathbb R^N)$}.
\]
Going back to the Poisson equation, we observe that
\[
n_{\mathrm e}-\rho_\eps=\sqrt\eps\nabla_x\cdot (\nabla_x\Psi_\eps)\]
where, by Proposition \ref{bdd}, $\nabla_x\Psi_\eps$ is bounded in
$L^\infty(0,T;L^2(\mathbb R^N))$.  Consequently, we establish the
following claim.

\begin{lem}\label{l:cv1}
  The sequence $\rho_\eps$ converges to $ n_{\mathrm e}=\rho $
  strongly in $L^\infty( 0,T ;H^{-1}(\mathbb R^N))$ and
  weakly-$\star$ in $\mathscr M^1 ([0,T]\times\mathbb R^N)$.  The
  limit $f$ is supported in $[0,T] \times\bar \Omega \times \mathbb R^N$. 
  The sequence $J^\eps = \int v f_\eps \ud v $ is bounded in $L^\infty(0,T ;L^1(\mathbb R^N))$; 
  it admits a subsequence which 
  converges, say weakly$-\star$ in $ \mathscr M^1([0,T]\times\mathbb R^N)$; the limit
  $J $ is divergence free, supported in $ [0,T] \times\bar\Omega $ 
  and may be written $ \int vf\ud v=n_{\mathrm e}W$ for some $W\in
  \mathscr M^1([0,T]\times\mathbb R^N)$.
\end{lem}

\noindent
{\bf Proof.}
Proposition \ref{bdd} tells us that $\lvert v \rvert^2 f_\eps$ is
bounded in $L^\infty( 0,T ;L^1(\mathbb R^N \times \mathbb R^N))$.
Hence, by using Cauchy-Schwarz' inequality, we get
\begin{equation}\label{toto}
\ds\int |J_\eps|\ud x\leq \ds\iint |v|\ \sqrt{f_\eps}\ \sqrt {f_\eps}\ud v\ud x 
\leq \left(\ds\iint \lvert v \rvert^2\ f_\eps\ud v\ud x\right)^{1/2} 
\left(\ds\iint f_\eps\ud v\ud x\right)^{1/2},
\end{equation}
which leads to the asserted uniform estimate on the current.
We can thus also assume $J_\eps\rightharpoonup J$ weakly$-\star$ in 
$\mathscr  M^1([0,T]\times\mathbb R^N)$. Furthermore, since the second order 
moment in $v$ of $f_\eps$ is uniformly bounded, we check that 
\[
\rho=\ds\int f\ud v,\qquad J=\ds\int v\ f \ud v.\] Note that
$\rho_\eps$ and $J_\eps$ satisfy \eqref{ch_cons}.  Letting $\eps$ go
to 0 yields
\[\partial_ t \rho +\nabla_x\cdot J=0=\partial_ t n_{\mathrm e}
+\nabla_x\cdot J=0+\nabla_x\cdot J=0.
\]
Thus, $J$ is divergence-free.  Finally, since $\lim_{|x|\rightarrow
  \infty}\Phi_{\mathrm e}(x)=+\infty$ and the second order moment in
$v$ of $f_\eps$ is uniformly bounded, $ \{ f_\eps , \, \eps > 0 \}$ is tight, and we can
write
\begin{align*}
  \ds \int_0^T\iint f_\eps\ud v\ud x \ud t = & \, \mathfrak{m} T =
  \ds\int_0^T\int \rho_\eps\ud x\ud t
  \\
  \xrightarrow[\eps \rightarrow 0]{} & \int_0^T \iint f\ud v\ud x \ud t =
  \mathfrak{m} T = \int_0^T\ds\int\rho\ud x\ud t =\ds T \int
  n_{\mathrm e}\ud x
  \\
  &\quad= \ds \int_0^T \ds\iint_{\overline\Omega
  } f\ud v\ud x\ud t + \ds \int_0^T
  \iint_{\mathbb R^N\setminus \overline\Omega
  } f\ud v\ud x \ud t
  \\
  &\quad =\ds\int_0^T\int_{\overline\Omega
  } n_{\mathrm e}\ud x\ud t =\ds
  \int_0^T\int_{\overline\Omega
  } \rho\ud x\ud t
  \\
  &\quad =\ds\int_0^T\iint_{\overline\Omega
  } f\ud v\ud x\ud t =\ds\int_0^T
  \int_{\Omega
  } n_{\mathrm e}\ud x\ud t
  \\
  &\quad =\ds\int_0^T\int_{\Omega
  } \rho\ud x\ud t =\ds\int_0^T
  \iint_{\Omega
  } f\ud v\ud x\ud t.
\end{align*}
It proves that $\mathrm{supp}(f)\subset [0,T]\times \overline\Omega
\times
\mathbb R^N$, and thus $\mathrm{supp}(J)\subset [0,T]\times
\overline\Omega
$.  In particular, we note that $f([0,T]
\times\partial\Omega\times
 \mathbb R^N)=0$, and $J([0,T]
 \times\partial\Omega
 )=0$.  \qed

In order to define the normal trace of $J$ over $\partial\Omega$ (that is 
the sphere $\p B(0,R) $ in the case \eqref{iso}), we shall use the theory introduced in \cite{CF}. 
As a consequence of the discussion above, we start by observing that $J$
belongs to the set $\mathscr D\mathscr M^{\mathrm{ext}}(\mathbb R^N)$ 
of extended divergence-measure fields over $\mathbb R^N$, see
\cite[Definition 1.1]{CF}.  Therefore, according to \cite[Theorem 3.1]{CF}, $J$
admits a normal trace $J\cdot \nu\big|_{\partial\Omega}
$ 
defined as a
continuous linear functional over $Lip(\gamma, \partial\Omega
)$,
$\gamma>1$ (see \cite[Equation (2.1)]{CF}) with
\[
\big\langle J \cdot \nu\big|_{\partial\Omega
},\phi\big\rangle=
\ds\int_{\Omega}
 \hat \phi \nabla_x \cdot J + \ds\int_{\Omega
} J \cdot
\nabla_x\hat \phi ,
\]
where the function $ \hat \phi \in Lip(\gamma, \Omega
)$ in the
right-hand side is an extension of $\phi \in Lip(\gamma, \partial\Omega
)
$.  However, by $\nabla_x\cdot J=0$ and the support property on $J$,
we can rewrite
\[\big\langle J\cdot \nu\big|_{\partial\Omega
},\phi\big\rangle= 0 +  \ds\int_{\mathbb R^N} J \cdot \nabla_x \hat\phi =
-\big\langle \nabla_x\cdot J, \hat\phi\big\rangle =0.\] Another way to
see this is to observe that the normal trace from $\Omega
$ must be
the same as the normal trace from $ \mathbb R^N \setminus \overline\Omega
$, which is clearly zero since $ J$ has support in $ \overline\Omega
$. Consequently,
$$
J \cdot \nu\big|_{\overline\Omega
} = 0 \quad \quad {\rm in} \quad [ C( 0, T ; Lip(\gamma, \overline\Omega
) )]^* .
$$
Remark that this is not a pointwise relation. In particular, it 
may happen that $ J_\eps^{\mathrm{init}} \cdot \nu\big|_{\overline\Omega
} = \rho_\eps V_\eps^{\mathrm{init}} \cdot
\nu\big|_{\overline\Omega
} $ is nonzero, but this does not prevent the
time integral of $J \cdot \nu\big|_{\overline\Omega } $ to vanish.
\medskip 

\subsection{Passing to the limit: modulated energy}

We now study the modulated energy
\[\mathscr H_{\mathscr V,\eps} = 
\ds\frac12\ds\iint |v-\mathscr V|^2\ f_\eps\ud v\ud x+
\ds\frac12\ds\int |\nabla_x\Psi_\eps|^2\ud x +\ds\frac1\eps\ds\iint
\Phi_{\mathrm e}\ f_\eps\ud v\ud x ,\] 
where all the terms integrated are nonnegative. Let us compute as follows
\[
\ds\frac{\ud }{\ud t}\mathscr H_{\mathscr V,\eps}= \ds\frac{\ud }{\ud
  t}\ds\iint \Big(\ds\frac{ \lvert \mathscr V \rvert^2}2\ - \mathscr
V\cdot v\Big)f_\eps\ud v\ud x = \ds\frac{\ud }{\ud t}\ds\int
\Big(\rho_\eps \ds\frac{\lvert \mathscr V \rvert^2}2\ - \mathscr
V\cdot J_\eps\Big) \ud x,
\]
by using Proposition \ref{bdd}.
We thus have
\[
\ds\frac{\ud }{\ud t}\mathscr H_{\mathscr V,\eps}=
\ds\int (\rho_\eps \mathscr V-J_\eps)\cdot \partial_t \mathscr V \ud x +
\ds\int \ds\frac{ \lvert \mathscr {V} \rvert^2}{2} \partial_t \rho_\eps \ud x 
- \ds\int \mathscr V\cdot \partial_t J_\eps\ud x
.\]

Here, we are assuming that the solution 
 $f_\eps$ of the Vlasov--Poisson system 
\eqref{VP1}--\eqref{VP2} is regular enough so that we can 
perform all the calculations that follow. 
Integrating the Vlasov equation, we obtain 
\begin{equation}\label{eqmt}
\partial_t J_\eps +\nabla_x\cdot \mathbb P_\eps +\ds\frac{1}{\sqrt\eps} \rho_\eps\nabla_x\Psi_\eps 
+ \ds\frac1\eps\rho_\eps\nabla_x\Phi_{\mathrm e}=0, 
\end{equation}
where we rewrite 
\[
\ds\frac{1}{\sqrt\eps} \rho_\eps\nabla_x\Psi_\eps= \ds\frac{ \rho_\eps
  -n_{\mathrm e}}{\sqrt\eps}\nabla_x\Psi_\eps + \ds\frac{n_{\mathrm
    e}}{\sqrt\eps}\nabla_x\Psi_\eps =
-\Delta_x\Psi_\eps\nabla_x\Psi_\eps + \ds\frac{n_{\mathrm
    e}}{\sqrt\eps}\nabla_x\Psi_\eps,
\]
and 
\[
\Delta_x\Psi_\eps\nabla_x\Psi_\eps = 
\nabla_x\cdot\big(\nabla_x\Psi_\eps\otimes \nabla_x\Psi_\eps\big) 
-\nabla_x\Big(\ds\frac{|\nabla_x\Psi_\eps|^2}{2}\Big).
\]
Combining these relations to the charge conservation \eqref{ch_cons}
and integration by parts, we arrive at
\[\begin{array}{l}
  \ds\frac{\ud }{\ud t}\mathscr H_{\mathscr V,\eps}=
  \ds\int (\rho_\eps \mathscr V-J_\eps)\cdot \partial_t \mathscr V \ud x +
  \ds\int J_\eps \cdot \nabla_x\ds \left( \frac{\lvert \mathscr V \rvert^2}{2} \right) \ud x\\[.4cm]
  \qquad\qquad\qquad-
  \ds\int D_x\mathscr V:(\mathbb P_\eps-\nabla_x\Psi_\eps\otimes\nabla_x\Psi_\eps)\ud x
  +\ds\int \nabla_x\cdot\mathscr V\ds\frac{|\nabla_x\Psi_\eps|^2}{2}\ud x
  \\[.4cm]
  \qquad\qquad\qquad+\ds\frac1\eps\ds\int \rho_\eps \mathscr V\cdot \nabla_x\Phi_{\mathrm e}\ud x 
  + 
  \ds\int \frac{n_{\mathrm e}}{\sqrt\eps} \mathscr V\cdot \nabla_x\Psi_\eps\ud x 
  ,
\end{array}\]
where $D_x\mathscr V$ stands for the jacobian matrix of the vector field 
$\mathscr V$. For the last integral, since $ n_{\rm e} $ is supported in 
$ \Omega $, we write it as 
$$
\int_\Omega \frac{n_{\mathrm e}}{\sqrt\eps} \mathscr V\cdot \nabla_x\Psi_\eps\ud x = 0 
$$
by integration by parts and using that $\nabla_x\cdot (n_{\rm
  e}\mathscr V)=0$ in $\Omega$ and the no-flux condition \eqref{CLEuler}.

Let us set
$$
\mathbb P_{\mathscr V,\eps} \ddef \ds\int (v-\mathscr V)\otimes
(v-\mathscr V)\ f_\eps\ud v = \mathbb P_\eps - \mathscr V\otimes
J_\eps - J_\eps\otimes \mathscr V + \rho_\eps \mathscr V\otimes
\mathscr V .$$ A direct substitution leads to
 \begin{equation}\label{est_en_mod}\begin{array}{lll}
\ds\frac{\ud }{\ud t}\mathscr H_{\mathscr V,\eps}&=&
\ds\int (\rho_\eps \mathscr V-J_\eps)\cdot\big( \partial_t \mathscr V 
+(\mathscr V\cdot \nabla_x)\mathscr V\big) \ud x 
\\[.4cm]
&&- \ds\int D_x\mathscr V:(\mathbb P_{\mathscr V, \eps}-\nabla_x\Psi_\eps\otimes\nabla_x\Psi_\eps)\ud x
\\[.4cm]
&& +\ds\int \nabla_x\cdot\mathscr V\ds\frac{|\nabla_x\Psi_\eps|^2}{2}\ud x
+\ds\frac1\eps\ds\int \rho_\eps \mathscr V\cdot \nabla_x\Phi_{\mathrm e}\ud x.
\end{array}\end{equation}
We shall use the shorthand notation $A\lesssim B$ when the inequality
$A\leq CB$ holds for some constant $C>0$, the value of which might
vary from a line to another. As a matter of fact, we can dominate the second and 
third integrals 
of the right-hand side by
\[
\|D_x\mathscr V\|_\infty \left(\ds\iint |v-\mathscr V|^2f_\eps\ud v\ud x 
+ \ds\int |\nabla_x\Psi_\eps|^2\ud x\right) 
\le \|D_x \mathscr V\|_\infty\ \mathscr H_{\mathscr V,\eps}.
\]
Let us distinguish 
the case of the isotropic potential in order to point out the difficulties.
When $\Phi_{\mathrm{ext}}$ is given by \eqref{iso}, we remind the reader  that $\Phi_{\mathrm e} $
is supported in $\{|x| \ge R\}$, radially symmetric and increasing in
$|x|$, see \eqref{defPhie}. Combining this with \eqref{tropcool} allows us to estimate the
last term in \eqref{est_en_mod} as follows:
\[\begin{array}{lll}
  \left|\ds\frac1\eps\ds\int \rho_\eps \mathscr V\cdot \nabla_x\Phi_{\mathrm e}\ud x\right| & =&
  \left|\ds\frac1\eps\ds\int_{ \{ |x|> R \} } \rho_\eps \ds\frac{\mathscr V\cdot x/|x| }{|x|-R}\ (|x|-R) 
    |\nabla_x \Phi_{\mathrm{e}}| \ud x\right|
  \\
  &\lesssim& 
  \ds\frac1\eps\ds\int_{\{ |x|> R \} } \rho_\eps \Phi_{\mathrm e}\ud x\ \left\| 
    \ds\frac{\mathscr V\cdot x/|x| }{|x|-R}\right\|_\infty\lesssim \mathscr H_{\mathscr V,\eps},
\end{array}\]
where we have used that $ \frac{\mathscr V \cdot x/|x| }{|x|-R} $ belongs to 
$L^\infty ( \R^N ) $ since $ \mathscr V $ is smooth, compactly supported, and 
$ \mathscr V \cdot \nu = 0 $ on $ \partial B(0,R) $. 
For a quadratic external potential \eqref{quadra}, we can proceed similarly
by using Lemma~\ref{controlpatate}. When dealing with a general potential, we made 
hypothesis H2) so that Lemma \ref{distaucarre} applies and \eqref{encorepluscool} allows us to estimate 
\[ \left|\ds\frac1\eps\ds\int \rho_\eps \mathscr V\cdot \nabla_x\Phi_{\mathrm e}\ud x\right| 
\lesssim \ds\frac1\eps\ds\int \rho_\eps \Phi_{\mathrm e}\ud x
\lesssim \mathscr H_{\mathscr V,\eps}.
\]
Therefore, we obtain 
\begin{equation} \label{Hve}
\ds\frac{\ud }{\ud t}\mathscr H_{\mathscr V,\eps}\lesssim \mathscr H_{\mathscr V,\eps} + r_\eps
\end{equation}
where we have set 
\[r_\eps \ddef \ds\int (\rho_\eps \mathscr
V-J_\eps)\cdot\big( \partial_t \mathscr V +(\mathscr V\cdot
\nabla_x)\mathscr V\big) \ud x .\] The Gr\"onwall lemma yields
\[\mathscr H_{\mathscr V,\eps}(t)
\leq e^{Ct}\left(\mathscr H_{\mathscr V,\eps}(0)+\ds\int_0^t e^{-C
    \tau} r_\eps(\tau)\ud \tau \right),\] for a certain constant
$C>0$. The assumption \eqref{query} on the initial data is that 
$\lim_{\eps\rightarrow 0}\mathscr H_{\mathscr V,\eps}(0)=0$. Hence, 
we are left with the task of proving that $\int_0^t r_\eps(\tau)\ud
\tau$ tends to $0$ as $\eps\rightarrow 0$. We have
\[\begin{array}{lll}
\ds\int_0^t r_\eps(\tau)\ud \tau&\xrightarrow [\eps\rightarrow 0]{} &
\ds\int_0^t \int (n_{\mathrm e} \mathscr V-J)\cdot\big( \partial_t \mathscr V 
+(\mathscr V\cdot \nabla_x)\mathscr V\big) \ud x\ud \tau
\\
&& \qquad
=
\ds\int_0^t \int_{\Omega}
( n_{\rm e}\mathscr V-J)\cdot\big( \partial_t \mathscr V
+(\mathscr V\cdot \nabla_x)\mathscr V\big) \ud x\ud \tau
\\&&\qquad 
= - \ds\int_0^t \int_{\Omega
} (n_{\rm e} \mathscr V-J)\cdot\nabla_x p \ud x\ud \tau
=0,
\end{array}\]
since $n_{\rm e} \mathscr V$ and $J$ are divergence free on $\Omega$ 
 and their normal trace vanish.
\qed

It is worth pointing out that  the regularity assumption of the 
sequence of solutions  $ f_\eps$ was only made to justify 
the computations leading to  \eqref{Hve}. If one 
consider less regular solutions,  we have to assume 
that  these solutions were constructed through a 
regularization procedure  and that the previous 
 calculations were done on these regularizations and hence 
\eqref{Hve} will still hold.

\subsection{Identification of the limit}
Let us observe that if the initial datum satisfies \eqref{query}, then
\eqref{bornener} holds true.  Let us first justify i): we shall show
that $ \int \rho_\varepsilon \chi \ud x \rightarrow \int n_{\mathrm e}
\chi \ud x $ uniformly on $ [ 0 ,T ] $ as $ \varepsilon \rightarrow 0
$ for any $\chi\in C^0_0(\mathbb R^N)$.  We start by observing that
\begin{equation}\label{unif}
 \left|\int \rho_\varepsilon(t,x) \chi(x) \ud x\right|\leq \mathfrak m \|\chi\|_\infty
 \end{equation}
 holds for any $\chi\in C^0_0(\mathbb R^N)$.
Next, consider $ \chi \in C^1_c(\mathbb R^N ) $. 
The charge conservation \eqref{ch_cons} yields
$$
\frac{\ud }{ \ud t} \int \rho_\varepsilon(t,x) \chi(x) \ud x = \int
\p_t \rho_\varepsilon \chi \ud x = - \int \nabla_x \cdot
J_\varepsilon \, \chi \ud x = \int J_\varepsilon \cdot \nabla_x \chi \ud x ,
$$
hence the uniform bound \eqref{toto} on $ J_\varepsilon $ implies a
uniform bound on $ \frac{\ud }{ \ud t}\int \rho_\varepsilon(t,x)
\chi(x) \ud x $ for $0\leq t\leq T$.  By virtue of the Ascoli-Arzel\`a
theorem, the set $\big\{t\mapsto \int \rho_\varepsilon(t,x) \chi(x)
\ud x,\ \varepsilon>0\big\} $ is therefore relatively compact in
$C([0,T])$ for any fixed $\chi \in C^1_c(\mathbb R^N )$.  This
property extends to any $\chi \in C^0_0(\mathbb R^N )$ by virtue of
\eqref{unif}.  Indeed, for any $\delta >0$, we can pick $\chi_\delta
\in C^1_c(\mathbb R^N )$ such that $\|\chi-\chi_\delta\|_\infty \leq \delta
/\mathfrak m$.  It follows that
\[
\ds \int \rho_\varepsilon(t,x) \chi(x) \ud x=\ds \int
\rho_\varepsilon(t,x) (\chi-\chi_\delta)(x) \ud x+ \ds \int
\rho_\varepsilon(t,x) \chi_\delta(x) \ud x
\]
where, owing to \eqref{unif}, the former integral is uniformly
dominated by $\delta $ and the latter lies in a compact set of
$C([0,T])$.  Therefore $\big\{t\mapsto \int \rho_\varepsilon(t,x)
\chi(x) \ud x,\ \varepsilon>0\big\} $ can be covered by a finite
number of balls with radius $2\delta $ in $ C([0,T] ) $.  Finally,
since $C_0^0(\mathbb R^N)$ is separable, we apply a diagonal argument
to extract a subsequence such that $ \int \rho_\varepsilon(t,x)
\chi(x)\ud x$ converges uniformly in $ C([0,T] ) $ for any element
$\chi$ of a numerable dense set in $C_0^0(\mathbb R^N)$.  By
uniqueness of the limit, we find
\[
\ds\lim_{\varepsilon \rightarrow 0}\ds\int \rho_\varepsilon(t,x)
\chi(x)\ud x=\int n_{\mathrm e} \chi \ud x.
\]
Going back to \eqref{unif}, we check that the convergence holds for
any $\chi\in C_0^0(\mathbb R^N)$.

The manipulations detailed in the previous Section prove ii). 
In order to establish 
iii), it is convenient to introduce the following functional: given
$\lambda$ a non negative bounded measure on $[0,T]\times \mathbb R^N$,
and $\mu$ a vector valued bounded measure on $[0,T]\times \mathbb
R^N$, we set
\[\mathscr K(\lambda,\mu) \ddef 
\sup_{\Theta}\left\{\ds\int \mu\cdot \Theta - \ds\frac12\ds\int
  \lambda |\Theta|^2 \right\}\] where the supremum is taken over
continuous functions $\Theta:[0,T]\times \mathbb R^N\rightarrow
\mathbb R^N$. According to \cite[Prop. 3.4]{VxBr}, we have:

\begin{lem}[\cite{VxBr}]\label{l:Vx}
  If $\mu$ is absolutely continuous with respect to $\lambda$,
  denoting by $\mathbb V$ the Radon-Nikodym derivative of $\mu$ with
  respect to $\lambda$, we have
\[
\mathscr K(\lambda,\mu)=\ds\frac12\ds\int \lambda |\mathbb V|^2\in [0,\infty],
\]
otherwise $ \mathscr K(\lambda,\mu)=+\infty$.
\end{lem}

Clearly $(\lambda,\mu)\mapsto \mathscr K(\lambda,\mu)$ is a convex and
lower semi--continuous (for the weak-$\star$ convergence)
functional. Let $ \eta : [0,T]\rightarrow [0,\infty)$ be a continuous
non negative function. Reasoning as in \cite{Bre}, we show that $J\in
L^\infty(0,T;L^2(\mathbb R^N))$ since
\[\begin{array}{lll}
  \mathscr K(\eta \rho_\eps,\eta J_\eps)&=& \ds\frac12\ds\int_0^T\ds\int_{\mathbb R^N} 
  \ds\frac{|J_\eps(t,x)|^2}{\rho_\eps(t,x)}\ \eta(t)\ud x\ud t 
  \\
  &=&
  \ds\frac12 \ds\int_0^T\ds\int_{\mathbb R^N} \ds\frac{1}{\rho_\eps(t,x)}\ 
  \left \lvert \ds\int_{\mathbb R^N} v\sqrt{f_\eps(t,x,v)}\ \sqrt{f_\eps(t,x,v)}\ud v\right \rvert^2\eta(t)\ud x\ud t 
  \\&\leq& 
  \ds\frac12 \ds\int_0^T\ds\iint_{\mathbb R^N\times \mathbb R^N} \lvert v \rvert^2 f_\eps(t,x,v)\ \eta(t)\ud v\ud x\ud t 
  \lesssim \ds\int_0^T \eta\ud t
\end{array}\]
becomes, as $\eps$ tends to 0
\[\mathscr K(\eta n_{\mathrm e} ,\eta J)\lesssim \|\eta\|_{L^1(0,T)}.\]
Reasoning the same way, we get
\[\begin{array}{lll}
  \mathscr K( \rho_\eps, J_\eps-\rho_\eps \mathscr V)
  &=&  \ds\frac12\ds\int_0^T\ds\int_{\mathbb R^N} \ds\frac{|J_\eps-\rho_\eps \mathscr V|^2}{\rho_\eps}\ud x\ud t
  \\& \leq &
  \ds\frac12 \ds\int_0^T\ds\iint_{\mathbb R^N\times \mathbb R^N} |v-\mathscr V|^2 f_\eps\ud v\ud x\ud t 
  \leq \ds\int_0^T  \mathscr H_{\mathscr V,\eps}\ud t.
\end{array}\]
It follows that $\mathscr K(n_{\mathrm e}, J-n_{\mathrm e} \mathscr V)=0$, which identifies 
the limit $J$ and ends the proof of iii). 

Finally, we can check that the initial data for the limit equation is
meaningful by establishing some time--compactness on the sequence
$J_\eps$. Let
\[
\mathscr W_R=\big\{ \Theta:[0,T]\times\mathbb R^N\rightarrow \mathbb
R^N,\ \Theta \text{ of class $C^1$},\ \mathrm{supp} (\Theta)\subset
[0,T]\times \overline\Omega, 
\ \nabla_x\cdot (n_{\rm e}\Theta)=0\big\},
\]
which is a closed subspace of the Banach space $C^1$ (endowed with the sup norm for 
the function and its first order derivatives). 
Multiplying \eqref{eqmt} by a function in $\mathscr W_R$, we shall get rid of the stiff terms. 
Indeed, for such a trial function $\Theta $, we deduce from \eqref{eqmt}
\begin{align}
 \label{eqmttheta}
 \ds\frac{\ud}{\ud t}\ds\int J_\eps \cdot \Theta\ud x 
 = 
 \ds \int J_\eps \cdot \partial_t \Theta\ud x 
 - \ds\int \Theta \cdot (\nabla_x\cdot \mathbb P_\eps) \ud x 
 - \ds \frac{1}{\sqrt\eps} \int \rho_\eps \Theta \cdot \nabla_x\Psi_{\eps} \ud x ,
\end{align}
since $ \Theta \cdot \nabla_x \Phi_{\rm e} = 0 $ pointwise in view of
the supports.  By using the estimates deduced from Proposition
\ref{bdd}, we observe that the first two terms are bounded in $
L^\infty(0,T) $. For the last one, we use the Poisson equation
\eqref{VP2psi} and integration by parts to infer
\begin{align*}
\frac{1}{\sqrt\eps} \int \rho_\eps \Theta \cdot \nabla_x\Psi_{\eps} \ud x 
= & \ 
\frac{1}{\sqrt\eps} \int n_{\mathrm e} \Theta \cdot \nabla_x\Psi_{\eps} \ud x 
- \int \Delta_x \Psi_{\eps} \Theta \cdot \nabla_x\Psi_{\eps} \ud x 
\\ = & \ 
0 + \int \nabla_x \Psi_{\eps} \cdot \nabla_x( \Theta \cdot \nabla_x\Psi_{\eps} ) \ud x 
\\ = & \ 
\frac12 \int \Theta \cdot \nabla_x ( \lvert \nabla_x\Psi_{\eps} \rvert^2 ) \ud x 
+ \sum_{1 \le j,k \le N} \int \partial_{x_j} \Psi_{\eps} \partial_{x_j} \Theta_k \partial_{x_k} \Psi_{\eps} \ud x
\end{align*}
where we have used that $n_{\rm e} \Theta (t,\cdot) $ is divergence free. 
For quadratic external potentials, an integration by parts shows that the 
first integral is zero (since $n_{\rm e} \Theta $ is divergence free).
In any cases, the right hand side can be dominated by $\|\nabla\Theta\|_{\infty}\|\nabla\Psi_\eps\|_{L^\infty(0,T;L^2(\mathbb R^N))}$ and it is thus 
 bounded in $ L^\infty(0,T) $. Reporting this into \eqref{eqmttheta} allows us to conclude that
\[
\ds\frac{\ud}{\ud t}\ds\int J_\eps \cdot \Theta\ud x
\text{ is bounded in $L^\infty( 0,T )$. }
\]
Since $\mathscr W_R$ is separable, we can boil down a diagonal
argument to justify that $J_\eps$ is relatively compact in
$C^0(0,T;\mathscr W_R'-\text{weak}-\star)$: we can assume that the
extracted subsequence is such that $\int J_\eps \cdot \Theta\ud x$
converges uniformly on $[0,T]$ for any $\Theta \in \mathscr W_R$.
\qed

\section{Asymptotic analysis of the Vlasov--Poisson--Fokker--Planck system}
\label{sec:VFP}

In this Section we state and prove a Theorem analogous
to Theorem \ref{thresultradial} when the basic equation is \eqref{FPV1}, which 
includes a Fokker--Planck operator, coupled with \eqref{VP2}.

For the well-posedness issues of the system \eqref{FPV1} coupled to \eqref{VP2}, we refer the reader
to \cite{Bo, Deg86}. The role of the external potential is precisely investigated in \cite{Dol}.
The associated moment system reads
\[\left\{\begin{array}{l}
\partial_t \rho_\eps+\nabla_x\cdot J_\eps=0,
\\
 \partial_t J_\eps+\nabla_x\cdot \mathbb P_\eps+\rho_\eps\nabla_x\Phi_\eps=- J_\eps,
\end{array}\right.
\]
where we still use the notation $J_\eps=\int vf_\eps\ud v$, $\mathbb P_\eps=
\int v\otimes vf_\eps\ud v$. 
As $\eps\rightarrow 0$, we expect as before that 
$ \rho_\eps \to n_{\mathrm e} = \mathbf 1_{\Omega} \Delta \Phi_{\rm ext}$ and that the behavior 
of the current is driven by the Lake Equation with friction
\begin{equation}
\tag{LE$_f$}\label{LE2}
\left\{ 
\begin{array}{l}
\partial_t V + V \cdot \nabla_x V + \nabla_x p =- V,
\\
\nabla_x\cdot ( n_{\rm e} V ) = 0.
\end{array}\right.
\end{equation}
If $ \Phi_{\rm ext} $ is quadratic as in  \eqref{quadra} (possibly isotropic), the domain $\Omega $ is an ellipsoid 
(possibly a ball) as in Section \ref{sec:patate} and \eqref{LE2} becomes the Incompressible Euler 
system with friction
\begin{equation}
\left\{ \label{Euler2}
\begin{array}{l}
\partial_t V + \nabla_x\cdot(V\otimes V) + \nabla_x p = -V,
\\
\nabla_x\cdot V = 0.
\end{array}\right.
\end{equation}
For a more general confining potential $ \Phi_{\rm ext} $, we make assumptions h1), h2), H1) and H2) as 
in Section \ref{sec:lake}. 
Since we work with finite charge data, the limit equation \eqref{LE2} holds in $ \Omega$, 
completed with the no flux boundary condition \eqref{CLEuler}, namely
\begin{equation*}
V(t,x)\cdot\nu(x)\Big|_{\p \Omega}=0.
\end{equation*}
Like in the previous section we associate with $V$, smooth solution of
\eqref{LE2}, a smooth compactly supported extension $\mathscr V$
defined on $[0,T]\times \mathbb R^N$ such that $  \mathscr V
\cdot\nu(x)\Big|_{\p \Omega}=0$.  


\medskip

We shall investigate this asymptotics in the specific case where the
``temperature'' $\theta=\theta_\eps$ goes to 0 as $\eps\rightarrow 0$.
In this context, we can derive an analog of Proposition \ref{bdd} that
accounts for the dissipation mechanisms induced by the Fokker--Planck
operator.

\begin{prop}\label{bdd2}
  The solution $(f_\eps,\Phi_\eps=\frac1\eps\Phi_{\mathrm
    e}+\frac{1}{\sqrt\eps}\Psi_\eps)$ of \eqref{FPV1}--\eqref{VP2}
  satisfies the following entropy dissipation inequality
\[\ds\frac{\ud}{\ud t}\left\{
\ds \frac12\iint \ds\lvert v \rvert^2\ f_\eps\ud v\ud x 
+ \ds\frac1\eps\ds\iint \Phi_{\mathrm e}\ f_\eps\ud v\ud x 
+ \theta_\eps \ds\iint f_\eps\ln(f_\eps)\ud v\ud x
+ \ds\frac12\ds\int |\nabla_x\Psi_{\eps}|^2\ud x\right\}=-\mathscr D_\eps\]
where we denote
\[
\mathscr D_\eps=\ds\iint |v\sqrt{f_\eps}+2\theta_\eps\nabla_v\sqrt{f_\eps}|^2\ud v\ud x\geq 0.
\]
Furthermore, the total charge is conserved
\[
\ds\iint f_\eps(t,x,v)\ud v\ud x=\ds\iint f_\eps(0,x,v)\ud v\ud x= \mathfrak{m}.\]
\end{prop}

\noindent
Uniform estimates are not directly included in this statement since
the function $z\mapsto z\ln (z)$ changes sign.  Nevertheless, we can
establish such uniform estimates.

\begin{coro}
\label{c:bdd2}
We assume that there exists some (large) $ \lambda > 1 $ such that
\begin{equation}
\label{integral}
 \int \exp( - \lambda \Phi_\ext ) \ud x < \infty .
\end{equation}
We suppose also that $0 < \eps \leq 1/ ( 8 \lambda) $ and $0 < \theta_\eps \le 1$. 
Let $f_\eps^{\mathrm{init}}:\mathbb R^N\times \mathbb
R^N\rightarrow [0,\infty)$ be a sequence of integrable functions
that satisfy the following requirements
\begin{equation}
\label{queryVFP}
\begin{array}{l}
\ds\iint f_\eps^{\mathrm{init}}\ud v\ud x=\mathfrak{m},
\\
\ds\sup_{0<\eps\leq 1/ ( 8 \lambda),\ 0<\theta_\eps \le 1}  \left\{
\ds\frac12\ds\iint
|v|^2\ f_\eps^{\mathrm{init}}\ud v\ud x+\theta_\eps \ds\iint
f_\eps^{\mathrm{init}} \lvert \ln (f_\eps^{\mathrm{init}}) \rvert \ud v\ud x
\right.
\\
 \qquad\qquad\qquad\qquad\qquad\qquad\left.
+ \ds\frac12\ds\int
|\nabla_x\Psi_\eps^{\mathrm{init}}|^2\ud x
+\ds\frac1\eps\ds\iint
\Phi_{\mathrm e}\ f_\eps^{\mathrm{init}}\ud v\ud x\right\}<\infty,
\end{array}
\end{equation}
with $$\Delta_x \Psi_\eps^{\mathrm{init}}=\ds\frac{1}{\sqrt\eps}\Big(n_{\mathrm{e}} 
-\ds\int f_\eps^{\mathrm{init}}\ud v\Big).$$
Let $0<T<\infty$ and let $(f_\eps,\Phi_\eps=\frac1\eps\Phi_{\mathrm e}+\frac{1}{\sqrt\eps}\Psi_\eps)$ 
be the associated solution of \eqref{FPV1}--\eqref{VP2}. Then, uniformly for $ 0 < \eps \le 1/ ( 8 \lambda) $ 
and $ 0 < \theta_\eps \le 1 $:
\begin{itemize}
\item [i)] $f_\eps(1+ \lvert v \rvert^2 +\theta_\eps \lvert \ln (f_\eps) \rvert)
+ \eps^{-1} \Phi_{\mathrm e}f_\eps$ is bounded in 
$L^\infty(0,T;L^1(\mathbb R^N\times \mathbb R^N))$ ,
\item [ii)] $\nabla_x\Psi_\eps$ is bounded in $L^\infty(0,T ;L^2(\mathbb R^N))$,
\item [iii)] $\mathscr D_\eps$ is bounded in $L^1(0,T )$.
\end{itemize}
\end{coro}

\begin{remark} 
\label{crapaud}
 In any dimension $ N \ge 1 $, \eqref{integral} is always true
 for quadratic potentials. If $N=1 $ or $N=2 $, hypothesis \eqref{integral} 
 is satisfied if h2) is, since $ \Phi_{\rm ext} (x) + \mathfrak{m} \Gamma (x) \to +\infty $ 
 when $ |x| \to + \infty $ and that $ \Gamma (x) = - |x| /2 $ or $ - \ln |x| / ( 2\pi ) $. 
 Therefore, hypothesis \eqref{integral} needs to be verified only for $N \ge 3 $.
\end{remark}

\noindent
{\bf Proof.} We first observe that hypothesis \eqref{integral} implies
$$
  \int \exp( - \lambda \Phi_{\rm e} ) \ud x < \infty .
$$
Indeed, we have $ \Phi_{\rm e} = \Gamma \star n_{\rm e} - C_* + \Phi_{\rm ext} \ge \Phi_{\rm ext} - C_* $. 
We write, for $h \geq 0 $,
\begin{align}
\label{tutu}
f_\eps\ln(f_\eps) \le f_\eps \lvert \ln(f_\eps)\rvert = & \ f_\eps\ln(f_\eps) 
- 2f_\eps\ln(f_\eps)\big( \mathbf 1_{e^{-h}\leq f_\eps\leq 1}
+\mathbf 1_{0\leq f_\eps <  e^{-h}}\big)
\nonumber \\
\leq& \ 
f_\eps\ln(f_\eps) +2h f_\eps+\ds\frac4e e^{-h/2} ,
\end{align}
and denote
$$
 E_\eps (f_\eps) = \frac12\iint \ds\lvert v \rvert^2\ f_\eps\ud v\ud x
 + \ds\frac1\eps\ds\iint \Phi_{\mathrm e}\ f_\eps\ud v\ud x 
+ \ds\frac12\ds\int |\nabla_x\Psi_{\eps}|^2\ud x. 
$$
We now use \eqref{tutu} with $h(x,v)= \lvert v \rvert^2 / (8 \theta_\eps)  + \Phi_{\mathrm e}(x)/(4\eps \theta_\eps)$ 
to infer
\begin{align}
\label{txtx} 
  \theta_\eps \iint f_\eps\ln(f_\eps) \ud v\ud x \le \theta_\eps \iint f_\eps 
  \lvert \ln(f_\eps)\rvert \ud v\ud x  \le 
  & \, \theta_\eps \iint f_\eps\ln(f_\eps) \ud v\ud x
  + \frac{\theta_\eps}{2} E_{\eps}(f_\eps) 
  \\ 
  \nonumber & \, + \theta_\eps \frac4e \iint \exp(- \lvert v \rvert^2 /(16 \theta_\eps ) 
  - \Phi_{\mathrm e}(x)/( 8\eps\theta_\eps) ) \ud v\ud x .
\end{align}
The last term is equal to
$$
 \theta_\eps \frac4e \int \exp(- \lvert v \rvert^2 /(16 \theta_\eps ) \ud v 
   \int \exp( -\Phi_{\mathrm e}(x)/( 8\eps\theta_\eps) )\ud x 
   \leq 
   \theta_\eps \frac4e \int \exp(- \lvert v \rvert^2 /16 ) \ud v 
   \int \exp( -\lambda \Phi_{\mathrm e}(x) )\ud x ,
$$
thus tends to zero as $\theta_\eps \rightarrow 0 $ (uniformly for $ 0 < \eps < 1/(8\lambda) $). 
Using the dissipation of the entropy given in Proposition \ref{bdd2},
we then infer
\begin{align*} 
E_{\eps}(f_\eps^\mathrm{init}) + \theta_\eps \iint f_\eps^{\mathrm{init}} 
  \lvert \ln(f_\eps^{\mathrm{init}}) \rvert \ud v\ud x 
\ge & \ 
E_{\eps}(f_\eps^\mathrm{init}) + \theta_\eps \iint f_\eps^{\mathrm{init}} \ln(f_\eps^{\mathrm{init}}) \ud v\ud x 
\\ \ge & \ 
E_{\eps}(f_\eps) + \theta_\eps \iint f_\eps \ln(f_\eps) \ud v\ud x 
\\ \ge & \ 
E_{\eps}(f_\eps) + \theta_\eps \iint f_\eps \lvert \ln(f_\eps) \rvert \ud v\ud x 
- \frac{\theta_\eps}{2} E_{\eps}(f_\eps)+ o_{ \theta_\eps \rightarrow 0 }(1) 
\end{align*}
and the conclusion follows since $ \theta_\eps \le 1 $.
\qed

Untill the end of the Section, we shall make hypothesis \eqref{integral}. 
Since we are dealing with the regime
\[0<\eps\ll 1,\qquad 0<\theta_\eps\ll 1, \] the estimates in
Proposition \ref{c:bdd2} do not provide $L^1$-weak compactness on the
particle distribution function and its moments; we still need to work
with convergences in spaces of finite measures.  The first step in the
investigation of the asymptotic behavior is summarized in the
following claim.

\begin{lem}
We make assumptions \eqref{integral} and \eqref{queryVFP}. 
  Up to a subsequence, we can assume that $f_\eps$ converges to $f$
  weakly--$\star$ in $\mathscr M^1((0,T)\times \mathbb
  R^N\times\mathbb R^N)$.  Then, $\rho_\eps$ converges to $n_{\mathrm
    e}=\int f\ud v$ in $L^\infty (0,T;H^{-1}(\mathbb R^N))$ and in
  $C^0(0,T ; \mathscr M^1(\mathbb
  R^N)-\text{weak--}\star)$. 
  Moreover, we can assume that $J_\eps\rightharpoonup J=\int vf\ud v$
  in 
  $\mathscr M^1([0,T]\times \mathbb R^N)$, the limit $J$ is
  divergence--free and supported in $[0,T]\times \bar{\Omega}$.
\end{lem}

\noindent
{\bf Proof.}
We follow the arguments of the previous Section.  We identify the
limit of $\rho_\eps$ by coming back to the Poisson equation
$\sqrt\eps\nabla_x\cdot\nabla_x\Psi_\eps=n_{\mathrm e}-\rho_\eps$.
The time compactness then appears as a consequence of the charge
conservation, together with the estimates on the current. 
We obtain the $L^\infty( 0,T;L^1(\mathbb R^N))$ estimate on $J_\eps$
as in \eqref{toto}.  Letting $\eps$ go to 0 in the charge conservation
equation, we obtain $\partial_t n_{\mathrm e}+ \nabla_x\cdot
J=0=\nabla_x\cdot J$. Still reproducing the arguments of the previous
section, based on the conservation of the total charge, we arrive at
the following conclusion:
\[\mathrm{supp}(f)\subset [0,T]\times \bar{\Omega} \times \mathbb R^N,\qquad 
\mathrm{supp}(J)\subset [0,T]\times \bar{\Omega}.\] Furthermore, $J$
belongs to the set $\mathscr D\mathscr M^
{\mathrm{ext}}(\mathbb
R^N)$, 
it admits a normal trace $J\cdot \nu\big|_{\p \Omega}$, which
actually vanishes. \qed

It remains to identify the limit $J$. As in the case of the pure
Vlasov--Poisson equation, the idea consists in introducing a suitable
functional intended to compare $f_\eps$ to the expected limit. Let
$\mathscr N_\eps:\mathbb R^N\rightarrow (0,\infty)$ be a given
function such that
\[
\ds\int \mathscr N_\eps\ud x= \mathfrak{m} =\ds\int n_{\mathrm e}\ud
x=\ds\iint f(0,x,v)\ud v\ud x 
\]
and let us set 
\[
M_{\mathscr V, \theta_\eps}(t,x,v)=\ds\frac{1}{(2\pi\theta_\eps)^{N/2}}\
\exp\left(-\ds\frac{|v-\mathscr V(t,x)|^2}{2\theta_\eps}\right) .\] A
natural candidate to replace the functional $ \mathscr H_{\mathscr
  V,\eps} $ would be the relative entropy of $f_\eps$ with respect to
$ n_{\mathrm e} (x) M_{\mathscr V, \theta_\eps}(t,x,v) $ associated with the
non-negative convex function $z\mapsto z\ln(z) - z+1$, namely
$$
\ds\iint \left(f_\eps\ln\left(\ds\frac{f_\eps}{n_{\mathrm e} M_{\mathscr V, \theta_\eps}}\right) 
- f_\eps + n_{\mathrm e} M_{\mathscr V, \theta_\eps} \right)\ud v\ud x ,
$$
but the first term is clearly meaningless since $ n_{\mathrm e} $ has
compact support.  Therefore, we introduce
\begin{equation}
\label{defNe}
\mathscr N_\eps(x)= \frac{\mathfrak{m}}{ Z_{\eps} } \exp\Big(-\ds\frac{\Phi_{\mathrm e}(x)}{ \eps \theta_\eps}\Big) , 
\quad \quad \quad {\rm where} \quad \ds 
Z_{\eps} = \int \exp\Big(-\ds\frac{\Phi_{\mathrm e}(y)}{ \eps \theta_\eps}\Big) \ud y ,
\end{equation}
and the following modulated functional
\begin{align}
\label{modulo}
\mathscr H^{\mathrm{FP}}_{\mathscr V,\eps} = & \ \theta_\eps \ds\iint
\left(f_\eps\ln\left(\ds\frac{f_\eps}{\mathscr N_\eps M_{\mathscr V,
        \theta_\eps}}\right) - f_\eps + \mathscr N_\eps M_{\mathscr V,
    \theta_\eps}\right)\ud v\ud x + \ds\frac12\ds\int
|\nabla_x\Psi_\eps|^2\ud x .
\end{align}
In fact, $\mathscr H^{\mathrm{FP}}_{\mathscr V,\eps}$ is up to the
term $ \frac12\int |\nabla_x\Psi_\eps|^2\ud x $, nothing but the
relative entropy of $f_\eps$ with respect to $\mathscr N_\eps
M_{\mathscr V, \theta_\eps} $ associated with the non-negative convex
function $ G: (0,+\infty ) \ni z \mapsto z \ln(z) - z+ 1 $.  This
implies in particular that the integrand in the first integral of
\eqref{modulo} is simply $ G( f_\eps) - G(N_\eps M_{\mathscr V,
  \theta_\eps} ) - G'(f_\eps) ( f_\eps -N_\eps M_{\mathscr V,
  \theta_\eps} )$, thus pointwise nonnegative, and vanishes only when
$ f_\eps = N_\eps M_{\mathscr V, \theta_\eps} $.  By definition of $
N_\eps $ and $ M_{\mathscr V, \theta_\eps} $ and using the fact $
\iint f_\eps \ud v\ud x = \mathfrak{m} = \iint \mathscr N_\eps
M_{\mathscr V, \theta_\eps} \ud v\ud x $ in view of our
normalizations, we infer
\begin{align}
\mathscr H^{\mathrm{FP}}_{\mathscr V,\eps} 
= & \ 
\theta_\eps \ds\iint f_\eps\ln(f_\eps) \ud v\ud x 
+ \ds \frac12 \iint \ds|v-\mathscr V|^2\ f_\eps\ud v\ud x 
+\ds\frac1\eps\ds\int \Phi_{\mathrm e}f_\eps\ud v\ud x 
\\ 
\nonumber & + \ds\frac12\ds\int |\nabla_x\Psi_\eps|^2\ud x 
+ \frac{1}{2} N \mathfrak{m} \theta_\eps \ln( 2\pi\theta_\eps) 
- \theta_\eps \mathfrak{m} \ln \left( \frac{\mathfrak{m}}{Z_\eps} \right)
\\ 
\nonumber
= & \ \mathscr H_{\mathscr V,\eps} + \theta_\eps \iint f_\eps\ln(f_\eps) \ud v\ud x 
+ \frac{1}{2} N \mathfrak{m} \theta_\eps \ln( 2\pi\theta_\eps) 
- \theta_\eps \mathfrak{m} \ln \left( \frac{\mathfrak{m}}{Z_\eps} \right).
\end{align}
This second expression of $ \mathscr H^{\mathrm{FP}}_{\mathscr V,\eps}
$ justifies the choice we have made for $ \mathscr N_\eps $. Actually,
for our purpose, the exact normalization $ \iint \mathscr N_\eps
M_{\mathscr V, \theta_\eps} \ud v\ud x = \mathfrak{m} $ is not
necessary, though natural in a modulated entropy argument, only the
fact that $ \ln( \mathscr N_\eps M_{\mathscr V, \theta_\eps} ) \approx
- \Phi_{\mathrm e}(x)/ ( \eps \theta_\eps) - |v-\mathscr V|^2/(2
\theta_\eps) $ is used. This is related to the fact that the
temperature $ \theta_\eps $ is small in the regime we are considering.

Let us now compare $ \mathscr H_{\mathscr V,\eps}$ and $ \mathscr
H^{\mathrm{FP}}_{\mathscr V,\eps} $ more precisely. As a first step,
note that, on the one hand,
$$
\theta_\eps \ln( 2\pi\theta_\eps) \to 0 
$$
when $ \theta_\eps \to 0 $; and on the other hand, that
$$
|\Omega | = \int_{\Omega} \exp\Big(-\ds\frac{\Phi_{\mathrm e}(y)}{
  \eps \theta_\eps}\Big) \ud y \le Z_{\eps} = \int
\exp\Big(-\ds\frac{\Phi_{\mathrm e}(y)}{ \eps \theta_\eps}\Big) \ud y
\le \int \exp\Big(-\lambda \Phi_{\mathrm e}(y) \Big) \ud y < +\infty
$$
if $ \theta_\eps \le 1 $ and $ \eps \le 1/(8\lambda) $, thus, as $ \theta_\eps \to 0 $,
$$
\theta_\eps \mathfrak{m} \ln \left( \frac{\mathfrak{m}}{Z_\eps} \right) \to 0 .
$$
The inequality \eqref{txtx} 
implies
\begin{align}
\label{comparaison}
\mathscr H_{\mathscr V,\eps} = & \ \mathscr H^{\mathrm{FP}}_{\mathscr V,\eps}
- \theta_\eps \iint f_\eps\ln(f_\eps) \ud v\ud x 
- \frac{1}{2} N \mathfrak{m} \theta_\eps \ln( 2\pi\theta_\eps )
+\theta_\eps \mathfrak{m} \ln \left( \frac{\mathfrak{m}}{Z_\eps} \right)
\nonumber \\ \le & \ 
\mathscr H^{\mathrm{FP}}_{\mathscr V,\eps} 
 - \theta_\eps \iint f_\eps\ln(f_\eps) {\bf 1}_{f_\eps \le 1 } \ud v\ud x + o_{\eps \rightarrow 0} (1) 
 \nonumber \\ \le & \ 
 2 \mathscr H^{\mathrm{FP}}_{\mathscr V,\eps} + o_{\eps \rightarrow 0} (1) .
\end{align}

Then, let us compute the time derivative of the modulated entropy 
$ \mathscr H^{\mathrm{FP}}_{\mathscr V,\eps} $. We get
\[\begin{array}{lll}
\ds\frac{\ud}{\ud t}\mathscr H_{\mathscr V,\eps}^{\mathrm{FP}}
&=&
\ds\frac{\ud}{\ud t}\left\{\theta_\eps \ds\iint f_\eps\ln(f_\eps) \ud v\ud x
+ \ds \frac12 \iint \ds|v-\mathscr V|^2\ f_\eps\ud v\ud x
+\ds\frac1\eps\ds\int \Phi_{\mathrm e}f_\eps\ud v\ud x\right.
\\
&&\qquad\left.
+
\ds\frac12\ds\int |\nabla_x\Psi_\eps|^2\ud x
\right\}
\\
&=&
\ds\frac{\ud}{\ud t}\left\{\theta_\eps \ds\iint f_\eps\ln(f_\eps) \ud v\ud x 
+ \ds\frac12\ds\iint \lvert v \rvert^2\ f_\eps\ud v\ud x 
+\ds\frac1\eps\ds\int \Phi_{\mathrm e}f_\eps\ud v\ud x
+
\ds\frac12\ds\int |\nabla_x\Psi_\eps|^2\ud x\right\}
\\
&&
\qquad +\ds\frac{\ud}{\ud t}\left\{-\ds\iint v\cdot \mathscr V f_\eps\ud v\ud x
+ \ds\frac12\ds\iint | \mathscr V|^2 f_\eps\ud v\ud x \right\}.
\end{array}\]
Bearing in mind the computation for proving Proposition \ref{bdd2}, we obtain
\[
\ds\frac{\ud}{\ud t}\mathscr H^{\mathrm{FP}}_{\mathscr V,\eps}
=
-\mathscr D_\eps +
\ds\frac{\ud}{\ud t}\left\{-\ds\int J_\eps\cdot \mathscr V\ud x 
+ \ds\frac12\ds\int\rho_\eps |\mathscr V|^2\ud x 
\right\}.
\]
Reasoning as in the previous section, and by using the moment equations, we are led to 
\[\begin{array}{lll}
\ds\frac{\ud}{\ud t}\mathscr H^{\mathrm{FP}}_{\mathscr V,\eps}
&=&-\mathscr D_\eps +\ds\int \mathscr V\cdot J_\eps\ud x
\\
&&
+
\ds\int (\rho_\eps \mathscr V-J_\eps)\cdot\big( \partial_t \mathscr V
+(\mathscr V\cdot \nabla_x)\mathscr V\big) \ud x 
\\&&- \ds\int D\mathscr V:(\mathbb P_{\mathscr V, \eps}-\nabla_x\Psi_\eps\otimes\nabla_x\Psi_\eps)\ud x
+\ds\frac1\eps\ds\int \rho_\eps \mathscr V\cdot \nabla_x\Phi_{\mathrm e}\ud x ,
\end{array}\]
using once again that $\nabla_x\cdot (n_{\rm e}\mathscr V)=0$ and the no-flux condition \eqref{CLEuler}. 
Let us set 
\[
\mathscr D_{\mathscr V,\eps}=\ds\iint \big|(v-\mathscr V)\sqrt{f_\eps}
+2\theta_\eps \nabla_v\sqrt{f_\eps}\big|^2\ud v\ud x\geq 0. 
\]
We rewrite
\[
\mathscr D_\eps= \mathscr D_{\mathscr V,\eps}-\ds\int \rho_\eps|\mathscr V|^2\ud x
+2\ds\int \mathscr V\cdot J_\eps\ud x.\]
Accordingly, we can reorganize terms as follows
\[\begin{array}{lll}
\ds\frac{\ud}{\ud t}\mathscr H^{\mathrm{FP}}_{\mathscr V,\eps}&=&-\mathscr D_{\mathscr V,\eps}+
\ds\int (\rho_\eps \mathscr V-J_\eps)\cdot\big( \partial_t \mathscr V 
+(\mathscr V\cdot \nabla_x)\mathscr V+\mathscr V\big) \ud x 
\\
&&
- \ds\int D\mathscr V:(\mathbb P_{\mathscr V, \eps}-\nabla_x\Psi_\eps\otimes\nabla_x\Psi_\eps)\ud x
+\ds\frac1\eps\ds\int \rho_\eps \mathscr V\cdot \nabla_x\Phi_{\mathrm e}\ud x .
\end{array}\]

\noindent
We can summarize the previous manipulations within the following inequality 
\begin{equation}\label{tata}
\ds\frac{\ud}{\ud t}\mathscr H^{\mathrm{FP}}_{\mathscr V,\eps}  +\mathscr D_{\mathscr V,\eps} \leq 
\ds\frac1\eps \ds\int \rho_\eps \mathscr V\cdot\nabla_x\Phi_{\mathrm e}\ud x+r_\eps
+
\ds\int D\mathscr V:(\mathbb P_{\mathscr V,\eps} - \nabla_x\Psi_\eps\otimes \nabla_x\Psi_\eps)\ud x,
\end{equation}
where, for any $0 <t \le T $,
\[\ds\int_0^t r_\eps\ud \tau=
\ds\int_0^t \ds\int (\rho_\eps \mathscr V-J_\eps)(\partial_\tau \mathscr V 
-\mathscr V\cdot \nabla_x\mathscr V+\mathscr V)\ud x\ud \tau\]
tends to 0 as $\eps\rightarrow 0$.
We wish to strengthen 
this result as follows.

\begin{lem}\label{l:estFP}
We make assumptions \eqref{integral}  and \eqref{queryVFP} and suppose that 
$\theta_\eps\rightarrow 0$ as $\eps\rightarrow 0$. We have
\[
\ds\frac{\ud}{\ud t}\mathscr H^{\mathrm{FP}}_{\mathscr V,\eps}  +\mathscr D_{\mathscr V,\eps} \lesssim 
 \mathscr H^{\mathrm{FP}}_{\mathscr V,\eps} + r_\eps
\]
where, for any $ 0 < t \le T $, $\lim_{\eps\rightarrow 0}\int_0^t r_\eps \ud \tau =0$.
\end{lem}

\noindent
{\bf Proof.}
We can also reproduce the arguments in the previous section used to estimate 
\[
\ds\frac1\eps \ds\int \rho_\eps \mathscr V\cdot\nabla_x\Phi_{\mathrm e}\ud x\lesssim 
\ds\frac1\eps \ds\int \rho_\eps\Phi_{\mathrm e}\ud x. 
\]
For the last term in \eqref{tata}, we have
\[\ds\int D\mathscr V:(\mathbb P_{\mathscr V,\eps} - \nabla_x\Psi_\eps\otimes \nabla_x\Psi_\eps)\ud x
\leq \|D\mathscr V\|_\infty \left(\ds\iint |v-\mathscr V|^2 f_\eps \ud v\ud x 
+ \ds\int | \nabla_x\Psi_\eps|^2\ud x\right),  
\]
so that, using \eqref{comparaison},
\begin{align*}
\ds\frac{\ud}{\ud t}\mathscr H_{\mathscr V,\eps}^{\mathrm{FP}}  +\mathscr D_{\mathscr V,\eps} 
\lesssim & \ 
\ds\frac1\eps \ds\iint f_\eps\Phi_{\mathrm e}
\ud v\ud x+ \frac12\ds\iint |v-\mathscr V|^2 f_\eps \ud v\ud x 
+ \frac12 \ds\int | \nabla_x\Psi_\eps|^2\ud x+
r_\eps 
\\ 
\lesssim & \ \mathscr H_{\mathscr V,\eps}^{\mathrm{FP}} + r_\eps + o_{\eps \rightarrow 0 }(1).
\end{align*}
It allows us to conclude by coming back to \eqref{tata}.
\qed

Let us now state our main result concerning the Vlasov-Poisson-Fokker-Planck system. 
We recall that we may work either with a quadratic potential $ \Phi_{\rm ext} $ (and then 
the domain $\Omega$) is an ellipsoid), or with a general potential where h1), h2), H1) and 
H2) are satisfied.

\begin{theorem} 
If $N \ge 3 $, we make assumption \eqref{integral}, that is we assume that there exists 
some (large) $ \lambda > 1 $ such that
$$
 \int \exp( - \lambda \Phi_\ext ) \ud x < \infty .
$$
Denote by $V$ the solution, on $ [0,T ] $, to the Lake Equation with friction \eqref{LE2} with the 
no-flux condition \eqref{CLEuler} given by Theorem \ref{Lakeexiste} and consider 
a smooth extension $\mathscr V$ to $V$. 
Let $f_\eps^{\mathrm{init}}:\mathbb R^N\times \mathbb R^N\rightarrow [0,\infty)$ be a sequence 
of integrable functions satisfying
$$
 \iint f_\eps^{\mathrm{init}}\ud v\ud x = \mathfrak{m} 
 \quad \quad \quad {\it and} \quad \quad \quad
 \mathscr H_{\mathscr V,\eps}^{\mathrm{FP}} (f_\eps^{\mathrm{init}} ) \to 0 ,
$$
where $ \mathscr H_{\mathscr V,\eps}^{\mathrm{FP}} $ is defined in \eqref{modulo}. 
Consider then the associated solutions $f_\eps$ of the Vlasov--Poisson--Fokker--Planck equation 
\eqref{FPV1}--\eqref{VP2}. 
Then, we have, as $ \eps \rightarrow 0 $ and $ \theta_\eps \to 0 $,
\begin{itemize}
\item[i)] $\rho_\eps$ converges to $n_{\mathrm e}$ in 
$C^0(0,T;\mathscr M^1(\mathbb R^N)-\text{weak}-\star)$;
\item[ii)] $\mathscr H_{\mathscr V,\eps}^{\mathrm{FP}} \rightarrow 0$ uniformly on $[0,T] $; 
\item[iii)] $J_\eps$ converges to $J$ in $\mathscr M^1([0,T]\times \mathbb R^N)$, 
where $J\big|_{[0,T]\times \Omega }=V$, $\nabla_x\cdot J=0$ and $J\cdot \nu(x)\big|_{\p \Omega}=0$.
\end{itemize}\end{theorem}

\begin{remark}
 We have seen (see Remark \ref{crapaud}) that the integrability assumption 
 $ \int \exp( - \lambda \Phi_\ext ) \ud x < \infty $ is automatically satisfied if 
 $N=1$, $2$ by h2) or for quadratic potentials. When $ N \ge 3 $, it is also 
 true if $ \Phi_{\rm ext} $ is convex and tends to $+\infty $ at infinity.
\end{remark}
\begin{remark}
  One may construct an admissible family of initial conditions
  following the lines of Remark 1.3. In particular, taking $G$ a
  normalized Gaussian, it is enough to choose $\theta_\eps$ and
  $\sigma_\eps$ such that $\theta_\eps \int f_\eps \ln f_\eps \to 0$,
  which imposes $\theta_\eps \ln \sigma_\eps\to 0$.
\end{remark}

\noindent
{\bf Proof.} It is clear that if $ \mathscr H_{\mathscr V,\eps}^{\mathrm{FP}} (f_\eps^{\mathrm{init}} ) \to 0 $, 
then \eqref{queryVFP} is satisfied. Item i) has already been discussed. 
Applying the Gr\"onwall lemma, we deduce readily that ii) holds from Lemma \ref {l:estFP}.
Coming back to \eqref{comparaison}, we infer that $\iint|v-\mathscr V|^2 f_\eps \ud v\ud x$ tends to 0.
Then, we appeal to Lemma \ref{l:Vx} to conclude that $J$ belongs to $L^\infty(0,T;L^2(\mathbb R^N))$ and 
that $J=n_{\mathrm e} V$. We can also justify some time--compactness as in the pure Vlasov--Poisson case.
%
%
%

\appendix

\section{Smooth solutions of the Lake Equations}
\label{sec:eq_lake}

\begin{theorem}\label{Lakeexiste}
 Let $ \Omega $ be a smooth ($ \p \Omega $ of class
$ C^{s+1} $ is enough) bounded open set in $ \R^N $,
$ \gamma $ be a real constant, $ s \in \N $ such that $s>1+N/2 $
and $ n_{\rm e} : \Omega \to \R $ in 
$ H^{s + 1 } $ such that $  \inf_\Omega n_{\rm e} > 0 $.
Let $V^{\mathrm{init}}: \Omega \rightarrow \mathbb R^N$ be a divergence free vector field
in $H^s$ satisfying the no flux condition $ V^{\mathrm{init}} \cdot \nu = 0 $
on $\p \Omega $. There exists $ T >0 $ and a unique solution
$ V \in L^\infty (0,T ;H^s(B(0,R)))$ of 
\begin{equation*}
\label{lake_f}
\left\{ \begin{array}{l}
\partial_t V + V \cdot \nabla_x V + \nabla_x p = - \gamma V , \\
\nabla_x\cdot ( n_{\rm e} V ) = 0 ,
\end{array} \right.
\end{equation*}
with the no flux condition \eqref{CLEuler}. Moreover, we have
\[
\ds
\sup_{0\leq t\leq T}\Big(
\|V(t)\|_{H^s} + \|\partial_t V(t)\|_{H^{s-1}} + \|\nabla_xp(t)\|_{H^s}
+ \|\partial_t\nabla_x p(t)\|_{H^{s-1}}\Big)\leq C(T ) \]
for some positive constant $C(T )$ depending on $ \gamma $, $T$, $ n_{\rm e} $ and the initial datum.
\end{theorem}

\noindent {\it Proof.} The scheme of proof is exactly the same as in \cite{RT}.
We shall denote $ \hat V \ddef n_{\rm e} V $, which is divergence free.
Applying $ \nabla \cdot( n_{\rm e} \, . ) $ to the equation, we see that
the pressure $p$ satisfies, for any $t$, the elliptic equation
\begin{align}
\label{Elliptik}
 - \nabla_x \cdot( n_{\rm e} \nabla_x p )
 = & \,
 \nabla_x \cdot \left( n_{\rm e} V \cdot \nabla_x V \right)
 =
 \nabla_x \cdot \left( \hat V \cdot \nabla_x \left( \frac{1}{n_{\rm e} }\hat V \right) \right)
 \nonumber \\
 = & \, 
 \hat V \cdot \nabla_x \left( \nabla_x \cdot \left( \frac{1}{n_{\rm e} }\hat V \right) \right) 
 + 
 \sum_{1 \le j,k \le N } \p_j \hat V_k \p_k \left( \frac{\hat V_j }{n_{\rm e}}  \right)
 \nonumber \\
 = & \, 
  \hat V \cdot \nabla_x \left( \hat V \cdot \nabla_x \left( \frac{1}{n_{\rm e} } \right) \right) 
  + 
 \sum_{1 \le j,k \le N } \p_j \hat V_k \p_k \left( \frac{\hat V_j }{n_{\rm e}}  \right) ,
\end{align}
where we have used that $\hat V $ is divergence free as well as the identity
$ \nabla_x \cdot ( \mathcal V \cdot \nabla_x \mathcal U ) -
\mathcal V \cdot \nabla_x ( \nabla_x \cdot \mathcal U ) =
\sum_{1 \le j,k \le N} \p_j \mathcal V_k \p_k \mathcal U_j $. 
We may further impose a suitable Neumann boundary condition for $p$ on $ \p \Omega $. 
We recall that for $ \sigma > N/2 $, $ H^\sigma $ is an algebra. 
Notice that if $ V \in H^s $, with $ s >  1 + N/2 $, then the right-hand
side of \eqref{Elliptik} is in $ H^{s-1} $ and
$$
 \left \lVert 
 \hat V \cdot \nabla_x \left( \hat V \cdot \nabla_x \left( \frac{1}{n_{\rm e} } \right) \right) 
  + 
 \sum_{1 \le j,k \le N } \p_j \hat V_k \p_k \left( \frac{\hat V_j }{n_{\rm e}}  \right) 
 \right \rVert_{H^{s-1}} \le C \lVert V \rVert_{H^s}^2 ,
$$
where $C$ depends on $ \inf_\Omega n_{\rm e} $ (which is assumed positive)
and the $ H^{s+1} $ norm of $ n_{\rm e} $.

Since $ n_{\rm e} $ is in $H^{s+1} $ and bounded away from zero and since the
boundary is assumed of class $ C^{s+1} $, it follows from classical elliptic
estimates that \eqref{Elliptik} endowed with the Neumann condition 
on $ \p \Omega $ has a unique solution $ p \in H^{s+1} ( \Omega )\Big|_\mathbb R $, enjoying the
estimate
\begin{equation}
\label{demi}
 \lVert p \rVert_{H^{s+1} } \le C \lVert V \rVert_{H^s}^2 ,
\end{equation}
where $C$ depends on $  \inf_\Omega n_{\rm e} $ and
the $ H^{s+1} $ norm of $ n_{\rm e} $. Assume now that
$ V $ is a smooth solution of \eqref{lake_f} and let us perform an $ H^s $ estimate.
For any $ \alpha \in \big(\N\cup\{0\}\big)^d $ with $ | \alpha | \le s $, we have
\begin{align*}
 \frac{\ud }{\ud t} \int_\Omega | \p^\alpha V |^2 \, \ud x
 = & \, - 2 \int_\Omega \p^\alpha V \cdot \p^\alpha \left( ( V \cdot \nabla_x ) V \right) \, \ud x
 - 2 \int_\Omega \p^\alpha V \cdot \p^\alpha \nabla_x p \, \ud x 
 - 2 \gamma \int_\Omega | \p^\alpha V |^2 \, \ud x
\end{align*}
Using classical commutator estimates, the Sobolev imbedding $ H^s \subset W^{1,\infty} $
and the $ H^{s+1} $ estimate \eqref{demi} on $p$, we then deduce
\begin{align*}
 \frac{\ud }{\ud t} \int_\Omega | \p^\alpha V |^2 \, \ud x
 \le &\, - 2 \int_\Omega \p^\alpha V \cdot \left( ( V \cdot\nabla_x ) \p^\alpha V \right) \, \ud x 
 + C ( \lVert V \rVert_{H^s} + \lVert V \rVert_{H^s}^2 + \lVert V \rVert_{H^s}^3 ) .
\end{align*}
We use integration by parts for the first integral (recall that $ V \cdot \nu = 0 $ on the
boundary), which then becomes
$ \int_\Omega ( \nabla_x \cdot V ) | \p^\alpha V |^2 \, \ud x \le C \lVert V \rVert_{H^s}^3 $. 
This yields
$$
 \frac{\ud }{\ud t} \int_\Omega | \p^\alpha V |^2 \, \ud x \le
 C ( \lVert V \rVert_{H^s} + \lVert V \rVert_{H^s}^2 + \lVert V \rVert_{H^s}^3 ) ,
$$
and it follows that, for some $ T_0 > 0 $ depending only on $\gamma $, $n_{\rm e} $
and $ V^{\mathrm{init}} $, we have $ \lVert V \rVert_{L^\infty ( 0,T_0 ; H^s )} \le 2
\lVert V^{\mathrm{init}} \rVert_{H^s } $. The conclusion of the theorem follows from a 
suitable viscous approximation where a careful treatment of boundary terms is needed, see 
\cite{Temam2}. \qed

\section{Construction of an extended divergence--free velocity}

\begin{lem}
\label{extension}
Let  $V \in L^\infty( 0,T ; H^s( B(0,R) , \mathbb R^N) ) $ be a 
divergence free vector field in $H^s$, with $s>1+N/2$, satisfying the no 
flux condition $ V \cdot \nu = 0 $ on $\p B (0,R) $.
There exists a solenoidal extension $\mathscr V$ of the vector field $V$ 
defined on the whole space and compactly supported. 
Namely, $ \mathscr V\in L^\infty(0,T;H^s(\mathbb R^N , \mathbb R^N)) $ and it satisfies:\\
i) \eqref{CondExtension},\\
ii) $\nabla_x \cdot \mathscr V =0$ in $\mathbb R^N$. 
\end{lem}

\noindent \emph{Proof.} Let us assume $ N =2 $ or $ N =3 $. Since $
\nabla_x\cdot V=0 $ in the ball $ B(0,R) $, which is convex, there
exists $ h \in L^\infty (0,T ; H^{s+1}(B(0,R), \mathbb R^{N}) ) $ such
that $ V = \nabla \times h $. 
Then, by standard extension results (see, e. g., \cite[Chapter I:
Theorem 2.1 p. 17 \& Theorem 8.1 p. 42]{LM}), there exists an
extension $ \tilde{h} \in L^\infty ( 0,T ; H^{s+1}(\mathbb R^N,
\mathbb R^{N}) )$ to $ h $.  Considering a cut-off function $ \chi \in
C^\infty_c(\mathbb R^N) $ such that $ \chi (x)= 1 $ for $x\in B(0,
3R/2 ) $ and denoting $ \mathscr V = \nabla \times ( \chi \tilde{h} )
$, we see that $ \mathscr V $ enjoys the desired properties.  For $N
\ge 4 $, the construction is similar but involves differential
forms. The arguments generalize to the case where $\Omega$ is an
ellipsoid. \qed

\begin{remark}
  In the case of a general potential, $n_{\rm e}$ is not uniform, and it is 
possible to construct an extension 
$ \mathscr V\in L^\infty(0,T;H^s(\mathbb R^N , \mathbb R^N)) $ which satisfies:\\
  i) \eqref{CondExtension},\\
  ii) $\nabla_x \cdot \left( n_{\rm e} \mathscr V\right) =0$ in $\mathbb R^N$.\\

  However, it requires further topological hypotheses on $\Omega$.
  Assuming H1, and assuming also that $\mathcal{K}$ is connected, and
  $\p \Omega$ has a finite number of connected components, we may
  apply \cite[Corollary 3.2]{Kato}: $ n_{\rm e} V $ is divergence free
  in the smooth domain $ \Omega $, hence we can construct a divergence
  free extension $ \mathscr{J} : [0,T] \times \R^N \to \R^N $ to $
  n_{\rm e} V $. Using a cut-off function, we may take $ \mathscr{J} $
  compactly supported in an arbitrary neighborhood of $ \mathcal{K} $,
  the latter can be chosen so that $ \Delta \Phi_{\rm ext} $ remains $
  > 0 $.  Finally, we set $ \mathscr V = \mathscr{J} / (\Delta
  \Phi_{\rm ext} ) $, which is well-defined even when $ \Delta
  \Phi_{\rm ext} $ vanishes.

\end{remark}

\section*{Acknowledgements}
Th. Goudon thanks Courant Institute in NYC, where a part of this work
has been done, for its hospitality. We are gratefully indebted to
S. Serfaty for many kind advices on the generalized Gauss problem and
the obstacle problem, and to A. Olivetti for providing figure 1.

\bibliographystyle{plain}
\bibliography{./VPE}

\end{document}